\documentclass[12pt]{article}
\usepackage{amssymb,amsmath}

\input epsf.tex

\title{Affine structures and non-archimedean analytic spaces}

\author {Maxim Kontsevich and Yan Soibelman}

\begin{document}
\maketitle

\newtheorem{defn}{Definition}
\newtheorem{thm}{Theorem}
\newtheorem{lmm}{Lemma}
\newtheorem{rmk}{Remark}
\newtheorem{prp}{Proposition}
\newtheorem{conj}{Conjecture}
\newtheorem{exa}{Example}
\newtheorem{cor}{Corollary}
\newtheorem{que}{Question}
\newtheorem{ack}{Acknowledgements}
\newcommand{\C}{{\bf C}}
\newcommand{\K}{{\bf k}}
\newcommand{\R}{{\bf R}}
\newcommand{\N}{{\bf N}}
\newcommand{\Z}{{\bf Z}}
\newcommand{\Q}{{\bf Q}}
\newcommand{\G}{\Gamma}
\newcommand{\A}{A_{\infty}}
\newcommand{\g}{{\bf g}}

\newcommand{\epi}{\twoheadrightarrow}
\newcommand{\mono}{\hookrightarrow}
\newcommand\ra{\rightarrow}
\newcommand\uhom{{\underline{Hom}}}
\renewcommand\O{{\cal O}}
\newcommand\nca{nc{\bf A}^{0|1}}
\newcommand{\epp}{\varepsilon}

\begin{abstract}
In this paper we propose a way to construct an analytic space
over a non-archimedean field,
starting with a real manifold with an affine structure 
which has integral monodromy. Our construction is motivated 
by the junction of Homological Mirror conjecture 
and geometric Strominger-Yau-Zaslow conjecture.
In particular, 
we glue 
from ``flat pieces" an analytic K3 surface. As a byproduct of our approach we obtain
an action of an arithmetic subgroup of the group $SO(1,18)$ by
piecewise-linear transformations on the $2$-dimensional
sphere $S^2$ equipped with naturally defined 
singular affine structure.
\end{abstract}

\section{Introduction}

\subsection{}

An integral affine structure on a
manifold of dimension $n$ is
given
by a torsion-free flat connection with the monodromy
reduced to $GL(n,{\Z})$.
There are two basic situations in which integral affine
structures occur naturally. One is the case of classical integrable systems 
described briefly in Section 3.
Most interesting for us is a class
of examples arising from analytic manifolds over non-archimedean fields which is
discussed in Section 4. It is motivated by
the approach to Mirror Symmetry suggested in [KoSo]. We recall it in
Section 5. 
From our point of view
manifolds with integral affine structure
appear in Mirror Symmetry in two ways.
One considers the Gromov-Hausdorff collapse
of degenerating families of Calabi-Yau manifolds.
 The limiting space can be interpreted either as a contraction (see Section 4.1) of an analytic
 manifold
 over a non-archimedean field of Laurent series $\C((t))$, or as a base of a fibration
 of a Calabi-Yau manifold by Lagrangian tori (with respect to the symplectic K\"ahler 2-form).
 On a dense open subset of the limiting space one gets two integral affine structures
 associated with two interpretations, the non-archimedean one and the symplectic one.
Mirror dual family of degenerating Calabi-Yau manifolds should have metrically the same
 Gromov-Hausdorff limit, with the roles of two integral affine structures interchanged.

Very interesting question arises: how to reconstruct these
families of Calabi-Yau manifolds from the corresponding manifolds
with integral affine structures?
This question was one of the main motivations for present work.

\subsection{}
Our approach to the reconstruction
of analytic Calabi-Yau manifolds from real manifolds with
integral affine structure
can be illustrated in the following toy-model
example. Let $S^1={\R}/{\Z}$ be a circle equipped with the
induced from ${\R}$ affine structure. We equip $S^1$ with
the canonical sheaf ${\cal O}^{can}_{S^1}$
of Noetherian ${\C}((q))$-algebras. By definition, for an open
interval $U\subset S^1$  algebra 
${\cal O}^{can}_{S^1}(U)$ consists of formal series
$f=\sum_{m,n\in {\Z}}a_{m,n} q^m z^n,\,\,\,a_{m,n}\in \C$ such that $\inf_{a_{m,n}\ne 0} (m+nx)> -\infty$. Here $x\in \R$
 is any point in a connected component of the pre-image of $U$ in $\R$,
 the choice of a different component $x\rightarrow x+k, \,\,k\in \Z$ corresponds to the substitution
 $z\mapsto q^k z$. The corresponding analytic space is the  Tate
elliptic curve $(E,{\cal O}_E)$, and there is a continuous
map $\pi:E\to S^1$ such that $\pi_{\ast}({\cal O}_E)=
{\cal O}^{can}_{S^1}$.

In the case of K3 surfaces one starts with $S^2$.
The corresponding integral
affine structure is well-defined on the  set
$S^2\setminus \{x_1,...,x_{24}\}\subset S^2$,
where $x_1,...,x_{24}$ are distinct points.
Similarly to the above toy-model example one can construct
the canonical sheaf 
${\cal O}_{S^2\setminus \{x_1,...,x_{24}\}}^{can}$ of 
algebras,
an open $2$-dimensional smooth analytic surface $X^{\prime}$
with the trivial canonical bundle (Calabi-Yau manifold),
and a continuous projection 
$\pi^{\prime}: X^{\prime}\to S^2\setminus \{x_1,...,x_{24}\}$
such that $\pi^{\prime}_{\ast}({\cal O}_{X^{\prime}})=
{\cal O}_{S^2\setminus \{x_1,...,x_{24}\}}^{can}$.
The problem is to find a sheaf  ${\cal O}_{S^2}$
whose restriction to $S^2\setminus \{x_1,...,x_{24}\}$ is locally isomorphic to
 ${\cal O}_{S^2\setminus \{x_1,...,x_{24}\}}^{can}$, an analytic compact K3 surface $X$, and a continuous
projection $\pi:X\to S^2$ such that $\pi_{\ast}({\cal O}_X)=
{\cal O}_{S^2}$. 
We call this problem (in general case) the
Lifting Problem and discuss it in Section 7.
 Unfortunately we do not know the conditions
one should impose on singularities of the affine
structure, so that the Lifting Problem would have a solution.
We consider a special case of
K3 surfaces in Sections 8-11. Here the solution is
non-trivial and depends on data which are not visible
in the statement of the problem. They are motivated
by Mirror Symmetry and consist, roughly speaking, of
an infinite collection of trees embedded into 
$S^2\setminus \{x_1,...,x_{24}\}$ with the tail vertices
belonging to the set $\{x_1,...,x_{24}\}$. The sheaf
${\cal O}_{S^2\setminus \{x_1,...,x_{24}\}}^{can}$ has
to be modified by means of automorphisms assigned to every
edge of a tree and then glued together with certain model
sheaf near each singular point $x_i$.

Informally speaking, we break $S^2\setminus \{x_1,...,x_{24}\}$ endowed with the 
sheaf ${\cal O}_{S^2\setminus \{x_1,...,x_{24}\}}^{can}$
into infinitely many infinitely small pieces and then glue them back together in a slightly deformed way.
The idea of such a construction was proposed several years ago independently by K.~Fukaya and the first author.
 The realization of this idea 
was hindered by a poor understanding of singularities of the Gromov-Hausdorff collapse and by the
 lack of knowledge of certain open Gromov-Witten invariants (``instanton corrections''). 
The last problem is circumvented here (and in fact solved)
with the use
 of some pro-nilpotent Lie group (see Section 10).

\subsection{}
The relationship between K3 surfaces and singular affine
structures on $S^2$ is of very general origin.
Starting with a projective analytic 
Calabi-Yau manifold $X$ over
a complete non-archimedean local field $K$  one can
canonically construct a PL manifold $Sk(X)$ called the
skeleton of $X$. If $X$ is a generic K3 surface then $Sk(X)$
is $S^2$. We discuss skeleta in
Section 6.6. The group of birational automorphisms
of $X$ acts on $Sk(X)$ by integral PL transformations.
For $X=K3$ we obtain an action of an arithmetic
subgroup of $SO(1,18)$ on $S^2$.
Further examples should come from Calabi-Yau manifolds with
large groups of birational automorphisms.

\subsection{}
We have already discussed the content of the paper. Let us
summarize it. The paper is naturally divided into three parts.
Part 1 is devoted to generalities on integral affine structures and
examples, including Mirror Symmetry. Motivated by string theory we use term
 A-model (resp. B-model) for examples arising in symplectic (resp. analytic)
 geometry.

In Part 2 we discuss the concept of singular
integral affine structure, including an affine version
of Gauss-Bonnet theorem. The latter implies that
if all singularities of an integral affine structure on $S^2$
are standard (so-called focus-focus singularities) then 
there are exactly $24$ singular points. Part 2 also
contains a statement of the Lifting Problem and discussion
of flat coordinates on the moduli space of complex Calabi-Yau
manifolds. We expect that under mild conditions on the
singular integral affine structure there exists a solution
of the Lifting Problem, which is unique as long as we fix
periods (see Sections 7.3 and 7.4 for more details).

Most technical Part 3 contains a solution of the Lifting Problem
for K3 surfaces. We construct the corresponding analytic K3 surface
as a ringed space. The sheaf of analytic functions is defined differently
near a singular point and far from the singular set.
It turns out that the ``naive" candidate for the sheaf on the complement
of the singular set has to be modified before we can glue it with
the model sheaf near each singular point. This modification
procedure involves a new set of data (we call them {\it lines}).
We also discuss the group of automorphisms of the canonical
sheaf which preserve the symplectic form. We use this group
in order to modify the ``naive" sheaf along each line.

The paper has two Appendices. First one contains some background
on analytic spaces, while the second one is devoted to Torelli
theorem.

{\it Acknowledgements}. We are grateful to Ilya Zharkov and Mark Gross for useful discussions. 
 Second author thanks Clay Mathematics
Institute for supporting him as a Fellow and IHES
for excellent research and living conditions.

\part{}
\section{$\Z$-affine structures}

\subsection{Definitions}

Let us recall that
an affine structure on manifold $Y$ (smooth, of dimension $n$)
 is given by a torsion-free flat
connection $\nabla$ on the tangent bundle $TY$.

 We will give below three equivalent definitions
of the notion of an integral affine structure.

\begin{defn}

An integral affine structure on $Y$ 
($\Z$-affine structure for short) is an affine structure 
$\nabla$ together
with a $\nabla$-covariant lattice of maximal rank 
$T^\Z=(TY)^{\Z}\subset TY$.

\end{defn}

It is easy to see that if $Y$ carries a $\Z$-affine structure
then for any point $y\in Y$ there exist small neighborhood
$U$, local coordinate system $(x_1,...,x_n)$ in $U$
such that $\nabla=d$
in coordinates $(x_1,...,x_n)$, and the lattice
$(T_xY)^{\Z}, x\in U$ is a free abelian group generated
by the tangent vectors $\partial/\partial x_i\in T_xY, 1\le i\le n$.
Let us call {\it $\Z$-affine} such a coordinate system in $U$
(sometimes we will call such $U$ a $\Z$-affine chart). 
For a covering of
$Y$ by $\Z$-affine charts the transition functions belong (locally) to 
$GL(n,\Z)\ltimes \R^n$. Explicitly, a change of
coordinates is given by the formula

$$ x_i^{\prime}=\sum_{1\le j\le n}a_{ij}x_j+b_i\,\,,$$
where $(a_{ij})\in GL(n,{\Z}), (b_i)\in {\R}^n$.

Hence, Definition 1 is equivalent
to the following

\begin{defn} A $\Z$-affine structure on $Y$ is given by a maximal
atlas of charts such that the transition functions belong locally 
to $GL(n,\Z)\ltimes \R^n$.

\end{defn}

In the above definition $Y$ is just a topological manifold, $C^{\infty}$-structure on it can be reconstructed
 canonically from $\Z$-affine structure.

We can restate the notion of $\Z$-affine
structure in the language
of sheaves of affine functions. 

We say that a real-valued  function $f$
on 
${\R}^n$ is $\Z$-affine if it has the form

$$f(x_1,\dots,x_n)=a_1 x_1+\dots +a_n x_n +b\,\,,$$
where $a_1,\dots,a_n\in \Z$ and $b\in \R$. 
We will denote by $Aff_{\Z,{\R}^n}$  the sheaf of functions
 on ${\R}^n$ which are locally $\Z$-affine.

\begin{defn} A $\Z$-affine structure (of dimension $n$)
on a Hausdorff topological space $Y$ is a subsheaf 
$Aff_{{\Z},Y}$ of the
sheaf of continuous functions on $Y$,
such that the pair $(Y,Aff_{{\Z},Y})$ is locally isomorphic
to  $({\R}^n,Aff_{\Z,{\R}^n})$.
\end{defn} 

Equivalence of the last two definitions follows from the
observation that a homeomorphism between 
two open domains in ${\R}^n$
preserving the sheaf $Aff_{\Z,{\R}^n}$ is given by the same 
formula $x^{\prime}=A(x)+b, A\in GL(n,{\Z}), b\in {\R}^n$ as
the change of coordinates between two $\Z$-affine coordinate
systems.

\subsection{Monodromy representation and its invariant}

With a given affine structure on $Y$ we can associate a flat
affine connection $\nabla^{aff}$ (see [KN]). The corresponding parallel
transport acts on tangent spaces by affine transformations.
For a $\Z$-affine structure the monodromy of $\nabla^{aff}$
belongs to $GL(n,\Z)\ltimes \R^n$, i.e. $\forall y \in Y$ we have a monodromy
representation
$$\rho:\pi_1(Y,y)\to GL(n,\Z)\ltimes {\R}^n\,\,.$$
 Alternatively,
we can define the monodromy representation by covering
a loop in $Y$ by $\Z$-affine coordinate charts and composing
the corresponding transition functions.  

Notice that a $\Z$-affine structure on $Y$ gives rise to a class
$$[\rho]\in H^1(Y,T^{\Z}\otimes {\R})=H^1(Y,T_Y^\nabla)\,\,\,,$$ 
where  $T_Y^\nabla\subset T_Y$ is the subsheaf of $\nabla$-flat sections\footnote{
Here we slightly abuse notations because $Y$ is not necessarily connected.}.  
De Rham representative of  class $[\rho]$
 is given by a 
differential $1$-form $\theta\in \Omega^1(Y,T_Y)$ such that
$\theta(v)=v$ for any tangent vector $v$. In affine coordinates
one has $\theta=\sum_i\partial/\partial x_i\otimes dx_i$.
Clearly $\nabla(\theta)=0$.

We will need later an explicit formula for the $\R$-valued pairing of $[\rho]$ with
 a closed singular 1-chain with coefficients in the
local system $(T^\ast)^\Z=(T^* Y)^\Z$, the dual covariant
lattice in $T^\ast Y$.
 With any singular $1$-chain $c$ with values in $(T^\ast)^\Z$ we
associate a real number $j(c)$ in the following way.
Suppose that $c$ is given by a continuous map
$\gamma:[0,1]\to Y$ and a section $\alpha\in \Gamma([0,1],\gamma^{\ast}(T^{\ast})^\Z)$.
Parallel transport via the connection $\nabla^{aff}$ gives rise to a map
$\overline{\gamma}: [0,1]\to T_{\gamma(0)}Y,\,\,\, \overline{\gamma}(0)=0$.
Let $\alpha_0=\alpha(0)\in (T_{\gamma(0)}^\ast)^{\Z}\subset T^*_{\gamma(0)}Y$. We define 
$j(c)=\langle \alpha_0,\overline{\gamma}(1)\rangle$.
We extend $j(c)$ to an arbitrary singular $1$-chain $c$
by additivity.
Then the class $[\rho]$ can be calculated as
$\langle [\rho], [c]\rangle =j(c)$ for any closed
$1$-chain $c\in C_1(Y, (T^{\ast})^{\Z})$.

\section{A-model construction}

\subsection{Integrable systems}

Let $(X, \omega)$ be a smooth symplectic
manifold of dimension $2n$, $B_0$ a smooth manifold
 of dimension $n$,
$\pi:X\to B_0$ a smooth map
with compact fibers,
such that $\{\pi^{\ast}(f),\pi^{\ast}(g)\}=0$ for any
$f,g\in C^{\infty}(B_0)$. Here $\{\cdot,\cdot\}$ denotes
the Poisson bracket on $X$.
We assume that $\pi$ is a submersion on  
an open dense subset $X^{\prime}\subset X$. Such a triple
$(X,\pi, B_0)$ is called an \emph{integrable system}. In applications
it is typically given by a collection of smooth functions
$(H_1,...,H_n)$ on $X$ (these functions are called Hamiltonians)
such that $\{H_i,H_j\}=0, 1\le i,j\le n$. Usually 
 first Hamiltonian $H=H_1$ is identified with the energy of mechanical
system.

Let us consider the case when $\pi$ is proper. It is a natural 
restriction, because in applications the energy $H_1$ is already
a proper map $H_1:X\to {\R}$.

Let $x\in B_0$ be a point such that the restriction
of $\pi$ to $\pi^{-1}(x)$ is a submersion. We call such points
{\it $\pi$-smooth}. According to Sard theorem $\pi$-smooth
points form an open dense subset of $B_0$.
The fiber $\pi^{-1}(x)$ is a compact 
Lagrangian submanifold of $X$.
The Liouville integrability theorem (see [Ar]) says that
$\pi^{-1}(x)$ is a disjoint union of finitely many tori
$T^n_{\alpha}$. Moreover, for each torus $T^n_{\alpha}$
there exists a local coordinate system
$(\varphi_1,...,\varphi_n, I_1,...,I_n)$ in a neighborhood 
$W_{\alpha}$
of $T^n_{\alpha}$ such that $\varphi_i\in {\R}/2\pi {\Z},
(I_1,...,I_n)\in {\R}^n$ and
$\omega=\sum_{1\le i\le n} dI_i\wedge d\varphi_i$.
These coordinates are called action-angle coordinates. The
map $\pi$ in action-angle coordinates is given
by the projection 
$(\varphi_1,...,\varphi_n, I_1,...,I_n)\mapsto (I_1,...,I_n)$.
There is an ambiguity in the choice of action-angle coordinates.
In particular action coordinates $I=(I_1,...,I_n)$ are defined
up to a transformation
$I^{\prime}=A(I)+b, A\in GL(n,{\Z}), b\in {\R}^n$.
Indeed, the free abelian group 
generated by $1$-forms $dI_i, 1\le i\le n$ in each cotangent space 
$T_x^{\ast}B_0$ admits
an invariant description. It is the free abelian group
generated by the restrictions of $1$-forms $\int_{\gamma}\omega$
to $T_x^{\ast}B_0$,
where $\gamma$ runs through closed singular $1$-chains
in $\pi^{-1}(x)\cap W_\alpha$. In this way we obtain a $\Z$-affine
structure on $\pi(W_{\alpha})$.

Let $B$ be the set of connected components of
fibers of $\pi$. Endowed with the natural topology it becomes
a locally compact Hausdorff space, projection from $X$ to $B$  will be denoted by the same letter $\pi$. The natural continuous
map $B\to B_0$ is a kind of ``ramified finite covering".
Let us define $B^{sm}\subset B$ as the set
of  connected components on which  $\pi$ is a submersion (i.e. the set of all Liouville tori). Then $B^{sm}$
is an open dense subset in $B$.
Hence it carries a $\Z$-affine structure given by the action
coordinates.

The singular part $B^{sing}=B\setminus B^{sm}$ 
consists of projections of singular
fibers. Typically the codimension of $B^{sing}$ is greater or equal to $1$.
The codimension $1$ stratum consists of the boundary of the image of $\pi$ and of the 
ramification locus of the map $B\to B_0$.
The structure of singularities of the integral affine structure in higher codimensions is less understood.
 It seems that the following property
is always satisfied:
\newline{\bf Fixed Point property }. For any $x\in B^{sing}$ there 
is a small neighborhood
$U$ such that the monodromy representation $\pi_1((U\setminus B^{sing})_\alpha)\to
GL(n,{\Z})\ltimes {\R}^n$ for any connected component $(U\setminus B^{sing})_\alpha$ of $U\setminus B^{sing}$
has a fixed vector in ${\R}^n$
in the natural representation by affine transformations.

We will discuss this property in Section 6 devoted
to compactifications.

\subsubsection{Cohomological interpretation of class $[\rho]$}

In Section 2.2 we introduced an invariant $[\rho]\in H^1(B^{sm}, T^{\Z}\otimes {\R})$ 
of a $\Z$-affine
structure. Here we will give an interpretation of $[\rho]$ for integrable systems.

Let us consider $X^{\prime}=\pi^{-1}(B^{sm})$ which is a
Lagrangian torus fibration over $B^{sm}$ (i.e. fibers are
Lagrangian tori such that the fiber over $x\in B^{sm}$ is isomorphic
up to a shift to the torus $T_x^*B^{sm}/(T_x^*B^{sm})^{\Z}\,$).

Any singular closed $1$-chain $c$ on $B^{sm}$ with values
in the local system 
$$(T^{\ast}_xB^{sm})^{\Z}\simeq H_1(T_x^*B^{sm}/(T_x^*B^{sm})^{\Z}, {\Z})$$
gives a $2$-chain $\overline{c}$ on $X^{\prime}$ with the
boundary belonging to a finite collection of fibers
$\pi^{-1}(x^{(i)}), 1\le i\le N$ of the fibration $\pi: X^{\prime}\to B^{sm}$.
Moreover, for every point $x^{(i)}$ the part of $\partial \overline{c}$ over $x^{(i)}$
is homologous to zero in $\pi^{-1}(x^{(i)})$. Therefore, there exists
a collection of $2$-chains $\overline{c}_i, 1\le i\le N$ supportred on $\pi^{-1}(x^{(i)}))$ 
such that
the $2$-chain $\overline{c}+\sum_{1\le i\le N}\overline{c}_i$ is
closed. In this way we obtain a group homomorphism
$J_s: H_1(B^{sm}, (T^{\ast})^{\Z})\to H_2(X^{\prime},{\Z})/H_2^{0}(X^{\prime},{\Z})$,
where $H_2^{0}(X^{\prime},{\Z})\subset H_2(X^{\prime},{\Z})$  denotes
the sum of images of $H_2(\pi^{-1}(y), {\Z})$ where $y\in B^{sm}$
 (it is enough to pick one base point $y$ for any connected component of $B^{sm}$).
 It is easy to see
that $\langle [\rho], [c]\rangle=\langle [\omega],J_s([c])\rangle$,
where $[\omega]$ is the class of the symplectic form $\omega$.

\subsection{Examples of integrable systems}

We describe here few examples  related to the
rest of the paper.

\subsubsection{Flat tori}

First example is the triple  $(X,\pi,B_0)$ 
where $X={\R}^{2n}/\Lambda,\,\, B_0={\R}^n/\Lambda^{\prime}$ are tori
(here $\Lambda\simeq {\Z}^{2n},\,\, \Lambda^{\prime}\simeq {\Z}^n$ 
are lattices), projection $\pi:X\to B_0$ is an affine map of tori,
and $X$ carries a constant symplectic form. 
Assuming that fibers of $\pi$ are connected we have
$B_0=B=B^{sm}$. The monodromy representation
is a homomorphism 
$\rho: \pi_1(B)\to {\R}^n\subset GL(n,{\Z})\ltimes {\R}^n$.
 Integral affine structure on $B$ depends on $n^2$ real parameters,
which are coefficients of an invertible $n\times n$ matrix
expressing a basis of the lattice $\Lambda^{\prime}\subset T_xB$
as a linear combination of generators of the lattice
$(T_xB)^{\Z}\subset T_xB$, where $x\in B$ is an arbitrary
point.

\subsubsection{Surfaces}

Let $(X,\omega)$ be a surface and $\pi:X\rightarrow B_0=\R$ be an arbitrary
 smooth proper function with isolated critical points. Then $(X,\pi,B_0)$
is an integrable system. Space $B$ of connected components of fibers is a
 graph, and $\Z$-affine structure on $B^{sm}\subset B$ gives a length element
 on edges of $B$.

\subsubsection{Moment map}

Consider a compact connected 
symplectic manifold 
$(X,\omega)$ of dimension $2n$ together with
a Hamiltonian action of the torus $T^n$. Then 
 one has an integrable system $\pi:X\to B_0$, where
$\pi$ is the moment map of the action and 
$B_0=(Lie(T^n))^{\ast}\simeq {\R}^n$. Furthermore, it is well-known 
 that $B=\pi(X)$ is a convex polytope and $B^{sm}$
is the interior of $B$.

\subsubsection{K3 surfaces}

Before considering this example let us remark that one can
define integrable systems in the case of complex manifolds.
More precisely, assume that  $X$ is a complex manifold of complex
dimension $2n$,  $\omega_\C$ is a holomorphic closed non-degenerate
$2$-form on $X$, $B=B_0$ is  a complex manifold of dimension $n$
and $\pi:X\to B$ is a surjective proper holomorphic map such that generic fibers of $\pi$ 
 are connected complex Lagrangian submanifolds of $X$.
With a complex integrable system one can associate a real one
by forgetting complex structures on $X$ and $B$ and taking
$\omega:=Re(\omega_\C)$ as a symplectic form on $X$. It is easy to see that
the image of the monodromy representation belongs to
$Sp(2n,{\Z})\ltimes {\R}^{2n}\subset GL(2n,{\Z})\ltimes {\R}^{2n}$.

Let $(X,\Omega)$  be a complex K3 surface  equipped
with a non-zero holomorphic 
2-form $\omega_\C=\Omega$ and 
$\pi:X\to {\C P}^1$  a holomorphic fibration
such that the  generic fiber of $\pi$ is an elliptic curve.
For example, $X$ can be represented
as a surface in ${\C P}^2 \times {\C P}^1$ given by a
general equation $F(x_0,x_1,x_2,y_0,y_1)=0$ of 
bidegree $(3,2)$ in homogeneous coordinates. Map $\pi$ is the projection 
to the second factor.
Holomorphic form $\Omega$ is given by 
$$\Omega=i_{Euler_x\wedge Euler_y}\frac{dx_0\wedge dx_1\wedge dx_2\wedge dy_0\wedge dy_1}{dF}\,\,,$$
where $Euler_p$ denotes the Euler vector field along coordinates $p=(x_i)$ or $(y_i)$.
Such an elliptic fibration gives an integrable system. Namely,
we set $X:=X({\C})$, $\omega:=Re(\Omega)$, $B:={\C P}^1\simeq S^2$.
Generically $B^{sing}$ 
is a set of $24=\chi(X)$ points in $S^2$.
Singularity of the affine structure near each of  $24$ points 
is well-known in the theory of integrable systems
where it is called focus-focus singularity (see e.g. [Au], [Zu]).
We will discuss it in Section 6.4. Here we give
 a short description of this singularity.
We take ${\R}^2$ with the standard integral affine structure, remove
the point $(x_0,0)$ on the horizontal axis. 
Then we modify the affine structure (and also the $C^\infty$-structure!)
on the ray $\{(x,0)\,|\,\,x>x_0\}$. New local
integral affine coordinates near points of this ray will be functions
$y$ and $x+\max(y,0)$ (see Figure 1). The monodromy of the resulting integral affine structure
 around removed singular point $(x_0,0)$ is given
by the transformation $(x,y)\mapsto (x+y,y)$.

\begin{figure}\label{figure1}
\centerline{\epsfbox{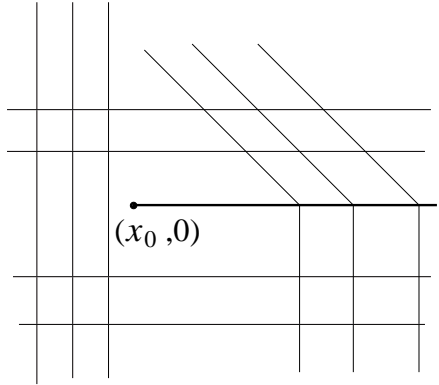}}
\caption{Focus-focus singularity. All lines are straight in the modified $\Z$-affine structure.}
\end{figure}

\subsection{Families of integrable systems and PL actions}

In many examples 
an integrable system depends on parameters. It often happens that
the parameter space
${\cal P}$ carries a natural foliation ${\cal F}$ such that the fundamental
group $\pi_1({\cal F}_p,p), p\in {\cal P}$ of any leaf acts on 
the base space $B_p$
of the corresponding torus fibration. This action
is given by piecewise-linear homeomorphisms
with integral linear parts.

Let us illustrate this phenomenon
in the case of the family of integrable systems associated
with a K3 surface discussed above.

Here the parameter space ${\cal P}$ has dimension $38$, which is
twice of the complex dimension of the space of polynomials $F$ 
modulo unimodular linear transformations.
On the other hand,  the miniversal family of  representations
(up to a conjugation)
$$\left\{\,\rho: \pi_1(S^2-\{24\mbox{ points}\})\ra SL(2,\Z)\ltimes \R^2\,\right\}$$
such that the monodromy around each puncture is conjugate to 
$\left(\begin{array}{cc} 1 & 1\\0& 1
\end{array}\right)$, 
 has dimension  $20$.

Thus, we obtain a foliation ${\cal F}$ of ${\cal P}$ of rank $18=38-20$. 
It is defined by the following property:
if we continuously vary parameters $p\in {\cal P}$ along leaves
of ${\cal F}$ then the conjugacy class of the 
monodromy representation $\rho$ remains unchanged.

Notice that in the local model 
described above we can
move the position $(x_0,0)$ at which we start the cut. 
Then we have on the sphere $S^2$ a  set of $24$
``worms'' (singular points, each of them can move in its preferred direction, which
is the line invariant under the local monodromy). One can show easily that any 
 continuous  deformation of $\Z$-affine structure satisfying Fixed Point property (see Section 3.1)
 and preserving the conjugacy class of $\rho$, corresponds to a movement of worms.
\footnote{Notice that in our example $rk({\cal F})=18$ is less than $24$.
 This means that there  are 6 constraints on moving worms.}

Moving ``worms''we get a canonical
identification of manifolds with integral affine structures
far enough from singular points. We will see later in Section 6.4 that we also have
a canonical PL identification of manifolds near singular points.
Therefore we obtain a local system along leaves of ${\cal F}$ 
with the fiber over $p\in {\cal P}$
being a manifold $B_p\simeq S^2$ with the above 
$\Z$-affine structure.
In this way we get a homomorphism
from $\pi_1({\cal F}_p,p)$ to $Aut_{\Z PL}(S^2)$, where $\Z PL$ denotes the 
group of integral PL transformations of $S^2$ equipped with the above
$\Z$-affine structure. We will return to this action in
Section 6.7 where it will be compared with another PL action
on the same space.

\section{B-model construction}

\subsection{$\Z$-affine structure on smooth points}

Here we are going to define an analog of the notion
of integrable system in the framework of rigid analytic
geometry. Roughly speaking, it is a triple $(X,\pi, B)$, where
$X$ is a variety defined over a non-archimedean field
(see [Be1] and Appendix A), $B$ a CW complex and
$\pi:X\to B$ a continuous map.
More precisely, let $K$ be a field with non-trivial
valuation, $X$ an irreducible algebraic variety over $K$ of 
dimension $n$, $f=(f_1,...,f_N)$ a collection
of non-zero rational functions on $X$. Then we have a multivalued
map 
$$X(\overline{K})\to [-\infty,+\infty]^N,\,\,
x\mapsto val_{\overline{K}}(f(x)):=\left(val_{\overline{K}}\left(f_1(x)),\dots
,val_{\overline{K}}(f_N(x)\right)\right)\,.$$
Here $\overline{K}$ is the algebraic closure of $K$,
and $val_{\overline{K}}$ denotes valuation 
on $\overline{K}$.

Let $\psi:[-\infty,+\infty]^N\to B$ be a continuous map
such that the composition $\pi=\psi\circ val(f)$ is single-valued.
Our map $\pi$ will always be of this form.
More generally, we can take $X$ to be a (not necessarily algebraic)
compact smooth $K$-analytic  space, and $\pi:X\to B$ be
a continuous map which factorizes as the composition
of the projection
$p_{\cal X}:X\to S_{\cal X}$ to the Clemens polytope 
$S_{\cal X}$ of some model ${\cal X}$
of $X$ and a continuous map 
$\pi':S_{\cal X}\to B$ (see  Section 4.2.3).

Now we would like to be more precise.
Let $K$ be a complete non-archimedean field, with valuation $val$
 and the corresponding 
 norm $|x|:=\exp(-val(x))\in \R_{\ge 0}$. 
 Before giving next definition we observe that there is
a canonical continuous map $\pi_{can}: ({\bf G}_m^{an})^n\to {\R}^n$
(see Section A.2 in Appendix A). Here
${\bf G}_m^{an}$ is a multiplicative group (considered as an analytic space
over $K$) and the restriction
of $\pi_{can}$ to $(\overline{K}^{\times})^n$ 
is given by the formula 
$$\pi_{can}(z_1,...,z_n)=(\log |z_1| ,\dots,\log |z_n| )\,\,.$$
The sheaf 
${\cal O}^{can}_{{\R}^n}:=(\pi_{can})_{\ast}({\cal O}_{({\bf G}_m^{an})^n})
$ of $K$-algebras
is called the {\it canonical} sheaf.

Let $X$ be a smooth
$K$-analytic space of dimension $n$, 
$\pi:X\to B$ a continuous map of $X$ into a Hausdorff 
topological
space $B$.

\begin{defn} We call a point $x\in B$ smooth
(or $\pi$-smooth) if there is a neighborhood
$U$ of $x$ such that the fibration $\pi^{-1}(U)\to U$ is 
isomorphic
to a fibration $\pi_{can}^{-1}(V)\to V$ 
for some open subset $V\subset {\R}^n$.
Here the isomorphism $\pi^{-1}(U)\simeq \pi_{can}^{-1}(V)$
is taken in the category of $K$-analytic spaces while
$U\simeq V$ is a homeomorphism.

In this case we will call $\pi$ (or the triple $(\pi^{-1}(U),\pi, U)$)
an analytic torus fibration.

\end{defn}

Let $B^{sm}$ denotes the set of smooth points of $B$.
It is a topological subspace
of $B$ (in fact a topological manifold of dimension $n$).

\begin{thm} The space $B^{sm}$ carries a sheaf of $\Z$-affine functions,
which is locally isomorphic to the canonical sheaf of $\Z$-affine
functions on ${\R}^n$.

\end{thm} 

{\it Proof.} We start with the following Lemma.

\begin{lmm} Let $V\subset {\R}^n$ be a 
connected open set,
$\varphi\in {\cal O}^{\times}_{{({\bf G}_m^{an})}^n}(\pi_{can}^{-1}(V))$ be an invertible
analytic function. Then the function $val_x(\varphi(x))$ is 
constant along fibers of $\pi_{can}$, and it is a pull-back
of a $\Z$-affine function on ${\R}^n$.

\end{lmm}

In order to prove Lemma we observe that any analytic function 
$\psi\in {\cal O}^{\times}_{{({\bf G}_m^{an})}^n}(\pi_{can}^{-1}(V))$
can be decomposed into Laurent series:
$$\psi=\sum_{I=(i_1,\dots,i_n)\in {\Z}^n}c_I z^I,\,\,\, c_I\in K$$ satisfying
certain convergence conditions (see Section A.2).

Then for a non-zero analytic function $\psi$ on $\pi_{can}^{-1}(V)$  we  introduce a real-valued function 
$Val(\psi)(x):=\inf_{I\in {\Z}^n}(val(c_I)-\langle I,x\rangle),
x\in V$.
It is a concave, locally piecewise-linear function on $V$.
It is easy to see that
\begin{description}
\item[a)] there is a dense open subset $V_1\subset V$ such that 
for any $x\in V_1$ the infimum
in the definition of $Val(\psi)$ is achieved
for a single multi-index $I$;

\item[b)] $Val(\psi_1\psi_2)=Val(\psi_1)+Val(\psi_2)$.
\end{description}
For an invertible function $\varphi$ we have
$0=Val(1)=Val(\varphi \varphi^{-1})=Val(\varphi)+Val(\varphi^{-1})$.
Since both $Val(\varphi)$ and $Val(\varphi^{-1})$ are concave, their sum can be
equal to zero iff they are both affine. Moreover they are both
$\Z$-affine since the linear part of $Val(\varphi)$
is given by the integer vector $I$ for some
single multi-index $I$.
Finally, observe that 
$val_x(\varphi(x))\ge\pi_{can}^*(Val(\varphi))(x), 
x\in \pi_{can}^{-1}(V)$. Therefore $val_x(\varphi(x))=Val(\varphi)(\pi_{can}(x))$
 for invertible $\varphi$. $\blacksquare$

Now we can finish the proof of the Theorem. The above
formula gives us a coordinate-free description of
$\pi_{can}^*(Val(\varphi)) $.
It is easy to see that any $\Z$-affine function on $V$ is of the form
$Val(\varphi)+c, c\in {\R}$ for some invertible $\varphi$ 
(in the case of ${\R}^n$ 
it suffices to take monomials as $\varphi$).
We can identify $\pi^{-1}(U)\to U$ with $\pi_{can}^{-1}(V)\to V$
for some small open $U\subset X$ and $V\subset {\R}^n$. Then we can define
$Val(\varphi)$ for any invertible $\varphi\in {\cal O}_X(\pi^{-1}(U))$
by the above formula.
Finally we define a sheaf of $\Z$-affine functions on $B^{sm}$ by taking
all functions of the form $Val(\varphi)+c, c\in {\R}$.
It follows from the above discussion that in this way we obtain
a $\Z$-affine structure on $B^{sm}$, which is locally isomorphic
to the standard one on ${\R}^n$. $\blacksquare$

We will denote by $Aff^{can}_{{\Z},B^{sm}}$ the sheaf of $\Z$-affine
functions constructed in the proof.

\subsection{Examples}

\subsubsection{Logarithmic map}

This is a basic example
$$\pi=\pi_{can}=\log |\cdot|:X=({\bf G}_m^{an})^n\to B_0=B={\R}^n$$
described in details in Appendix A. For any algebraic (or analytic)
subvariety $Z\subset ({\bf G}_m^{an})^n$ of dimension $m\le n$ its image
 $\pi(Z)$  is a non-compact piecewise-linear closed subset of $\R^n$
 of real dimension $m$. Smooth points for $\pi_{|Z}$ are dense in $\pi(Z)$.

              In particular, if $Z$ is a curve then  $\pi(Z)$
               is a graph in $B$ with straight edges having rational
               directions. 
               One can try to make a dictionary which translates the properties
 of the algebraic variety $Z\subset {\bf G}_m^n$ to the properties
 of the PL set $\pi(Z^{an})$ which is the closure
of  $\pi(Z(\overline{K}))$ in $\R^n$. This circle of ideas
 is a subject  of the so-called  ``tropical geometry" 
 (see e.g. [Mi]).

\subsubsection{Tate tori}

Let $\rho:{\Z}^n\to (K^{\times})^n$ be a group homomorphism
such that the image of the composition $val\circ \rho:{\Z}^n\to {\R}^n$
is a rank $n$ lattice in ${\R}^n$. Group $(K^{\times})^n$ acts
by translations on the analytic space $({\bf G}_m^{an})^n$. Restriction
of this action to ${\Z}^n$ (via $\rho$) is discrete and cocompact.
The quotient is a $K$-analytic space $X$ called {\it Tate torus}.
There is an obvious map $\pi:X\to B:={\R}^n/(val\circ\rho)({\Z}^n)$.
All points of $B$ are smooth. The space $X$ depends on
$n^2$ parameters taking values in $K^{\times}$(cf. with the flat tori example
in Section 3.2.1).

\subsubsection{Clemens polytopes and their contractions}

For any smooth projective variety $X$ of dimension $n$,
 and  and a snc model ${\cal X}$ of it (see Appendix A)
 we have a canonical projection to the corresponding Clemens polytope
 $$p_{{\cal X}}:X^{an}\ra S_{{\cal X}}\,\,.$$
All interior points of $n$-dimensional simplices of $S_{{\cal X}}$
 are $p_{{\cal X}}$-smooth, although there might be other smooth
points too. More generally, one can compose projection $p_{{\cal X}}$
 with a continuous surjection $\pi':S_{{\cal X}}\epi B$ where
 $B$ is a finite CW complex and map $\pi'$ is a cell map for some cell subdivision
 of $S_{{\cal X}}$. We assume that fibers of the composition 
 $\pi:=\pi'\circ p_{{\cal X}}: X^{an}\to B$ are connected.
 This seems to be the most general case of maps from projective varieties
 over complete local fields to CW complexes relevant for our purposes.

\subsubsection{Curves}

Let $X/K$ be a connected smooth projective curve of genus $g>1$.
After passing to a finite extension $K^{\prime}$ of $K$ we
may assume that $X$ has a canonical model 
${\cal X}$ with stable reduction. The graph $\Gamma^{\prime}$
corresponding to the special fiber ${\cal X}_0$ is a retraction
of $(X\otimes_K K^{\prime})^{an}$. The quotient graph 
$\Gamma=\Gamma^{\prime}/Gal(K^{\prime}/K)$ is a retraction
of the  analytic curve $X^{an}$ (see [Be1]).
We define $B:=\Gamma$. Then $B^{sm}$ is a complement
to a finite set. As in Section 3.2.2, a  $\Z$-affine structure on a graph
 is the same as 
a length element  (i.e. a metric). Therefore $\Gamma$ is a metrized
graph. Notice also that the maximal number of edges of the graph
corresponding to a genus $g$ curve is $3g-3$, which is the dimension
of the moduli space of genus $g$ curves.

Notice that if in Section 4.2.1 subvariety $Z$ is a curve then its projection
 is a noncompact metrized graph with unbounded  edges corresponding to punctures
 ${\overline{Z}}\setminus Z$.
\subsubsection{K3 surfaces}

Here we will describe a particular case of the construction from Section
4.2.3 (a contraction
 of a Clemens polytope).

Let field $K$ be ${\C}((t))$ and $X\subset {\bf P}^3_K$ be
a formal family of complex K3 surfaces given by the equation

$$x_0x_1x_2x_3+tP_4(x_0,x_1,x_2,x_3)=0\,\,,$$
where $P_4$ is a generic homogeneous polynomial of degree four,
and $t$ is a formal parameter.

The special fiber at $t=0$ of this family is singular, it is given by the equation
$x_0x_1x_2x_3=0$. Let us denote by $\widetilde{{\bf P}^3}$ the blow-up
of the total space of the trivial ${\bf P}^3$-bundle over $Spec(\O_K)$
at $24$ points $p_{\alpha}, 1\le \alpha\le 24$ of the special fiber,
where each $p_{\alpha}$ is a solution of the equation
$$P_4(x_0,x_1,x_2,x_3)=0,\,\, x_i=x_j=0,\,\, 1\le i<j\le 4\,\,.$$
The closure ${\cal X}$ of $X$ in $\widetilde{{\bf P}^3}$ is a model
with simple normal crossings.
The associated Clemens polytope $S_{\cal X}$ has $28$ vertices.
Four of them correspond to coordinate hyperplanes $x_i=0$ in 
${\bf P}^3$, and $24$ other correspond to divisors sitting at
the pre-images of the points $p_{\alpha}$. 
Therefore $S_{\cal X}$ is the union of the boundary
$\partial\Delta^3$ of the standard $3$-simplex $\Delta^3$ with
$24$ copies of the standard $2$-simplex $\Delta^2$.
Those $24$ triangles $\Delta_{\alpha}^2, 1\le \alpha\le 24$ are
decomposed into six groups of four triangles in each. All triangles
from the same group have a common edge, which is identified
with an edge of $\partial\Delta^3$ (tetrahedron with $24$ ``wings").
As we mentioned in the previous example, there is a continuous
map $p:X^{an}\to  S_{\cal X}$. We are going to construct $B$ as
a retraction of $S_{\cal X}$.

In order to do this we observe that for an edge $e\subset \Delta^2$
and a point $a\in e$
one has the canonical retraction $p_{a,e}:\Delta^2\to e$.
Namely, let us identify the edge $e$ with the interval $[-1,1]$
of the real line, so that $a$ is identified with the point
$a=(a_0,0)$, and $\Delta^2$ is bounded by $e$ and the segments
$0\le y\le 1-|x|$. Then we define $p_{a,e}$ by the formulas (see Figure 2)
$$\begin{array}{llcl}
(x,y) & \mapsto & (x+y,0)\,,\,\,\, & x+y\le a_0\,\,; \\
(x,y) & \mapsto & (x-y,0)\,,\,\,\, & x-y\ge a_0\,\,; \\
(x,y) & \mapsto & (a_0,0)\,,\,\,\, & \mbox{ otherwise.}
\end{array}$$

\begin{figure}\label{figure2}
\centerline{\epsfbox{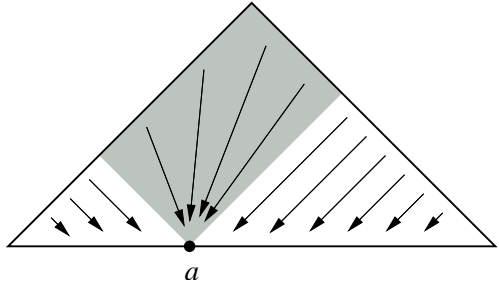}}
\caption{Triangle contracted to one side. The dashed area maps to point $a$.}
\end{figure}

Now we choose a point $q_{ij}, 0\le i<j\le 3$ in the interior
of each edge $e_{ij}$ of $\partial\Delta^3$ (here $i,j$
are identified with the vertices of $\partial\Delta^3$).
There are four ``wings"
$\Delta_{\alpha}^2$ having $e_{ij}$ as a common edge. Then we retract
each $\Delta_{\alpha}^2$ to $e_{ij}$ by the map $p_{q_{ij},e_{ij}}$.
This gives us a retraction $\pi'=\pi'_{(q_{ij})}:S_{\cal X}\to \partial\Delta^3$.
Let $\pi:=p_{\cal X}\circ \pi'_{(q_{ij})}:X^{an}\to B$ be the composition
of the projection $p_{\cal X}:X^{an}\to S_{\cal X}$ with the above
retraction.
One can show that all points of $B:=\partial\Delta^3$ are $\pi$-smooth
except of the chosen six points $q_{ij}, 0\le i<j\le 3$. 
According to Theorem 1 we obtain a $\Z$-affine structure on 
$S^2\setminus \cup_{1\le i<j\le 3}\{q_{ij}\}$. One can show that the local
monodromy around each point $q_{ij}$ is conjugate to the matrix
$$\left(\begin{array}{cc} 1 & 4\\0& 1
\end{array}\right)\,\,.$$
We skip the computations here.

\subsection{Stein property}

A $K$-analytic space $X$ is called {\em Stein} if the natural map
  $$X\to Spec^{an}(\Gamma(X,{\cal O}_X))$$ is a homeomorphism. Here $\Gamma(X,{\cal O}_X)$
 is considered as a topological $K$-algebra. This definition is equivalent to the standard one.
Let us call the projection $\pi:X\to B$ {\em Stein} if
for any $b\in B$ there exists  a fundamental systems of neighborhoods $U_i$ of $x$ 
such that $\pi^{-1}(U_i) \subset X$ is a Stein domain.
 If $\pi$ is Stein then we can reconstruct $(X,{\cal O}_X)$ and $\pi$ from the space $B$ endowed with 
the sheaf
 $\pi_*({\cal O}_X)$ of topological $K$-algebras.

\begin{prp} Let $B$ be a contraction of Clemens polytope $S_{{\cal X}} $ of 
some model $\cal{X}$ of $X$
as in Section 4.2.3, and $\pi$ a Stein map.
Then $B^{sm}$ is dense in $B$.

\end{prp}

{\it Proof.}{\footnote{We thank to Ofer Gabber
for suggesting the proof below}} 
 It suffices to prove that $n$-dimensional cells are dense in $B$, where $n=\dim X$.
 For any open  $U\subset B,\,\,U\ne \emptyset$ 
 we have $H_c^n(U,\pi_{\ast}(\Omega^n_X))\simeq  H_c^n(\pi^{-1}(U),\Omega^n_X)$.
 
The last group is nontrivial, because for any non-empty open
 $V\subset X^{an}$ the integration map
 $\int: H_c^n(V,\Omega^n_X)\to K$ is onto.
 Therefore $\dim(U)\ge n$.
$\blacksquare$

All the examples in Sections 4.2.1--4.2.5 (except Section 4.2.3) have Stein property.

\section{$\Z$-affine structures and mirror symmetry}

\subsection{Gromov-Hausdorff collapse of Calabi-Yau manifolds}

We recall that a Calabi-Yau metric on a complex manifold $X$ is a
K\"ahler metric with vanishing Ricci curvature. If such a metric exists
then $c_1({TX})=0\in H^2(X,{\R})$ and hence the class of the
 canonical  bundle $\bigwedge^{\dim X}(T^*X)$
  is torsion in $Pic(X)$. According to the famous Yau theorem,
for any compact K\"ahler manifold $X$ such that $c_1({TX})=0\in H^2(X,{\R})$, and
any K\"ahler class $[\omega]\in H^2(X,{\R})$ there exists a unique
Calabi-Yau metric $g_{CY}$ with the class $[\omega]${\footnote{Notice that there is a discrepancy in terminology.
 In algebraic situation one usually calls Calabi-Yau a projective variety with the
 trivial canonical class in  $Pic(X)$, 
and the polarization is not considered as a part of data.}}.
Up to now,
there is no explicitly known non-flat Calabi-Yau metric on a compact
manifold. 

In Mirror Symmetry one studies the limiting behavior of $g_{CY}$
as the complex structure on $X$ approaches a  ``cusp"
in the moduli space of complex structures (``maximal degeneration").
Well-known conjecture of Strominger, Yau and Zaslow (see [SYZ])
claims a torus fibration structure of Calabi-Yau manifolds
near the cusp. A metric approach to the maximal degeneration (see [GW],
[KoSo]) explains the structure of such Calabi-Yau manifolds
in terms of their Gromov-Hausdorff limits. We recall this
picture below following [KoSo].

We start with the definition of a maximally degenerating family
of algebraic Calabi-Yau manifolds.

Let ${\C}_t^{mer}=\{f=\sum_{n\ge n_0}a_n t^n\}$ be the field of germs
at $t=0$ of meromorphic functions in one complex variable, and
${X}_{mer}$ be an algebraic $n$-dimensional Calabi-Yau manifold 
over ${\C}_t^{mer}$ (i.e. ${X}_{mer}$ is a smooth projective manifold
over ${\C}_t^{mer}$ with the trivial canonical class: $K_{{X}_{mer}}=0$).
We fix an algebraic non-vanishing
volume element $\Omega\in \Gamma ({X}_{mer},K_{{X}_{mer}})$.
The pair $({X}_{mer},\Omega)$ defines a 1-parameter 
analytic family of complex
Calabi-Yau manifolds $(X_t,\Omega_t), 0<|t|<\epsilon$, for  some $\epsilon>0$.

Let $[\omega]\in H^2_{DR}({X}_{mer})$ be a cohomology class in the ample cone.
Then for every $t$, such that $0<|t|<\epsilon$ it defines a K\"ahler class
$\omega_t$ on $X_t$. We denote by
$g_{X_t}$ the unique Calabi-Yau metric
on $X_t$ with the K\"ahler class $[\omega_t]$.

 It follows from the resolution of singularities,
that as $t\to 0$ one has 

$$\int_{X_t}\Omega_t\wedge \overline{\Omega}_{t}=C(\log|t|)^m|t|^{2k}(1+o(1))$$
for some $C\in {\C}^{\times}, k\in {\Z}, 0\le m\le n=\dim\,({X}_{mer})$.

\begin{defn} We say that ${X}_{mer}$ has maximal degeneration at $t=0$
if in the formula above we have $m=n$.

\end{defn}

Let us rescale the Calabi-Yau metric: 
$g_{X_t}^{new}=g_{X_t}/diam(X_t,g_{X_t})^{1/2}$.
In this way we obtain a  family of Riemannian manifolds
$X_t^{new}=(X_t,g_{X_t}^{new})$ of  diameter $1$.

\begin{conj} If ${X}_{mer}$ has maximal degeneration at $t=0$ then
 $$diam(X_t,g_{X_t})=(\log|t|)^{-1}\exp(O(1))$$ and
there is a limit $(B,g_{B})$ of $X_t^{new}$
in the Gromov-Hausdorff metric as $t\to 0$, such that:
\begin{description}

\item[a)] $(B,g_{B})$ is a compact metric space,
which contains a smooth oriented 
Riemannian manifold $(B^{sm},g_{B^{sm}})$ of dimension $n$
as a dense open metric subspace.
The Hausdorff dimension of $B^{sing}=B\setminus B^{sm}$
is less  than or equal to $n-2$.

\item[b)] $B^{sm}$ carries a $\Z$-affine structure.

\item[c)] The metric $g_{B^{sm}}$ has a potential. 
This means that it is locally given 
in affine coordinates by a symmetric matrix
$(g_{ij})=(\partial^2 F/\partial x_i\partial x_j)$, where $F$ is a
smooth function (defined modulo adding an affine function).

\item[d)] In affine coordinates the metric volume element is constant, i.e.
 $$\det(g_{ij})=\det(\partial^2 F/\partial x_i\partial x_j)=const$$ 
 (real Monge-Amp\`ere equation).
\end{description}

\end{conj}

There is a more precise conjecture (see [KoSo] for the details) which says
that outside of $B^{sing}$ the space $X_t^{new}$ is metrically close
to a torus fibration with flat Lagrangian fibers 
(integrable system). This torus fibration can be canonically
reconstructed (up to a locally constant twist) from the limiting data a)-d).

Conjecture 1 holds for abelian varieties (since $B=B^{sm}$
is a flat torus in this case). It
is non-trivial for K3 surfaces (see [GW] for the proof). In
 3-dimensional case  there is now
 a substantial progress  (see [LYZ]).

\begin{defn} A Monge-Amp\`ere manifold is a triple $(Y,g,\nabla)$,
where $(Y,g)$ is a smooth Riemannian manifold with the metric $g$, and
$\nabla$ is
a flat connection on $TY$ such that:
\begin{description}
\item[a)] $\nabla$ defines an affine structure on $Y$.

\item[b)] Locally in affine coordinates $(x_1,...,x_n)$
 the matrix $(g_{ij})$
 of $g$ is given by $(g_{ij})=(\partial^2 F/\partial x_i\partial x_j)$
 for some smooth real-valued function $F$.
 
\item[c)] The Monge-Amp\`ere equation
$\det(\partial^2 F/\partial x_i\partial x_j)=const$ is satisfied.
\end{description}
\end{defn}

The following easy Proposition is well-known.

\begin{prp} For a given Monge-Amp\`ere manifold $(Y,g_Y,\nabla_Y)$
there is a canonically defined
dual Monge-Amp\`ere manifold $(Y^{\vee},g_Y^{\vee},\nabla_Y^{\vee})$
such that $(Y,g_Y)$ is identified
with $(Y^{\vee},g_Y^{\vee})$ as Riemannian manifolds, and
the local system $(T{Y^{\vee}},\nabla_Y^{\vee})$ 
is naturally isomorphic to the local system dual to  
$(TY,\nabla_Y)$ (dual local system is constructed
via the metric $g_Y$).

\end{prp}

\begin{cor} If $\nabla_Y$ defines 
an integral affine structure on $Y$ with the covariantly constant
lattice $(TY)^{\Z}$
then $\nabla_Y^{\vee}$ defines an integral affine structure on 
$Y^{\vee}$ such that for all $x\in Y^{\vee}=Y$ the lattice
$(T_xY^{\vee})^{\Z}$ is  dual to $(T_xY)^{\Z}$ with respect to the Riemannian
metric $g_Y$  on $Y$.

\end{cor}

We will call {\it integral} a Monge-Amp\`ere manifold with 
$\Z$-affine structure.

In Mirror Symmetry one often has a so-called dual family of Calabi-Yau
manifolds associated with the given one. There is no general
definition of the dual family, but there are many examples.
The following Conjecture (see [KoSo]) formalizes  Strominger-Yau-Zaslow picture of 
 Mirror Symmetry:

\begin{conj} Smooth parts of Gromov-Hausdorff limits
of dual families
of Calabi-Yau manifolds are dual integral Monge-Amp\`ere manifolds.

\end{conj}

One can say that Monge-Amp\`ere manifolds with integral affine structures are real
analogs of Calabi-Yau manifolds. 
Conversely, having an integral Monge-Amp\`ere manifold
$(Y,g_Y,\nabla_Y, (TY)^{\Z})$ one can construct a torus fibration
 $TY/(TY)^{\Z}\to Y$. It is easy to see that
the total space of this fibration
 is in fact a Calabi-Yau manifold (typically non-compact as $Y$ is non-compact too).
Rescaling the covariant lattice we can make fibers small (of the size
$O((\log|t|)^{-1})$). As we already mentioned, the extended version
of Conjecture 1 says that this torus fibration is  close (after a locally constant twist)
to $X_t^{new}$ outside of a ``singular" subset.

\subsubsection{K3 example}
In the case of collapsing K3 surfaces the corresponding  intergal Monge-Amp\`ere 
 manifold has an explicit description.

Let $S$ be a complex surface endowed with a holomorphic
non-vanishing volume form $\Omega_S$, and $\pi:S\to C$ be a holomorphic
fibration over a complex curve $C$, such that fibers of $\pi$
are non-singular elliptic curves.

We define a metric $g_C$ on $C$ as the K\"ahler metric associated
with the $(1,1)$-form $\pi_{\ast}(\Omega_S\wedge \overline{\Omega}_S)$.
Let us choose (locally on $C$) a basis $(\gamma_1,\gamma_2)$
in $H_1(\pi^{-1}(x),{\Z}), x\in C$. We define two closed 1-forms on $C$ by the
formulas

$$\alpha_i=Re\left(\int_{\gamma_i}\Omega_S\right),\,\,\, i=1,2\,\,.$$

It follows that $\alpha_i=dx_i$ for some
functions $x_i, i=1,2$. We define a $\Z$-affine structure on $C$,
and the corresponding connection $\nabla$,
by saying that $(x_1,x_2)$ are $\Z$-affine coordinates (compare with 3.2.4).
One can check directly that $(C,g_C,\nabla)$ is a Monge-Amp\`ere
manifold. In a typical example of elliptic fibration
of a K3 surface, one gets $C={\C P}^1\setminus \{x_1,...,x_{24}\}$,
where $\{x_1,...,x_{24}\}$ is a set of distinct $24$ points in 
${\C P}^1$. M.~Gross and P.~Wilson (see [GW]) proved that there exists
 a family of K3 surfaces with Calabi-Yau metrics collapsing to $S^2\simeq \C P^1$
 with the intergal Monge-Amp\`ere structure described above.

\subsection{Non-archimedean picture for the space $B$}

Here we would like to formulate a conjecture which relates the 
Gromov-Hausdorff limit with non-archimedean geometry, thus giving
a pure algebraic description of  $\Z$-affine structure on $B^{sm}$.
Let $\overline{{\C}_t^{mer}}=\cup_{m\ge 1}{\C}_{t^{1/m}}^{mer}$ 
be the algebraic closure of  ${\C}_t^{mer}$.
We denote by $\pi_{mer}:X(\overline{{\C}_t^{mer}})\to B$ the map which
associates the limiting point (in Gromov-Hausdorff metric) 
of points $x(t^{1/m})\in X_{t^{1/m}}({\C})$ as $t^{1/m}\to 0$.

Let ${K}={\C}((t))$ be the field of Laurent formal
series. Then, by extending scalars we obtain an algebraic
Calabi-Yau manifold ${X}$ over ${K}$.
We denote by ${X}^{an}$ the corresponding smooth
${K}$-analytic space.

\begin{conj} The map $\pi_{mer}$ is well-defined and extends by
continuity to the map $\pi:{X}^{an}\to B$.
The set $B^{sm}$ (defined as the maximal open subset of $B$ on which the limiting metric is smooth)
 coincides with the set of
$\pi$-smooth points. Two $\Z$-affine structures on $B^{sm}$, one coming from the collapse picture, another coming from
 non-archimedean picture, coinside with each other.

\end{conj}

Also we make the following conjecture (or better a wish, because it is based on 
a very thin evidence):
\begin{conj} Map $\pi$ is Stein.
\end{conj}

\part{}

\section{Compactifications of $\Z$-affine structures}

\subsection{Properties of compactifications}

Assume that we are given a non-compact
manifold $B^{sm}$ with a $\Z$-affine 
structure.
We would like to ``compactify" it, i.e. to find a compact
Hausdorff topological space $B$ such that $B^{sm}\subset B$ is
an open dense subset. We do not require an extension of the
$\Z$-affine structure to $B$.
 The question is: what kind of properties one
should expect from such a compactification? 
We cannot give a complete
list of such properties at the moment. Instead, we formulate two of them
and illustrate the notion of compactification in PL case.
Similarity between examples in Sections 3.2 and 4.2 suggests
 that the  class of singularities which appear
 in integrable systems should be more or less the same as  the class of singularities
 appearing in non-archimedean geometry. 

Let $x\in B^{sing}:=B\setminus B^{sm}$ 
be a singular point of some compactification
of $B^{sm}$. Then we require the following
\newline{\bf Finiteness property}. There is a fundamental system of   neighborhoods  
$U\subset B$ of $x$
such that the number of connected components 
of $U\cap B^{sm}$ is finite.

Let $U\cap B^{sm}=\sqcup_{1\le i\le N}U_i$ be the 
disjoint union of the connected components. Let us pick a
point $x_i\in U_i$  and consider a continuous path 
$\gamma:[0,1]\to B$ such that $\gamma(0)=x_i, \gamma(1)=x,
\gamma([0,1))\subset U_i$. Using the affine structure
we can canonically
lift this path to a path
$\gamma':[0,1)\to T_{x_i}B$. We assume that the lifted path $\gamma'$ extends to 
time $t=1$ and is analytic at $t=1$
(it is a technical assumption,
 helping to avoid
 some pathologies). Then we require the following
\newline{\bf Independence property}. Path $\gamma$ with the properties as above exists, and  
point $\gamma'(1)\in T_{x_i}B$
does not depend on the choice of $\gamma$.

Independence property implies the existence of a fixed vector for
the monodromy representation restricted to
$\pi_1(U\cap B^{sm}, x_i)$ (this implies the Fixed Point property from Section 3.1).

\subsection{PL compactifications}

Let $V$ be a finite set, $S\subset 2^V$ belongs to 
the set of $(n+1)$-element subsets of $V$. Then we have
a $n$-dimensional  simplicial complex $B=\cup_{Y\in S}\Delta^Y\subset \Delta^V$.

Let us choose a $\Z$-affine structure on $n$-dimensional faces of $B$
which is compatible with the standard affine structure,  and
consider all $(n-1)$-dimensional faces which enjoy the 
following property: they belong to exactly two $n$-dimensional
faces. For any two such $n$-dimensional faces $\sigma$ and $\tau$
we choose a $\Z$-affine structure on $\sigma\cup \tau$ which is
compatible with the already chosen $\Z$-affine structures
on $\sigma$ and $\tau$ (such a choice is equivalent to a
choice of $\Z$-affine structure in a neighborhood of the 
$(n-1)$-dimensional face $\sigma\cap \tau$). In this way we obtain
 a $\Z$-affine structure on the union $U$ of the interior points of all $n$-dimensional simplices
 and also  the interior points of $(n-1)$-dimensional faces belonging to exactly two top-dimensional cells.

\begin{prp} There exists (and unique) maximal extension of this $\Z$-affine structure
to an open subset $U_{\max}\subset B$ containing $U$.

\end{prp}

{\it Proof.} Let us proceed inductively by codimension of faces. 
The induction step reduces to the obvious remark that the extension
 of the standard $\Z$-affine structure on ${\R}^n\setminus L$ to a neighborhood
 of point $p\in L$ in $\R^n$ is unique in the case when $L\subset \R^n$
 is an affine subspace, $\dim L\le n-2\,\,$. $\blacksquare$

It is easy to see that $B^{sm}:=U_{\max}$ with $\Z$-affine structure on it, compactified  by $B$
 satisfies both Finiteness and Independence properties.

We introduce PL compactifications both as a ``toy model'', and also (as we hope, 
see Conjecture 6 in Section 6.3) as
 a sufficiently representative class for applications.
 In this case we can try to formulate additional desired properties. One of goals is to find a good
 substitution for the algebro-geometric
 notion of a canonical singularity (which is, morally,  a singularity of a non-collapsing limit of a family 
of Calabi-Yau 
manifolds with fixed K\"ahler class).

For a large class of maximally degenerating Calabi-Yau manifolds there is a proposal
 by several authors (see [GS] and [HZh])
 for a PL compactification $B$ conjecturally related
 to the Gromov-Hausdorff limit. Space $B$ is topologically a sphere $S^n$,
 it carries two dual cell decompositions. Each of these decompositions is
 identified with the boundary $\partial P_1$ or $
\partial P_2$ of a
 convex $(n+1)$-dimensional polytope. Moreover, on each $n$-dimensional face
 of each polytope we have a $\Z$-affine structure compatible with 
the natural affine structure. The assumption is that for any two open $n$-cells
 $U, U'$ from the first and the second cell decompositions, two induced 
$\Z$-affine structures on $U\cap U'$ coincide. This gives a $\Z$-affine structure
 on $B\setminus \left(Sk_{n-1}\cap Sk_{n-1}'\right)$ where 
$Sk_{n-1}$ and $ Sk_{n-1}'$
 are $(n-1)$-skeletons of two CW-structures.

\subsection{Some conjectures about singular sets}

Our conjectures are in fact rather ``wishes", i.e. 
they are desired properties of $B^{sing}=B\setminus B^{sm}$.
For simplicity we
assume that $B^{sing}$ is a stratified set (say, CW complex)
of dimension less or equal than $n-1$.

\begin{conj} We have a decomposition 
$B^{sing}=B^{sing}_{n-1}\cup B^{sing}_{\le n-2}$, where
$B^{sing}_{\le n-2}$ consists of strata of dimension
less or equal than $n-2$, $B^{sing}_{n-1}$ is the union
of strata of dimension $n-1$, and locally
near every point $x\in B^{sing}_{n-1}$ the $\Z$-affine structure
is modeled by the ``book"
$\cup_{i\in I}{\R}^{n-1}\times {\R}_{\ge 0}$.
Here $I$ is a finite set, all half-spaces have a common
plane ${\R}^{n-1}\times \{0\}$ and $x$ belongs to this
plane. $\Z$-affine structure on $B^{sm}=\sqcup_{i\in I} {\R}^{n-1}\times {\R}_{> 0}$ is the natural one.

\end{conj}

This conjecture gives a local model for a singular
$\Z$-affine structure at a singular component of codimension one.
Let us discuss the case of higher codimension.
We start with the following definition.

\begin{defn} A $\Z$-affine structure with singularities on $B$
is given by:\begin{enumerate}
\item a closed subset $B^{pre-sing}\subset B$ of a compact space $B$;

\item a $\Z$-affine structure on the open set
$B\setminus B^{pre-sing}$ .
\end{enumerate}

\end{defn}

One can think about closed set of ``potential singularities'' $B^{pre-sing}$ as containing the actual set of
singularities $B^{sing}$).

\begin{defn} A continuous path $\gamma(t), t\in [0,1]$ in the space of $\Z$-affine
structures with singularities on a given compact space $B$ is given by:
\begin{enumerate}
\item a continuous path $B_t^{pre-sing}$ 
in the space of all compact subsets of $B$,
\item a $\Z$-affine structure on $B\setminus B_t^{pre-sing}$ for all $t\in [0,1]$
\end{enumerate}

Notice that for each $t_0\in (0,1)$ and $x_0\in B\setminus B_{t_0}^{pre-sing}$
we can choose neighborhoods $U_{t_0}$ of $t_0$ and $U_{x_0}$
of $x_0$ such that $U_{x_0}\subset B\setminus B_t^{pre-sing}$ for all
$t\in U_{t_0}$. Then we require that:
\begin{enumerate}
\item[3.] if $U_{t_0}$ and $U_{x_0}$ are sufficiently small
then the induced $\Z$-affine structure on $U_{x_0}$
does not depend on $t\in U_{t_0}$.
\end{enumerate}

\end{defn}
Notice that in the case when the homotopy type of $B\setminus B_t^{pre-sing}$
 remains unchanged the representation $\rho_t:\pi_1\left( B\setminus B_t^{pre-sing}\right)\to
GL(n,\Z)\ltimes \R^n$ stays the same.

We are going to give an example of a non-trivial path in the next
subsection. We expect that  singularities which appear in 
 the collapse of Calabi-Yau manifolds satisfy the following

\begin{conj} If $B^{sing}=B^{pre-sing}$ is of codimension
at least two in $B$, then there is a continuous path  $\gamma(t)$ in the
space of $\Z$-affine structures with singularities which connects
a given structure with the one coming from a PL compactification, and such that for all 
 $t\in [0,1]$ we have $codim(B_t^{pre-sing})\ge 2$ and $\gamma(t)$ has Finiteness and Independence properties.

\end{conj}

\subsection{Standard singularities in codimension two}

Let us remove the angle $\left\{(x,y)\in\R^2\,|\,\,0<x<y\,\right\}$ from ${\R}^2$. 
After that we 
identify sides of the angle by the affine transformation
$(x,y)\mapsto (x+y,y)$. In this way we introduce a new
$\Z$-affine structure on ${\R}^2\setminus \{(0,0)\}$ with
the monodromy around $(0,0)$ given by the unipotent matrix (see Figure 3)
$$\left(\begin{array}{cc} 1 & 1\\0& 1
\end{array}\right)\,\,.$$

\begin{figure}\label{figure3}
\centerline{\epsfbox{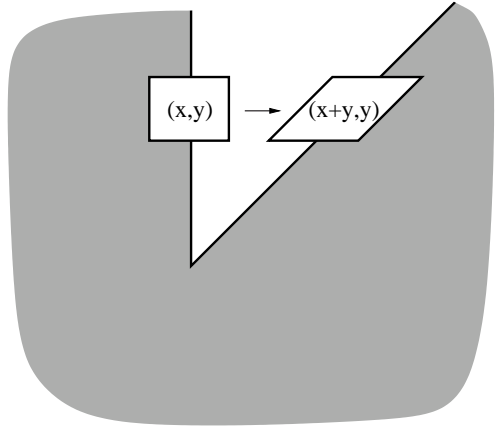}}
\caption{Glue both sides of the dashed area. White parallelograms are identified.}
\end{figure}

This $\Z$-affine structure does not admit a continuation
to ${\R}^2$. Therefore we obtain a $\Z$-affine structure
with singularities on ${\R}^2$. We will call 
{\it standard} the singularity at $(0,0)$.

Equivalently, we can describe this 
 $\Z$-affine structure on ${\R}^2\setminus \{(0,0)\}$
by taking a
cut along the ray $\{(x,0)\,|\,x> 0\}$ in ${\R}^2$ and 
glue the standard $\Z$-affine structure above and below the cut
by means of the affine transformation $(x,y)\mapsto (x+y,y)$
(see Figure 1 in Section 3.2.4). In this description it is clear that
 we can start the cut at \emph{arbitrary} point $(x_0,0)$ on the $x$-axes.
  The resulting singularity will be also called the standard one.

\begin{rmk}
We can vary a position of $(x_0,0)$, thus obtaining
a continuous path in the space of $\Z$-affine structures with
singularities in ${\R}^2$.

More generally, suppose that $B$ is equipped with a $\Z$-affine structure
which has standard singularities at points $p_1,\dots,p_m$.
Then we can slightly move each point $p_i$ in the direction invariant
under the local monodromy around $p_i$. This gives a new $\Z$-affine 
structure which is {\Z}PL-isomorphic to the initial one. 
\end{rmk}

Standard singularity is
called focus-focus singularity in the theory
of integrable systems (see [Zu]). In non-archimedean
geometry it appears as a singular value 
of some map $f:X^{an}\to {\R}^2$, where $X$ is an algebraic 
surface in  $3$-dimensional affine space
${\bf A}_K^3$ (see Section 8).

Let us consider the Cartesian product of 
${\R}^2\setminus \{(0,0)\}$ equipped with
the above $\Z$-affine structure with the standard (non-singular)
$\Z$-affine structure on ${\R}^{n-2}$. Let us choose
a continuous function $f(z_1,...,z_{n-2})$ and start the cuts
at all points $(f(z_1,...,z_{n-2}), 0,z_1,...,z_{n-2})$. 
This means that we introduce
the standard non-singular $\Z$-affine 
structure in the region $y\ne 0$ 
as well as
in the region $(y=0, x<f(z_1,...,z_{n-2}))$.
 Near points $(y=0, x>f(z_1,...,z_{n-2}))$
we introduce a modified $\Z$-affine structure
by declaring functions 
$$(y, x+\max(y,0), z_1,\dots,z_{n-2})$$
 to be $\Z$-affine coordinates.
This gives an example of a ``curved" singular set $B^{sing}$ of codimension 2.
Since  function $f$ can be approximated by PL functions,
the above $\Z$-affine structure can be deformed to
a PL one.

\subsection{$\Z$-affine version of Gauss-Bonnet theorem}

Let $B$ be a connected compact oriented topological surface,
$B^{sing}\subset B$ a finite set. Assume that 
$B^{sm}=B\setminus B^{sing}$ carries a $\Z$-affine structure
such that for any $x\in B^{sing}$ there exists a small neighborhood
$U$ such that $U=\cup_{i\in I}U_i$ where $I$ is a finite set
and each $U_i$ is affine equivalent to a germ of an angle in ${\R}^2$,
with $x$ being the apex of each angle.

The aim of this section is to define a map
$i_{loc}:B^{sing}\to {\frac{1}{12}}{\Z}$ (which depends only on $\Z$-affine
structure near $B^{sing}$) and prove the following result
(a kind of Gauss-Bonnet theorem).

\begin{thm} The following equality holds

$$\sum_{x\in B^{sing}}i_{loc}(x)=\chi(B)\,\,,$$
where $\chi(B)$ is the Euler characteristic of $B$.

\end{thm}

We start with the construction of $i_{loc}$. Let us denote by
$\widetilde{SL(2,{\Z})}$ the pre-image of $SL(2,{\Z})$ in the universal
covering $\widetilde{SL(2,{\R})}$ of the group $SL(2,{\R})$.
The group $\widetilde{SL(2,{\R})}$ contains 
$\pi_1(SL(2,{\R}))\simeq {\Z}$. Let $u$ be a generator of the
latter (it belongs also to $\widetilde{SL(2,{\Z})}$).

We have an exact sequence of groups

$$1\to {\Z}\to \widetilde{SL(2,{\Z})}\to PSL(2,{\Z})\to 1\,\,.$$
Notice that $PSL(2,{\Z})$ is a free product ${\Z}/2*{\Z}/3$.
Moreover, in the above exact sequence ${\Z}$ is embedded into
the center of $\widetilde{SL(2,{\Z})}$. Notice that $u$ is the image
of $2\in {\Z}$. One can choose representatives $a_2,a_3$ 
of the standard generators of $PSL(2,{\Z})$ in such a way that
$\widetilde{SL(2,{\Z})}$ is generated by $u,a_2,a_3$ subject
to the relations $a_2^2=a_3^3, a_2^4=a_3^6=u$.
This gives a homomorphism of groups 
$\phi:\widetilde{SL(2,{\Z})}\to {\Z}$ such that $\phi(a_2)=3,
\phi(a_3)=2, \phi(u)=12$. Dividing by $12$ we obtain a
homomorphism 
$i:\widetilde{SL(2,{\Z})}\to {\frac{1}{12}}{\Z}$ such that $i(u)=1$.

Let us consider a topological $S^1$-bundle $E$ over $B$, such that
the fiber over $x\in B$ is the union of all affine rays outcoming
of $x$. Then the restriction of $E$ to $B^{sm}$ is just
the spherical bundle. Let us pick $x_0\in B^{sm}$ and
remove it from $B$ together with small neighborhoods of
all points $B^{sing}$. We denote by $B_1$ the topological space
obtained in this way.
We can trivialize the tangent bundle over $B_1$ 
(we choose a $C^{\infty}$-trivialization, compatible with $SL(2,{\R})$-structure) in such a way that it
extends to a continuous trivialization
of $S^1$-bundle $E$ over $B\setminus\{x_0\}$. Let 
$\alpha\in \Omega^1(B_1)\otimes sl(2,{\R})$ be a $1$-form
defined by means of the affine structure $\nabla$ on $B_1$.
Then $d\alpha+{\frac{1}{ 2}}[\alpha,\alpha]=0$ and
$\alpha$ defines a flat connection on a
trivial $\widetilde{SL(2,{\R})}$-bundle on $B_1$.  It gives a
a monodromy representation 
$\pi_1(B_1)\to \widetilde{SL(2,{\Z})}$ defined up to
 a conjugation. Composing
it with the homomorphism $i$ we obtain a homomorphism
$\pi_1(B_1)\to {\frac{1}{12}}{\Z}$. Since ${\frac{1} {12}}{\Z}$
is an abelian group, the latter homomorphism is a composition
$\pi_1(B_1)\to H_1(B_1,\Z)\to {\frac{1}{12}}{\Z}$.
Let us pick small circles $[\gamma_x]\in H_1(B,\Z)$ for each
$x\in B^{sing}$. Then the above homomorphism gives us a number
denoted by $i_{loc}(x)\in {\frac{1}{12}}{\Z}$.

{\it Proof of the Theorem}. Let us pick up a small circle 
$[\gamma_{x_0}]\in H_1(B_1,\Z)$
around $x_0$. Then $\sum_{x\in B^{sing}}[\gamma_x]+[\gamma_{x_0}]=0$.
The monodromy around $x_0$ can be easily
computed via the winding number of the induced vector
field (section of $E_{|\gamma_{x_0}}$) and is equal to $-\chi(B)u$.
Applying homomorphism $i$ we obtain the result. $\blacksquare$

\begin{cor} Suppose that the monodromy for each point 
$x\in B^{sing}$ is the standard one (see Section 6.4).
Then one has two possibilities:

a) $B^{sing}=\emptyset$ and $B=B^{sm}$ is
a  $2$-dimensional torus;

b) the set $B^{sing}$ consists of $24$ distinct points
on the sphere $S^2$.

\end{cor}

{\it Proof.} It is easy to see that for each point $x\in B^{sing}$
one has $i_{loc}(x)={\frac{1}{12}}$. Then from Gauss-Bonnet
theorem one deduces that $\chi(B)=2-2g$, where $g$ is the genus
of the Riemann surface $B$. Then we have
$${\frac{|B^{sing}|} {12}}=2-2g\,\,.$$
Since LHS is non-negative we conclude that either $g=1$ or
$g=0$. In the first case $|B^{sing}|=0$ and we have a $\Z$-affine
structure on  a 
torus. In the second case we have $|B^{sing}|=24$ and $g=0$. $\blacksquare$

\begin{rmk} This corollary was proved in [LeS] by different
methods.

\end{rmk}

Similarly, for the affine structure with the monodromy
at each point conjugate to 
$$\left(\begin{array}{cc} 1 & 4\\0& 1
\end{array}\right)$$
one has $i_{loc}(x)=\frac{1}{3}$, and we have six singular points
on $S^2$ (see Section 4.2.5).

\subsection{Skeleton of a non-archimedean Calabi-Yau variety}

Let $X=X/K$ be a smooth proper algebraic variety over
a non-archimedean field $K$, 
$\dim\,X=n$ and $\Omega \in \Gamma(X,\Omega_X^n)$
be a non-zero top degree form on $X$.
We will associate canonically with the pair $(X,\Omega)$ 
a piecewise-linear
compact space $Sk(X,\Omega)\subset X^{an}$ such that
$\dim_{\R}Sk(X,\Omega)\le n$.

Let us assume for simplicity that $K={\C}((t))$ and $X$ is defined over
${\C}_{t}^{mer}\subset K$. 
Analytic space $X^{an}$ contains a dense subset $X_{Div}$ of
divisorial points corresponding to irreducible components
of special fibers of all snc models ${\cal X}$ of $X$ (see Appendix A):

$$ X_{Div}=\cup_{\cal X} i_{\cal X}(V_{S_{\cal X}})\,\,\,,
$$
where $V_{S_{\cal X}}$ is the set of vertices of the Clemens polytope
$S_{\cal X}$.

Top degree form $\Omega$ gives rise to a map
$\psi_{\Omega}:X_{Div}\to {\bf Q}$. Namely, if
$p:{\cal X}\to Spec({\C}_t^{mer})$ is a snc model and
$D\subset {\cal X}_0$ is an irreducible divisor of the special fiber then
we define

$$ \psi_{\Omega}(D)={{ord_D(\Omega\wedge dt/t)}\over{{ord_D(p^{\ast}(t))}}}\,\,.
$$
Here $\Omega\wedge dt/t$ is a meromorphic top degree form on
${\cal X}/{\C}$.

It is easy to show that $\psi_{\Omega}(D)$ depends only on the
point $i_{\cal X}(D)\in X^{an}$. Function $\psi_{\Omega}$ is (globally)
bounded from below. 

\begin{defn} A divisorial point $i_{\cal X}(D)$ is called essential
if $$\psi_{\Omega}(D)=\inf_{x\in X_{Div}}\psi_{\Omega}(x)\,\,.$$

\end{defn}

\begin{defn} Skeleton $Sk(X,\Omega)$ is the closure in $X^{an}$
of the set of essential points.\footnote{Our notion of a skeleton should not be mixed with the one
 introduced in [Be3]. The latter is related to the Clemens polytope $S_{\cal X}$ of a snc model $\cal X$.}

\end{defn}

Let ${\cal X}$ be a snc model. We will explain how to describe
$Sk(X,\Omega)$ in terms of ${\cal X}$ and $\Omega$. In fact
it is a nonempty simplicial subcomplex of $i_{\cal X}(S_{\cal X})$.

Let us call {\it ${\cal X}$-essential}
a divisor $D_i\subset {\cal X}_0$ such that
$$\psi_{\Omega}(D_i)=\min_{D_j\in {\cal X}_0}\psi_{\Omega}(D_j)\,\,.$$
A nonempty collection $D_{i_1},...,D_{i_l}$ of divisors in ${\cal X}_0$
is called {\it ${\cal X}$-essential} if all $D_{i_k}$ are 
${\cal X}$-essential, the intersection 
$D_{i_1}\cap D_{i_2}\cap...\cap D_{i_l}$ is non-empty and does not 
belong to the closure of the divisor of zeros of $\Omega$
is ${\cal X}\setminus {\cal X}_0$.

\begin{thm} The skeleton $Sk(X,\Omega)$ is the image under
$i_{\cal X}$ of the subcomplex $Sk({\cal X},\Omega)\subset S_{\cal X}$
consisting of simplices corresponding to ${\cal X}$-essential
collections of divisors.

\end{thm}

{\it Sketch of the proof.} Notice that for any snc model
${\cal X}$ the subset $S_{\cal X}({\bf Q})\subset S_{\cal X}$ 
consisting of points with rational barycentric coordinates is
mapped by $i_{\cal X}$ into $X_{Div}$.
Namely, we can modify ${\cal X}$ by blowing up at nonempty
intersections of irreducible components  of the special fiber
and then continue this process indefinitely. Divisorial points
obtained in this way exhaust all points of 
$i_{\cal X}(S_{\cal X}({\bf Q}))$.

We will prove that the set of essential points in $X^{an}$ coincides
with $i_{\cal X}(S_{\cal X}({\bf Q})\cap Sk({\cal X},\Omega))$.
First of all, a direct computation shows that $\psi_{\Omega}$
being restricted to $i_{\cal X}(S_{\cal X}({\bf Q}))$ achieves
its absolute minimum on 
$i_{\cal X}(S_{\cal X}({\bf Q})\cap Sk({\cal X},\Omega))$.
Secondly, another  straightforward computation shows
that the latter set  does not change under blow-ups
of first and second type (see Section A.5 in Appendix A). This concludes the proof.
$\blacksquare$

For Calabi-Yau manifold $X$ we will denote $Sk(X,\Omega)$ simply by $Sk(X)$,
 as there exists only one (up to a scalar) non-zero top-degree form $\Omega$ on $X$
 and $Sk(X,\lambda \Omega)=Sk(X,\Omega)\,\,\,\forall \lambda\in K^\times$.

One can prove that
the PL space $Sk(X,\Omega)$
is in fact a birational invariant. Moreover the group
$Aut^{brt}(X)$ of birational automorphisms of $X$ acts
on the skeleton by  ${\Z}PL$ transformations.
In order to obtain non-trivial examples of such actions
we need Calabi-Yau manifolds with
large groups of birational automorphisms. An 
example of a ${\Z}PL$-action is considered in the next subsection.

\subsection{K3 surfaces and ${\Z}PL$-actions on $S^2$}

\subsubsection{Integrable systems}

Recall that in Section 3.3 we constructed a $38$-dimensional space ${\cal P}$
parameterizing  integrable systems $(X,\omega)\to B$ with $B\simeq S^2$. 
The space ${\cal P}$ carries a codimension $20$ foliation ${\cal F}$ corresponding
to small deformations of integrable systems which do not change
the invariant $[\rho]$ of the local system 
$\rho: \pi_1(B^{sm})\to SL(2,{\Z})\ltimes {\R}^2$.
We explained that the fundamental group of a leaf of ${\cal F}$
acts by PL homeomorphisms of $S^2$.
Here we are going to give a (partial) description of ${\cal P}$ and ${\cal F}$
in cohomological terms using Torelli theorem (see Appendix B).

An algebraic polarized K3 surface $X/{\C}$ elliptically fibered 
over ${\C}P^1$, equipped with a holomorphic volume form
$\Omega$ can be encoded by the data
$(\Lambda, (\cdot, \cdot), [\omega], [\Omega], [\gamma], {\cal K}_X)$, where
\begin{enumerate}

\item $(\Lambda, (\cdot, \cdot), \C [\Omega], {\cal K}_X)$ is a K3 period data;

\item $[\omega],[\gamma]\in \Lambda$, $\Omega\in \Lambda\otimes {\C}\,\,$;

\item $[\omega]\in {\cal K}_X,\,\,\,\gamma\in \partial  {\cal K}_X,\,\,\, ([\omega], [\Omega])=([\gamma], [\Omega])=([\gamma], [\gamma])=0\,\,$;

\item $\gamma$ is a non-zero primitive lattice vector.
\end{enumerate}
Here $[\omega]$ is the class of polarization (projective
 embedding)  of $X$, $[\gamma]$ is dual to the class of generic
 fiber of the elliptic fibration $\pi: X\to {\C}P^1$. 

Perhaps one can express in cohomological terms the fact that
$\pi$ has exactly $24$ critical values. The latter is an open
condition.

Let $L\subset H_2(X, {\Z})$ be a subgroup consisting of homology classes
which can be represented by cycles which are projected into graphs in $B^{sm}$
(such cycles are circle fibrations over graphs).
When we move along a leaf of ${\cal F}$ then the pairing
of $Re([\Omega])$ with $L$ remains unchanged (see Section 3.1.1). Clearly
$L\subset [\gamma]^{\perp}$, and moreover, one can check that
$L=[\gamma]^{\perp}\simeq {\Z}^{21}$. The pairing
with $Re([\Omega])$ gives a map 
$\Lambda_{2,18}:= [\gamma]^{\perp}/{\Z}[\gamma]\to {\R}$, where
$\Lambda_{2, 18}$ is the following even 
unimodular lattice of signature $(2,18)$:
$$\Lambda_{2, 18}=\left(\begin{array}{cc} 0 & 1 \\ 1 & 0\end{array}\right)
\oplus
\left(\begin{array}{cc} 0 & 1 \\ 1 & 0\end{array}\right)
\oplus \left(-E_8\right)\oplus \left(-E_8\right)\,\,,$$
where $-E_8$ is the Cartan matrix for Dynkin diagram 
$E_8$ taken with the minus sign.

The functional $(Re[\Omega], \cdot)$ on $\Lambda_{2, 18}$ can be represented
as $(v_{Re[\Omega]},  \cdot)$ where 
$v_{Re[\Omega]}\in \Lambda_{2,18}\otimes {\R}$ is a vector with the strictly 
positive square norm. One can show 
 that the (non-Hausdorff) space of leaves of ${\cal F}$ is canonically  identified
with the set $\{v\in \Lambda_{2,18}\otimes {\R}|(v,v)>0\}/Aut (\Lambda_{2,18} )$.

The fundamental group of the leaf ${\cal F}_v$ corresponding 
to a vector $v\in \Lambda_{2,18}\otimes {\R}$ maps onto
the group $\Gamma_v\subset Aut (\Lambda_{2,18})$. This group is 
(up to a conjugation) the stabilizer in 
$(Aut(\Lambda_{2,18}), (\cdot, \cdot)_{2,18}, v)$
of the cone
$K_v$, which is a connected component of
the set 
$$\{w\in \Lambda_{2,18}\otimes {\R}|(w,v)=0, (w,w)>0\}\setminus
\bigcup_{\gamma\in \Lambda_{2,18}, (\gamma,\gamma)=-2, (\gamma,v)=0}H_{\gamma}$$
and $H_{\gamma}\in \Lambda_{2,18}\otimes {\R}$ is the hyperplane
orthogonal to $\gamma$ (cf. Appendix B).
Let us denote by $Aut_{{\Z}PL,v}(S^2)$ the group
of piecewise-linear transformations of $S^2$ with 
integer linear part. Index $v$ signifies the dependence of $\Z PL$-structure on $S^2$ on $v$.

\begin{conj} The homomorphism $\pi_1({\cal F}_v)\to Aut_{{\Z}PL,v}(S^2)$
arising from the monodromy of the local system along
the leaf ${\cal F}_v$ (see Sections 3.3, 6.4) is equal
to the composition

$$\pi_1({\cal F}_v)\twoheadrightarrow \Gamma_v\to Aut_{{\Z}PL,v}(S^2)\,\,,$$
where the homomorphism $\phi_v: \Gamma_v\to Aut_{{\Z}PL,v}(S^2)$ is
uniquely determined by this property.

\end{conj}

One can consider the whole moduli space ${\cal M}_{44}$ of $\Z$-affine structures on $S^2$ with $24$ standard singularities.
 This space is a Hausdorff orbifold (with a natural $\Z$-affine structure!)
 of dimension $44$, and it carries a foliation of codimension $20$ as before.
 It seems that using our main result (Theorem 5 in Part III) together with certain natural assumption 
(see Conjecture 11 in Section 11.6) one can show that
 the action by $\Z PL$ transformations of $S^2$ of the fundamental group of leaves of the foliation on the larger space
 ${\cal M}_{44}$ is again reduced to the action of $\Gamma_v$.

\subsubsection{Analytic surfaces}

Let $X=(X_t)_{t\to 0}$ be a maximally degenerate 
K3 surface over the field ${\C}_t^{mer}$ (see Section 5.1).
We denote by $\Lambda_X$ the quotient group
$[\gamma_0]^{\perp}/{\Z}[\gamma_0]$ where $[\gamma_0]\in H_2(X_t,\Z)$ is the vanishing
cycle. Then $\Lambda_X\simeq \Lambda_{2,18}$.
Let us assume  that the monodromy acts trivially on $\Lambda_X$.

We define a natural homomorphism $\rho_X:\Lambda_X\to ({\C}_t^{mer})^{\times}$ by the
formula

$$\rho_X([\gamma])=\exp\left(2\pi i {\int_{\gamma}\Omega_t\over{\int_{\gamma_0}\Omega_t}}\right),
\,\,\,[\gamma]\in [\gamma_0]^{\perp}\,\,.$$

One can give a more abstract definition of $\rho_X$ in terms
of the variation of Hodge structure. It is easy to see
that $\left(val_{{\C}_t^{mer}}\circ \rho_X\right)([\gamma])=(v_X, [\gamma])$ where
$v_X\in \Lambda_X$ is a vector such that $(v_X,v_X)>0$, and
$val_{{\C}_t^{mer}}$ is the standard valuation on the field
${\C}_t^{mer}\subset {\C}((t))$.

Let $X^{an} $ be the corresponding analytic K3 surface over
the field $K={\C}((t))$.
We have an analytic torus fibration 
over $S^2\setminus \{x_1,...,x_{24}\}$ which can be extended
to a continuous map $X^{an}\to S^2$. Let us call
such an extension a {\it singular} analytic torus fibration with standard singularities.

\begin{conj} For any analytic K3 surface $X^{an}/K$ admitting an analytic torus fibration
$X^{an}\to S^2$ with standard singularities, one can define intrinsically the lattice
$\Lambda_{X^{an}}$ and the homomorphism $\rho_{X^{an}}:\Lambda_{X^{an}}\to K^{\times}$.

\end{conj}

Notice that for K3 surfaces any birational automorphism
is biregular. Hence the
group of birational automorphisms $Aut^{brt}(X)$ acts
by a {\Z}PL-transformations of the sphere $S^2$ which is
equipped with a singular 
$\Z$-affine structure (see Section 6.6), i.e.
we have a homomorphism 
$$Aut^{brt}(X)=Aut(X)\to Aut_{{\Z}PL,v_X}(Sk(X^{an}, \Omega))\simeq
Aut_{{\Z}PL,v_X}(S^2)\,\,.$$

\begin{conj} 1) The image $\Gamma_{\rho_X}$ of $Aut(X)$ in $Aut(\Lambda_X, \rho_X)$
is a subgroup of $\Gamma_{v_X}$ where
$v_X:=val_K\circ \rho_X: \Lambda_X\to {\R}$.

2) The homomorphism $Aut(X)\to Aut_{{\Z}PL,v_X}(S^2)$ is conjugate to the 
restriction to $\Gamma_{\rho_X}$ of the homomorphism $\phi_{v_X}$
defined in the previous subsection.

\end{conj}

\subsubsection{Lattice points}

Let us consider the special case when vector $v$ is a lattice vector, i.e. $v\in \Lambda_{2,18}$.
 In A-model picture it corresponds to the integrality of the class $[\omega]$ of symplectic 2-form.
 In B-model this
 means that the non-archimedean field $K$ has valuation in $\Z\subset \R$.
 In terms of $\Z$-affine structures it means that the monodromy of the affine connection is reduced
 to $SL(2,\Z)\ltimes \Z^2$.
Group $\Gamma_v$ is a subgroup (and also a quotient group)
 of an arithmetic subgroup in the Lie group $SO(1,18)$. 
Also in this case there is a $\G_v$-invariant notion
 of a point with integer coordinates on $B\simeq S^2$, as well of 
points with coordinates in 
$\frac{1}{N}\Z$ for any integer $N\ge 1$.
 The number $M_{v,N}$ of such points is finite. It is not hard to see that
 $M_{v,N}=Area_v+2={{(v,v)}\over {2}} N^2+2$
 where $Area_v$ is the area of $B$ with a $\Z$-PL structure corresponding to $v$. 
This is analogous to the Riemann-Roch formula
 $rk\,\left(\Gamma(X_\C,L^{\otimes N})\right)=\int_X{{c_1(L)^2}\over{2}}+2$ for an ample line bundle $L$ on a complex K3
 surface $X_\C$.

 The action of $\Gamma_v$ on $S^2$ gives rise to a homomorphism
  $\Gamma_v\to S_{M_{v,n}}$ where $S_{M_{v,n}}$ is the symmetric group. Also the action
  gives a homomorphism from $\Gamma_v$ to the mapping class group
 $\pi_1({\cal{M}}_{0,M_{v,N}}^{unord})$, the fundamental group of the moduli space of genus zero complex curves
 with  $M_{v,N}$ unordered distinct marked points. The last group is closely related to the braid group.
 The conclusion is that we have constructed homomorphisms from arithmetic
 groups to a tower of braid groups.

One can  deduce from Torelli theorem an interpretation of $\Gamma_v$ as a quotient group
 of the fundamental group of a neighborhood $U$ of a cusp in 19-dimensional moduli space of polarized 
 complex algebraic
 K3-surfaces, where 
vector $v$ corresponds to the polarization. Therefore the homomorphism
  $\Gamma_v\to S_{M_{v,n}}$ gives a finite covering $U'$ of $U$. One may wonder 
whether there exists
 a line bundle over $U'$ whose direct image to $U$ coinsides with the direct image
 of the sheaf $L^{\otimes N}$ from the universal family of K3 surfaces
(this question is in spirit of some ideas of  Andrey Tyurin, see e.g. [Tyu]).

\subsection{Further examples}

There are many families of Calabi-Yau varieties with huge groups 
of birational automorphisms.
The following example we learned from D.Panov and D.Zvonkine.
 For any real numbers $l_1,\dots,l_n>0$ we can consider the space 
of planar $n$-gons with the length of edges equal to $l_1,\dots,l_n$,
 modulo the group of orientation-preserving motions.
 This space can be identified with the space of solutions
 of the following system of equations
  $$\sum l_i z_i=0,\,\,\,\sum l_i z_i^{-1}=0$$
where $(z_1:\dots :z_n)\in \C P^{n-1}$ 
is a point satisfying the reality condition
 $|z_i|=1, \,\,i=1,\dots,n $.
Hence we obtain a singular subvariety of $\C P^{n-1}$ 
of codimension $2$, depending on 
parameters $l_1,\dots, l_n$. 
One can check that this variety is birationally isomorphic
 to a non-singular Calabi-Yau variety.
 For any proper set $I\subset \{ 1,\dots,n\}$, $2\le |I|\le n-2$
 we have a birational involution
$\sigma_I$ defined by the formula
$$\sigma_I^*(z_i)=\left\{\begin{array}{ll}c/z_i &\mbox{ if $i \in I$}\\
 z_i &\mbox{ if $i \notin I$}\end{array}\right.$$
where $c:=\frac{\sum_{i \in  I} l_i z_i}{\sum_{i\in I} l_i/z_i}$.

We do not know at the moment the structure of the 
group $G_n$ generated by
involutions $\sigma_I$. One can obtain easily explicit formulas
 for the action of $G_n$ by piecewise-linear homemorphisms of 
 $S^{n-3}$.
Length parameters $l_i$ should be replaced by elements of a 
non-archimedean field $K$
with ``generic'' norms $\lambda_i=val_K(l_i)\in \R$. 
Denote by $\zeta_i,\,\,i=1,\dots,n$ real
 variables which have the meaning of valuations of variables $z_i\in K$.
 Sphere $S^{n-3}$ is obtained in the following way. In $\R^n$ we consider the intersection
of two subsets:
$$\{(\zeta_1,\dots\zeta_n)\,|\,\,\,\min_i(\lambda_i+\zeta_i)\mbox
{ is achieved at least twice}\}$$
and
$$\{(\zeta_1,\dots\zeta_n)\,|\,\,\,\min_i(\lambda_i-\zeta_i)\mbox
{ is achieved at least twice}\}$$
and then take the quotient by the action of $\R$: 
$$(\zeta_1,\dots\zeta_n)\ra (\zeta_1+c,\dots\zeta_n+c)$$
corresponding to the projectivization. For appropriately chosen
$(\lambda_1,\dots,\lambda_n)$ we obtain a set which is the union of $S^{n-3}$
 with several ``wings" going to infinity.
 The action of the involution $\sigma_I$ is obtained from algebraic formulas from above,
in which one replace non-archimedean variables by real ones,
 addition by minimum and multiplication (division) by addition (subtraction).

\section{$K$-affine structures}

\subsection{Definitions}

Let $B^{sm}$ be a manifold with $\Z$-affine
structure. The sheaf of $\Z$-affine functions 
$Aff_{\Z}:=Aff_{{\Z},B^{sm}}$ gives rise to an exact
sequence of sheaves of abelian groups

$$ 0\to {\R}\to Aff_{\Z}\to (T^{\ast})^{\Z}\to 0\,\,.$$

Let $K$ be a complete non-archimedean  field with a valuation
map $val$. We give two equivalent definitions of a $K$-affine structure on $B^{sm}$ compatible with a 
given $\Z$-affine structure.

\begin{defn} A  $K$-affine structure on $B^{sm}$
compatible with the given $\Z$-affine structure
is a sheaf $Aff_K$ of abelian groups
on $B^{sm}$, an exact sequence of sheaves

$$0\to K^{\times}\to Aff_K\to (T^{\ast})^{\Z}\to 0\,\,,$$
together with a homomorphism $\Phi$ of this exact sequence
to the exact sequence of sheaves of abelian groups
$$0\to {\R}\to Aff_{\Z}\to (T^{\ast})^{\Z}\to 0\,\,,$$
such that  $\Phi=id$ on $(T^{\ast})^{\Z}$ and
$\Phi=val$ on $K^{\times}$.

\end{defn}

Since $B^{sm}$ carries a $\Z$-affine structure, we
have an associated  $GL(n,{\Z})\ltimes {\R}^n$-torsor on $B^{sm}$, 
whose fiber
over a point $x$ consists of all $\Z$-affine coordinate systems
at $x$.

\begin{defn} A $K$-affine structure on $B^{sm}$
compatible with the given $\Z$-affine structure is
a $GL(n,{\Z})\ltimes (K^{\times})^n$-torsor 
on $B^{sm}$ such that the application of $val^{\times n}$ to
$(K^{\times})^n$ gives the initial $GL(n,{\Z})\ltimes {\R}^n$-torsor.

\end{defn}

Equivalence of two definitions from above is obvious in local $\Z$-affine coordinates.
 The reason is that the 
 set
of automorphisms of the exact sequence of groups

$$0\to K^{\times}\to K^{\times}\times {\Z}^n\to {\Z}^n\to 0$$
identical on $K^\times$ coincides with 
the group $GL(n,{\Z})\ltimes (K^{\times})^n$.

Finally, we can formulate  the Fixed Point Property
for $K$-affine structures (see  Section 3.1 for
$\Z$-affine case):
\newline{\bf Fixed Point Property for $K$-affine structures}.
In the notation of the end of Section 3.1,
for any $b\in B^{sing}$ and sufficiently
small neighborhood $U$ of $b$ the lifted monodromy representation
$\pi_1(U)\to GL(n,{\Z})\ltimes (K^{\times})^n$ has fixed vectors
in $K^{\times n}$, and the ${\R}$-affine span of the corresponding
(under the valuation map) vectors in ${\R}^n$ coincides with the
set of fixed points of the monodromy representation
$\pi_1(U)\to GL(n,{\Z})\ltimes {\R}^n$.

\subsection{$K$-affine structure on smooth points}

Starting from this section till the end of the paper (except of the Section 11.7) we will assume the following
\newline {\bf Zero Charactersistic Assumption}. 
 $K$ is a complete non-archimedean local field such that
its residue field has characteristic zero.

Let $X$ be a 
$K$-analytic manifold of dimension $n$ and   we are given a continuous
map $\pi: X\to B$, where $B$ is a topological space.
 Then $B^{sm}$ carries a $\Z$-affine
structure (Theorem 1). Suppose that there is an open $K$-analytic
submanifold $U\subset X$ such that 
$\pi^{-1}(B^{sm})\subset U$ and there is a nonwhere vanishing
analytic form $\Omega\in \Gamma(U,\Omega_X^n)$.
We are going to define a $\Z$-affine function $Val(\Omega)$
similarly to the definition of the function $Val(\varphi)$
in Section 4.1. 
Namely, in local coordinates $(z_1,...,z_n)$ we consider the expression
$\varphi:=\Omega/\bigwedge_{1\le i\le n}(dz_i/ z_i)$. This is an invertible function, and we define
$Val(\Omega)$ as   $Val(\varphi)$. 
The independence on the choice of coordinates follows from the following lemma
\begin{lmm} Let $(z_i)_{i=1,\dots,n},\,\,(z'_i)_{i=1,\dots,n}$ be two systems of invertible coordinates
 on $\pi^{-1}(U)$ for some connected open $U\subset B^{sm}$. Then 
$$\left|\left({\textstyle \bigwedge_{1\le i\le n} (dz_i/z_i)}\right)/
\left({\textstyle \bigwedge_{1\le i \le n }
(dz'_i/ z'_i)}\right)\right|_x=1
\,\,\,\forall x\in\pi^{-1}(U)\,\,.$$ 
\end{lmm}
{\it Proof:} By Lemma 1 from Section 4.1 we know that $z_i'$ as any invertible function can be written
 in form $c_i z^{I^{(i)}} (1+o(1))$ for some nonzero $c_i
\in K$ and a multi-index $I^{(i)}\in\Z^n$. Vectors $I^{(1)},\dots, I^{(n)}$
form a basis of $\Z^n$, as follows from the condition that $z_1',\dots,z_n'$ form a coordinate system.
 Therefore, after applying the change of coordinates $z_i\mapsto c_i z^{I^{(i)}}$
  preserving form $\bigwedge_i dz_i/z_i$
 up to sign,
 we may assume that $z_i'=(1+o(1)) z_i$. The Jacobian matrix of the transformation $(z_i)\to (z_i')$
 is the identity matrix plus terms of size $o(1)$. 
Therefore its determinant has norm equal to 1. $\blacksquare$

Now we make the following
\newline{\bf Constant Norm Assumption}. The function $Val(\varphi)$
is locally constant.

\begin{thm} If the Constant Norm Assumption is satisfied then
there is a $K$-affine structure on $B^{sm}$ compatible
with the  $\Z$-affine structure $Aff_{{\Z},B^{sm}}^{can}$
(see Section 4.1).

\end{thm}

{\it Proof.} 
Let us write in local coordinates
 $\Omega=\varphi (z_1,...z_n)\bigwedge_{1\le i\le n}{dz_i\over z_i}$.
Define residue $Res(\Omega)\in K$ as the constant term $\varphi_0$
 in the Laurent expansion
 $\varphi (z_1,...z_n)=\sum_{I\in\Z^n} \varphi_I z^I$.
It is easy to see  that $Res(\Omega)$ does not depend (up to a sign) on the choice of local coordinates.
 For non-vanishing everywhere $\Omega$ satisfying Constant Norm Assumption we have 
 $\exp(-Val(\varphi))=|\varphi|=|\varphi_0|$.
   Therefore we have $Res(\Omega)\ne 0$.

Let us return to the proof of the Theorem. 
Let $F$ be the sheaf
of abelian groups $F\subset \pi_{\ast}({\cal O}_{X}^{\times})$
consisting of $f$  
such that $Val(f)=0$.
Then we have an exact sequence of sheaves

$$0\to K^{\times}/\O_K^{\times}\to 
\pi_{\ast}({\cal O}_{X}^{\times})/F\to (T_X^{\ast})^{\Z}\to 0\,\,\,,$$
where $\O_K$ denotes the constant sheaf with the fiber being the
ring of integers of $K$.
Indeed we embed $K^{\times}/\O_K^{\times}$ into 
$\pi_{\ast}({\cal O}_{X}^{\times})/F$ as constant functions.
The projection $\pi_{\ast}({\cal O}_{X}^{\times})/F\to (T_X^{\ast})^{\Z}$
assigns to the function $f$ the linear part of the corresponding
$\Z$-affine function $Val(f)$.

Notice that if $U\subset B^{sm}$ is a  connected domain
then any $f\in \Gamma(U,F)$ can be written (non-canonically) as
$f=a(1+r)$, where $a\in \O_K^{\times}$ and $r=o(1)$ in $\pi^{-1}(U)$.

We define an epimorphism of sheaves $p_{\Omega}:F\epi \O_K^{\times}$
by formula
$$p_{\Omega}(f)=p_{\Omega}(a(1+r))=a\,\exp\left({Res(\Omega\,\log(1+r))\over 
{Res(\Omega)}}\right)\,\,.$$

Here $\exp$ and $\log$ are understood as infinite 
convergent series (in order to make sense of them we use
Zero Characteristic Assumption).

It is easy to see that $p_{\Omega}$ is well-defined.
Then the exact sequence of sheaves

$$
1\to K^{\times}\to 
\pi_{\ast}({\cal O}_{X}^{\times})/\ker(p_{\Omega})\to (T_X^{\ast})^{\Z}\to 1
$$
defines a $K$-affine structure on $B^{sm}$ compatible
with  $Aff_{{\Z},B^{sm}}^{can}$.
This concludes proof of the Theorem. $\blacksquare$

Notice that the above proof gives an explicit construction
of the $K$-affine structure. We will denote it
by $Aff_{K,B^{sm}}^{\Omega}$. It is easy to see that this $K$-affine
structure does not change if we make a rescaling
$\Omega\mapsto c\Omega, c\in K^{\times}$.

\subsection{Lifting Problem}

Let $K$ be as in Section 7.2, $B\supset
B^{pre-sing}$ be a space with singular $\Z$-affine structure 
(see Section 6.3), and an extension
 of $\Z$-affine structure on $B\setminus B^{pre-sing}$ 
to a $K$-affine structure satisfying fixed point property (see 7.1).
 We assume that $\Z$-affine structure cannot be extended to a larger open
 set $U\supset B\setminus B^{pre-sing} , U\ne B\setminus B^{pre-sing}$.  
Slightly abusing notation we will denote $B\setminus B^{pre-sing}$ simply
by $B^{sm}$.
We want to have a  $K$-analytic space $X$, meromorphic
non-zero top degree form $\Omega$ and a continuous proper (and maybe also Stein) map
$\pi:X\to B$ such that:
\begin{enumerate}
\item $B^{pre-sing}$ coincides with $B^{sing}$, and $\Z$-affine structure on $B^{sm}$ arising from the projection $\pi$
coincides with 
the given one;

\item the restriction $\Omega_{|\pi^{-1}(B^{sm})}$ is
a nowhere vanishing analytic form which satisfies 
the Constant Norm
Assumption;

\item the $K$-affine structure on $B^{sm}$ arising from
the pair $(X,\Omega)$ coincides with the initial one.

\end{enumerate}
We call the problem of finding such data {\it Lifting Problem}.

\begin{rmk} If a solution of the Lifting Problem exists  then
$B^{sm}$ is orientable. Indeed, $Res(\Omega)$ is locally
a constant defined up to a sign which depends on the orientation
of $B^{sm}$. Global choice of the constant gives an
orientation. For oriented $B^{sm}$ we can rescale
$\Omega$ canonically in such a way that $Res(\Omega)=1$

\end{rmk}
{\bf Question}. What restrictions on the behavior of the $K$-affine
structure near $B^{pre-sing}=B^{sing}$ 
should we impose in order to guarantee the existence
of a solution of the Lifting
Problem?

Let $B=B^{sm}$ be a flat torus (see Section 3.2.1).
Then the Lifting Problem has a solution (canonical
up to rescaling of $\Omega$) for any
compatible $K$-affine structure. More precisely, the
groupoid of Tate tori and isomorphisms between them
is equivalent to the groupoid of $K$-affine structures
on real flat tori.

In Sections 8-11 we are going to discuss a solution
of the Lifting Problem 
for K3 surfaces. In that case $B^{sing}\ne \emptyset$.

If we restrict ourselves only to the smooth part  $B^{sm}$ 
(i.e. we allow non-compact $X$) then there
is a canonical solution of this ``reduced" Lifting
Problem. In other words one can construct a smooth
$K$-analytic space $X^{\prime}$ with an analytic top degree
form $\Omega^{\prime}$ and a map 
$\pi^{\prime}:X^{\prime}\to B^{sm}$ satisfying the 
above conditions 1--3.
Let us explain this construction assuming that $B^{sm}$
is oriented.

First of all we notice that the orientation of $B^{sm}$ gives
a reduction to 
$SL(n,{\Z})\ltimes (K^{\times})^n$ of the structure group of the torsor
defining the $K$-affine structure. The reduced group
naturally acts by automorphisms of the fibration
$\pi_{can}:({\bf G}_m^{an})^n\to {\R}^n$ preserving the form
$\bigwedge_{1\le i\le n}{dz_i\over z_i}$. The action on $({\bf G}_m^{an})^n$
is induced from the action on monomials. Namely, the inverse to an element 
$(A,\lambda_1,...,\lambda_n)\in SL(n,{\Z})\ltimes (K^{\times})^n$
acts on monomials as
$$z^{I}=z_1^{I_1}\dots z_n^{I_n}\mapsto 
\left({\textstyle\prod_{i=1}^n}\lambda_i^{I_i}\right)\,\,z^{A(I)}\,\,\,.$$ 
The action
of the same element on ${\R}^n$ is given by the similar formula 
$$x=(x_1,\dots ,x_n)\mapsto A(x)-(val(\lambda_1),\dots,val(\lambda_n))\,\,\,.$$

Let $B^{sm}=\cup_{\alpha}U_{\alpha}$ be an open covering
by coordinate charts $U_{\alpha}\simeq V_{\alpha}\subset {\R}^n$
such that for any $\alpha, \beta$ we are given elements
$g_{\alpha,\beta}\in SL(n,{\Z})\ltimes (K^{\times})^n$ satisfying the
$1$-cocycle condition for any triple
$\alpha,\beta,\gamma$. Then the space $X^{\prime}$ is obtained
from $\pi_{can}^{-1}(V_{\alpha})$ by gluing by means of the transformations
$g_{\alpha,\beta}$. The form $\bigwedge_{1\le i\le n}{dz_i\over z_i}$
gives rise to a nowhere vanishing analytic top degree form $\Omega^{\prime}$
on $X^{\prime}$.
Thus we have obtained a solution of the reduced Lifting Problem.
The sheaf $\pi_{\ast}({\cal O}_{X^{\prime}}):={\cal O}_{B^{sm}}^{can}$
is called the {\it canonical sheaf}.

In the case $B^{pre-sing}\ne \emptyset$ this solution seems to be
a ``wrong" one, i.e. it cannot be extended to a solution
$\pi:X\to B$, where $X$ and $B$ are compact. In the case of K3 surfaces
we will show later how to modify it in order to obtain a ``true" solution
of the Lifting Problem.

\subsection{Flat coordinates and periods}

Here we are going to discuss a relation between $K$-affine structures
and so-called flat coordinates on the moduli space
of complex structures on Calabi-Yau manifolds.
We assume the picture of collapse from Section 5.1.

\subsubsection{Flat coordinates for degenerating complex Calabi-Yau manifolds}
Let ${X}_{mer}=(X_t)_{t\to 0}$ be a maximally degenerating
 algebraic Calabi-Yau manifold of dimension $n$ over ${\C}_t^{mer}$. 
 We denote by $B$ the Gromov-Hausdorff limit
of our family (see Conjecture 1, Section 5.1). Its connected
oriented open dense part $B^{sm}$ carries a $\Z$-affine structure
with the covariant lattice $T^{\Z}$. 

Recall that according to the picture of collapse presented in
Section 5.1 there is a canonical isotopy class
of embeddings from a torus bundle $p:X_t^{\prime}\to B^{sm}$
to the complex manifold $X_t$ for all sufficiently small $t\ne 0$.
Let us denote by $[\gamma_0]\in H_n(X_t^{\prime}, {\Z})$ the
fundamental class of the fiber of $p$.
This is the homology class of a singular chain in $X_t^{\prime}$
which projects to a point by $p$. 

Let $H_n^{\le 1}(X_t^{\prime}, {\Z})\subset H_n(X_t^{\prime}, {\Z})$ 
be 
the subgroup generated by homology classes of chains which are
projected into graphs in $B^{sm}$.
It follows from the definition that we have an epimorphism
$$
J_a: H_1(B^{sm}, {\textstyle \bigwedge ^{n-1}} T^{\Z})\epi H_n^{\le 1}(X_t^{\prime}, 
{\Z})/{\Z}[\gamma_0]
$$
similar to the homomorphims $J_s$ defined in the symplectic case (see Section 3.1.1). 
The following formula defines a homomorphism of groups
$$P:H_n^{\le 1}(X_t^{\prime}, {\Z})/{\Z}[\gamma_0]\to ({\C}_t^{mer})^{\times},\,\,\,\,
[\gamma]\mapsto \exp\left(2\pi i{\int_{[\gamma]}\Omega_t\over \int_{[\gamma_0]}\Omega_t}\right)\,\,.$$

We will call $P$ the {\it period map}.
Notice that 
${\Z}[\gamma_0]:=H_n^{\le 0}(X_t^{\prime}, {\Z})\subset 
H_n^{\le 1}(X_t^{\prime},{\Z})$ is a low degree part of the limiting
Hodge filtration on the homology of Calabi-Yau manifold $X_t$.
Non-zero complex numbers

$$\exp\left(2\pi i{\int_{[\gamma_i]}\Omega_t\over \int_{[\gamma_0]}\Omega_t}\right)\,\,,$$
where $\gamma_i$ is a set of generators of 
$H_n^{\le 1}(X_t^{\prime},{\Z})/H_n^{\le 0}(X_t^{\prime}, {\Z})$
are called {\it flat coordinates} in Mirror Symmetry (see e.g. [Mor]).
Those are local coordinates near a point close to the ``cusp" of
the moduli space of complex structures (local Torelli theorem).

The orientation of $B^{sm}$ gives rise to an isomorphism
$\bigwedge^{n-1}T^{\Z}\simeq (T^{\ast})^{\Z}$. Therefore, combining maps $J_a,P$ and the
 above isomorphism
   we obtain
a homomorphism $$\widetilde{P}:H^1(B^{sm},(T^{\ast})^{\Z})\to \left({\C}_t^{mer}\right)^{\times}\,\,.$$

\subsubsection{Non-archimedean periods}
Let ${X}^{an}$ be a smooth analytic Calabi-Yau manifold
associated with $X_{mer}$. Assuming the equivalence
of Gromov-Hausdorff and non-archimedean pictures of collapse
presented in Section 5 we have a continuous map 
$\pi:{X}^{an}\to B$. It gives a $K$-affine structure on
$B^{sm}$. The corresponding exact sequence

$$0\to K^{\times}\to Aff_K\to (T^{\ast})^{\Z}\to 0$$
represents a class in 
$H^1(B^{sm},T^{\Z}\otimes K^{\times})\simeq 
Ext^1((T^{\ast})^{\Z},K^{\times})$. Pairing with this class gives
another homomorphism
$$P':H_1(B^{sm},(T^{\ast})^{\Z})\to K^{\times}=H_0(B^{sm},K^{\times})\,\,.$$

\begin{conj} Homomorphism $P'$ is equal to the composition of 
 $\widetilde{P}$ with the embedding $\left({\C}_t^{mer}\right)^\times\mono K^\times$.

\end{conj}

\part{}

We fix field $K$   satisfying Zero Characteristic  Assumption.

 Let $B$ be a compact oriented surface, $B^{sing}\subset B$
a finite set, and $Aff_K=Aff_{K,Y}$ a sheaf defining a $K$-affine
structure on $Y:=B^{sm}=B\setminus B^{sing}$.
We assume that all singularities of the underlying ${\Z}$-affine
structure are standard (see Section 6.4), and local monodromy
around each $b\in B^{sing}$ acts on $(K^{\times })^2$ with a fixed
point (see the Fixed Point property at the end of Section 7.1).
Main result of Part 3 of the paper can be formulated such as follows.

\begin{thm} There exist a compact $K$-analytic surface $X^{an}$,
a top degree analytic form $\Omega=\Omega_{X^{an}}$ and a continuous
proper Stein map $\pi: X^{an}\to B$ such that the set of $\pi$-smooth points
coincides with $Y$ and the induced $K$-affine structure coincides
with the one given by $Aff_K$.

In other words, the triple $(X^{an}, \pi, \Omega)$ is a solution
of the Lifting Problem.

\end{thm}

 By Stein property it suffices to construct the 
sheaf ${\cal O}_B=\pi_{\ast}({\cal O}_{X^{an}})$ of $K$-algebras on $B$.
We will see that outside of the finite singular set $S=\{x_1,\dots,x_{24}\}$
the sheaf ${\cal O}_B$
is locally isomorphic to ${\cal O}_Y^{can}$.
In the next section we will describe the local model for
the sheaf ${\cal O}_B$ near each singular point. It will be glued
together with a modification of the canonical sheaf
${\cal O}_Y^{can}$. This modification depends on the data
called {\it lines}. Appearance of lines is motivated by
Homological Mirror Symmetry (see [Ko], [KoSo])
\footnote{The main idea is that $X$ is a component of the moduli space
 of certain objects (skyscrapper sheaves) in the derived category $D^b(Coh(X))$.
 These objects correspond to $U(1)$-local systems on Lagrangian tori in the 
Fukaya category of the mirror dual symplectic manifold.}.
Roughly speaking, lines correspond (for mirror dual
K3 surface) to collapsing holomorphic discs with boundaries
belonging to fibers of the dual torus fibration (see Section 5.1 and [KoSo]).
Such ``bad" fibers are Lagrangian tori, but they do not
correspond to objects of the Fukaya category (A-branes
in terminology of physicists). There are infinitely many such
fibers and hence infinitely many lines. We will axiomatize
this piece of data in Section 9. Subsequently, with each line $l$
we will associate an automorphism of the restriction
of ${\cal O}_Y^{can}$ to $l$. This will give us the above-mentioned
modified canonical sheaf.

\section{Model near a singular point}

Here we will construct an analytic torus fibration corresponding to
 standard singularity (see Sections 3.2.4 and 6.4).

Let
$X\subset {\bf A}^3$ be the algebraic surface given by equation 
$ (\alpha\beta-1)\gamma=1$ in coordinates $(\alpha,\beta,\gamma)$,
and $X^{an}$ be the corresponding analytic space.
We define a continuous map 
$f: X^{an}\to {\R}^3$ by the formula
$f(\alpha,\beta,\gamma)=(a,b,c)$ 
where $a=\max(0,\log|\alpha|_p), b=\max (0, \log|\beta|_p), 
c= \log|\gamma|_p=-\log|\alpha\beta-1|_p$. Here 
$|\cdot|_p=\exp(-val_p(\cdot))$ denotes the multiplicative seminorm
corresponding to the point $p\in X^{an}$
(see Appendix A).

\begin{prp} The map $f$ is proper. Moreover

a) Image of $f$ is homeomorphic to $\R^2$.

b) All points of the image except of $(0,0,0)$ are $f$-smooth.

\end{prp}

{\it Proof.} Here is the plan of the proof.
\begin{enumerate}
\item We define three open domains $T_i,\,\,i=1,2,3$ in three copies of the standard
 two-dimensional analytic torus $({\bf G}_m^{an})^2$, 
and  continuous maps $\pi_i:T_i\to {\R}^2$ such that 
all points of the image $U_i=\pi_i(T_i)$ are $\pi_i$-smooth
(i.e. each $\pi_i$ is an analytic torus fibration).
 Domains $U_i$ cover $\R^2\setminus\{(0,0)\}$.

\item For each $i, 1\le i\le 3$ we  construct
an open embedding  $g_i:T_i\hookrightarrow X^{an}$.

\item We  construct an embedding $j: {\R}^2\hookrightarrow {\R}^3$ 
such that  each open set $U_i$ is homeomorphically
identified with $f(g_i(T_i))$ and $j((0,0))=(0,0,0)$.
 Moreover, $\pi_i$-smooth points are
mapped into $f$-smooth points.
\end{enumerate}

The Proposition will follow from 1)-3).

Let us describe the constructions and formulas.
We start with open sets $U_i, 1\le i\le 3$.
Let us fix a number $0<\varepsilon<1$ and define
$$\begin{array} {lll}
U_1 & = & \{(x,y)\in {\R}^2|x<\varepsilon |y|\,\}\\ 
U_2 & = & \{(x,y)\in {\R}^2|x>0, y<\varepsilon x\,\}\\ 
U_3 & = & \{(x,y)\in {\R}^2|x>0, y>0\} 
\end{array}$$
Clearly ${\R}^2\setminus\{(0,0)\}=U_1\cup U_2\cup U_3$. 
We define also a slightly modified domain $U_2'$ as 
$\{(x,y)\in {\R}^2|x>0, y<\frac{\varepsilon}{1+\varepsilon} x\,\}$.

We define $T_i:=\pi_{can}^{-1}(U_i)\subset ({\bf G}_m^{an})^2, i=1,3$ and
 $T_2:=\pi_{can}^{-1}(U_2')\subset ({\bf G}_m^{an})^2$.
Then the projections $\pi_i: T_l\to U_l$ are given by the formulas
$$
\pi_i(\xi_i,\eta_i)=\pi_{can}(\xi_i,\eta_i)=
(\log|\xi_i|, \log|\eta_i|), \,\,\,i=1,3\,\,,$$ 
$$\pi_2(\xi_2,\eta_2)=\left\{\begin{array}{ll}(\log|\xi_2|, \log|\eta_2|)& \mbox{  if }|\eta_2|<1\\
(\log|\xi_2|-\log|\eta_2|, \log|\eta_2|) & \mbox{ if }
|\eta_2|\ge 1
\end{array}\right.\,\,. $$
 In these formulas $(\xi_i,\eta_i)$ are  coordinates
on $T_i, 1\le i\le 3$.

We define inclusion $g_i: T_i\mono X, 1\le i\le 3$ by the following formulas:

$$\begin{array}{lll}
g_1(\xi_1,\eta_1) & = & ({1\over {\xi_1}}, \xi_1(1+\eta_1), {1\over {\eta_1}})\\
g_2(\xi_2,\eta_2) & = & ({1+\eta_2\over {\xi_2}}, \xi_2, {1\over {\eta_2}})\\
g_3(\xi_3,\eta_3) & = & 
({1+\eta_3\over {\xi_3\eta_3}}, \xi_3\eta_3, {1\over {\eta_3}})
\end{array}$$

Let us decompose $X^{an}=X_-\cup X_0\cup X_+$
 according to the sign of $\log |\gamma|_p$ where $p\in X^{an}$
 is a point. It is easy to see that 
$$\begin{array}{lll}
f(X_-) & = & \{\,(a,b,c)\in {\R}^3\,|\,c<0,a\ge 0, b\ge 0, \,ab(a+b+c)=0\,\}\\
f(X_0) & = & \{\,(a,b,c)\in {\R}^3\,|\,c=0,a\ge 0, b\ge 0, \,ab=0\, \}  \\
f(X_+) & = &  \{\,(a,b,c)\in {\R}^3\,|\,c>0,a\ge 0, b\ge 0,\, ab=0\, \}
\end{array}$$

From this explicit description we see that $f$ is proper
and the image of $f$ is
homeomorphic to ${\R}^2$. 

Let us consider the 
embedding $j: {\R}^2\to {\R}^3$ given by  formula 
$$j(x,y)=\left\{\begin{array}{lll}
(-x\,,\, \max(x+y,0)\,,\, -y\,) & 
\mbox{ if } & x\le 0\\
 (\,0\,,\, x+\max(y,0)\,,\, -y\,) & \mbox{ if } & x\ge 0
\end{array}\right.$$
One can easily check that the image of $j$ coincides with
the image of $f$, $j\circ \pi_i=f\circ g_i$ and 
$f^{-1}(j(U_i))=g_i(T_i)$ for all $1\le i\le 3$.
This concludes the proof of Proposition. $\blacksquare$

We can derive more from explicit formulas given in the proof.

Let us denote by $\pi: X^{an}\to {\R}^2$ the map $j^{(-1)}\circ f$.
It is an analytic torus fibration outside of point $(0,0)$.
The induced $\Z$-affine
structure on ${\R}^2\setminus \{(0,0)\}$ is in fact
the standard singular $\Z$-affine structure described in Sections 3.2.4 and 6.4,
 as follows immediately from formulas for projections $\pi_i,\,i=1,2,3$.

Let us introduce another sheaf ${\cal O}^{can}$ on ${\R}^2\setminus \{(0,0)\}$.
 It is defined as $(\pi_i)_*\left(\O_{T_i}\right)$ in each domain $U_i$,
 with identifications 
$$\begin{array}{llcl}
(\xi_1,\eta_1) & = & (\xi_2,\eta_2) & \mbox{ on } U_1\cap U_2 \\
(\xi_1,\eta_1) & = & (\xi_3,\eta_3) & \mbox{ on } U_1\cap U_3 \\
(\xi_2,\eta_2) & = & (\xi_3\eta_3,\eta_3) & \mbox{ on } U_2\cap U_3
\end{array}$$

Let us consider the direct image sheaf $\pi_{\ast}({\cal O}_{X^{an}})$.
It is easy to see that on the sets $U_1$ and $U_2\cup U_3$ this
 sheaf is canonically isomorphic to  ${\cal O}^{can}$.
 The isomorphism is given by the identification  
  of coordinates $(\xi_1,\eta_1)$ on $U_1$, and
 of coordinates $(\xi_2,\eta_2)$ and $(\xi_3,\eta_3)$ on  $U_2\cup U_3$.
 Therefore on the intersection $U_1\cap (U_2\cup U_3)$ we identify two 
copies of  the canonical sheaf by certain automorphism $\varphi$ of 
${\cal O}^{can}$
which preserves
one coordinate (namely, the  coordinate $\eta$). 
We will develop the theory of such transformations
and their analytic continuations in Section 11. The explicit formulas
 for $\varphi$ is
 $$\varphi(\xi,\eta)=\left\{\begin{array}{cll}
(\xi(1+\eta),\eta) & \mbox{ on } & U_1\cap U_2 \\
(\xi(1+1/\eta),\eta) & \mbox{ on } & U_1\cap U_3 
\end{array}
\right.$$

We would like to say now few words about analytic volume forms.
Notice that each $T_i\subset ({\bf G}_m^{an})^2$  carries a
nowhere vanishing  top degree
analytic form given by the formula 
$\Omega_i={d\xi_i\wedge d\eta_i\over{\xi_i\eta_i}}$.
Then a  straightforward calculation shows that 
 $\Omega^{T_i}$ is the pullback under $g_i$
 of nowhere vanishing on $X^{an}$
analytic top degree form

$$\Omega=-\gamma\,d\alpha\wedge d\beta\,\,.$$
Form $\Omega$ satisfies Constant Norm Assumption, hence it gives
 a $K$-affine structure on $\R^2\setminus \{(0,0)\}$.
On the other hand, the sheaf ${\cal O}^{can}$ of algebras is also
endowed with top-degree form $\Omega^{can}$, equal 
to ${d\xi_i\wedge d\eta_i\over{\xi_i\eta_i}}$ in local coordinates.

\begin{lmm} The $K$-affine structure on $\R^2\setminus \{(0,0)\}$
associated with  $\Omega$ coincides with the one 
associated with $\Omega^{can}$.
\end{lmm}
{\it Proof:} Using definitions from Section 7.2 one sees immediately that the statement of the Lemma 
  follows from the equality 
$p_{\Omega^{can}}(1+\eta)=1$, which is straightoforward:  
$p_{\Omega^{can}}(1+\eta)=exp\left(Res(\Omega^{can}\,\log(1+\eta))\right)=1\in {\cal O}_K^\times$. 
$\,\blacksquare$

In all the definitions and formulas in this section on
 can shift domains $U_i,_,i=1,2,3$ by vector $(x_0,0)\in\R^2$
 for arbitrary
 $x_0\in\R$, thus 
 giving a map $X^{an}\to \R^2$
 with  singularity at the point $(x_0,0)$.

Finally, we denote $\pi_{\ast}({\cal O}_{X^{an}})$ by ${\cal O}^{model}_{\R^2}$.
This will be our model for the sheaf ${\cal O}_B$ near each point of 
the singular set $B^{sing}$.

\section{Lines on surfaces}

In this section we are going to describe axiomatically the notion
of collection of lines on a surface.

\subsection{Data}
\begin{description}
\item[\bf a)] A compact oriented
surface $B$, a finite subset
$B^{sing}\subset B$.

\item[\bf b)] A $\Z$-affine structure on $Y=B^{sm}=B\setminus B^{sing}$ with the standard
singularities near each $b\in B^{sing}$.

\item[\bf c)] A set ${\cal L}$ of {\it lines}. 
With each line $l\in {\cal L}$ there is an associated  continuous map
$f_l:(0,+\infty)\to Y$. We assume that ${\cal L}$ is decomposed
into a disjoint union of two subsets 
${\cal L}={\cal L}_{in}\sqcup {\cal L}_{com}$. Lines belonging to
${\cal L}_{in}$ are called {\it initial}, while those in ${\cal L}_{com}$
are called {\it composite}. We assume that for any  $l\in {\cal L}$
 there exists a continuous extension $f_l:[0,+\infty)\to B$ such that
$f_l(0)\in B^{sing}$ if $l\in {\cal L}_{in}$ and $f_l(0)\in Y=B^{sm}$ if $l\in {\cal L}_{com}$.

\item[\bf d)] A collection of covariantly constant nowhere vanishing
integer-valued $1$-forms
$\alpha_l\in \Gamma((0,+\infty),f_l^{\ast} ((T^\ast)^\Z),
l\in {\cal L}$.
We assume that for $l\in{\cal L}_{in}$  in the
standard coordinates $(x,y)$ near singular point $f_l(0)$ 
we have: $f_l(t)=(0,t)$ or $f_l(t)=(0, -t)$ for all sufficiently small
$t>0$, and $\alpha_l(t)=\pm f_l^{\ast}(dy)$.

\item[\bf e)] A map ${\cal L}\to {\cal L}\times {\cal L},\,\,\,
l\mapsto (p_{left}(l),p_{right}(l))$ (the letter $p$ stands
for ``parent": one can think about these
lines as ``generating $l$ in a collision").

\end{description}

Notice that since the form $dy$ is invariant with respect to the monodromy,
the  condition in {\bf d)} is coordinate-independent. The covector $\alpha_l(t)$ will be
called a {\it direction covector} of $l$ at time $t$. It gives rise
to a half-plane 
$$P_{l,t}^{(0)}=\{v\in T_{f_l(t)}Y|\langle \alpha_l(t), v\rangle>0\}\,\,.$$

\subsection{Axioms}

To every $l\in {\cal L}_{in}$  we assign a pair 
$\left(f_l(0), sgn\left(\alpha_l(0)\right)\right)\in B^{sing}\times \{\pm 1\}$,
where   $sgn(\alpha_l(0))$ is a choice of sign
in $\pm f_l^{\ast}(dy)$ (see data {\bf d)} in the previous subsection).
In this way we obtain a map $r:{\cal L}_{in}\to B^{sing}\times \{\pm 1\}$.
\begin{description}
\item[\bf Axiom 1.] Map $r$ is one-to-one.
\end{description}

Let $U\subset Y$ be a simply-connected domain, and line $l$ intersects $U$.
Let $I\subset \R_+$ be an interval such that $f_l(I)\subset U$.
Then there exists a covariantly constant
closed non-zero $1$-form $\beta_U$ in $U$ (with constant integer coefficients), such that $f_l^{\ast}(\beta_U)=
\alpha_l$, when both sides are restricted to $I$.

\begin{description}
\item[\bf Axiom 2.] 
For any $t_1,t_2\in I$ one has 
$$\int_{f_l(t_1)}^{f_l(t_2)}\beta_U=t_2-t_1\,\,.$$
\end{description}

Let $l_1,l_2\in {\cal L}$, $t_1,t_2>0$ satisfy the condition
$f_{l_1}(t_1)=f_{l_2}(t_2)=x\in Y$. 
In this case we say that lines $l_1$ and $l_2$ have a collision at $x$
at the times $t_1$ and $t_2$ respectively.

\begin{description}
\item[\bf Axiom 3.] 
Under the above assumptions there are only two possibilities:

{\bf 3a)} either $l_1=l_2$ and $t_1=t_2$, or

{\bf 3b)} covector $\alpha_{l_1}(t_1)$ is not proportional to $\alpha_{l_2}(t_2)$.
Then we may assume that $\alpha_{l_1}(t_1)\wedge\alpha_{l_2}(t_2)>0$.
Under these conditions we require that for any coprime 
positive integers $n_1,n_2$  there exists
a unique line $l\in {\cal L}$ such that $l_1=p_{left(l)}, l_2=p_{right}(l)$,
$f_l(0)=x$ and $\alpha_l(0)=n_1\alpha_{l_1}(t_1)+n_2\alpha_{l_2}(t_2)$.

\end{description}

In other words, $l_1$ and $l_2$ are ``parents of $l$", and the direction
covector of $l$ at the intersection point is  a primitive integral linear
combination of those for $l_1$ and $l_2$ (see Figure 4).

\begin{figure}\label{figure4}
\centerline{\epsfbox{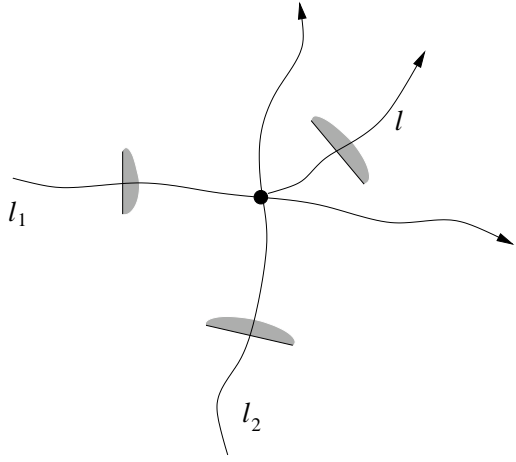}}
\caption{Line $l$ and its two parents $l_1, l_2$. Dashed half-planes are  domains in  tangent planes
 where
 $1$-forms $\alpha$ take positive values.}
\end{figure}

\begin{description}
\item[\bf Axiom 4.]  For every line $l\in {\cal L}_{com}$ there exist $l_1$ and
$l_2$ such that they satisfy the condition {\bf 3b)}.

\item[\bf Axiom 5.]   For any $x\in Y$ there are no more than two
pairs $(l,t)\in {\cal L}\times (0,+\infty)$ such that $x=f_l(t)$. 
In other words, there are no more than two lines intersecting at a point
in $Y$.
\end{description}

Let $l_1,l_2,t_1,t_2, x$ mean the same as in the  Axiom 3, 
 and assume that $\alpha_{l_1}(t_1)\wedge\alpha_{l_2}(t_2)>0$.
 Let us consider
the set ${\cal L}_{(x)}$ of germs of all $l\in {\cal L}_{com}$ starting at $x$
(i.e. such that $f_l(0)=x$).

\begin{description}
\item[\bf Axiom 6.]  
 For any finite subset ${\cal L}^{\prime}\subset{\cal L}_{(x)}$
there is an orientation preserving homeomorphism of a
neighborhood of $x$ onto a neighborhood of $(0,0)\in \R^2$ such that:

{\bf 6a)} Germs of oriented curves which are images of $l_1$ and $l_2$ get
transformed into the germs at $(0,0)$ of coordinate axes $(x,0)$ and $(0,y)$ 
respectively.

{\bf 6b) } Germ of the image of $l\in {\cal L}^{\prime}$ gets transformed
into the germ of the ray $\{(n_1t, n_2t)\,|\,t>0\}$ where
$\alpha_l(0)=n_1\alpha_{l_1}(t_1)+n_2\alpha_{l_2}(t_2)$.
\end{description}

Figure 5 illustrates this axiom.

\begin{figure}\label{figure5}
\centerline{\epsfbox{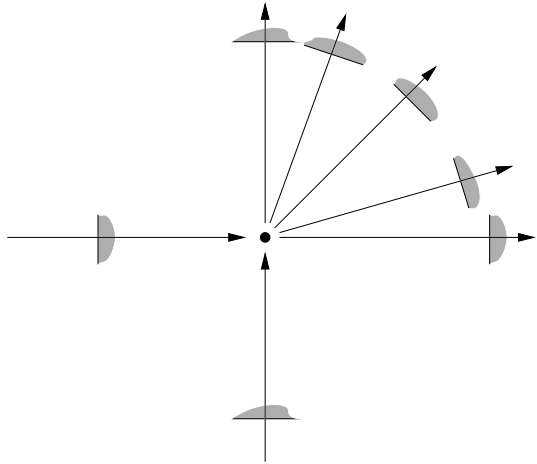}}
\caption{Two intersecting lines and some of  new lines obtained as a result of collision.  All lines are
straightened by a homeomorphism of $\R^2$. }
\end{figure}

\begin{description}
\item[\bf Axiom 7.]
 Let $p_i$ denotes either $p_{left}$ or $p_{right}$.
Then for any $l\in {\cal L}$ there exists $N\ge 1$ such that
if the line $p_1(p_2(\dots p_N(l)\dots)$ is well-defined
then it belongs to ${\cal L}_{in}$.
\end{description}

This axiom says that
any composed line $l\in {\cal L}_{com}$ appears as a
result of finitely many collisions. The tree of ancestors of a given line form a tree embedded in $B$, see Figure 6.

\begin{figure}\label{figure6}
\centerline{\epsfbox{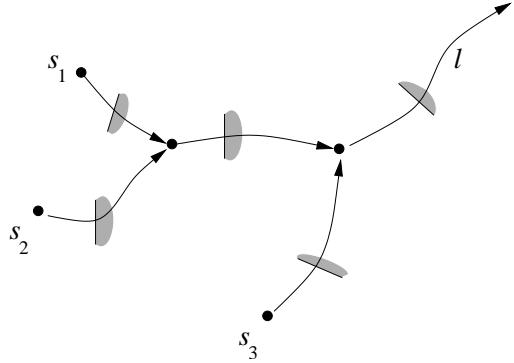}}
\caption{Tree of ancestors of line $l$ starting from $3$ singular points $s_1,s_2,s_3\in B^{sing}$.}
\end{figure}

\subsection{Example: gradient lines}

Here we offer a construction of the set of lines satisfying the above
axioms.

Let us use the standard $\R^2$ as a model around each $b\in B^{sing}$
in order to fix a structure
of smooth manifold on the whole surface $B$. Let $\widetilde{Y}$
denotes the covering of $Y$ such that the fiber over $y\in Y$ 
is $(T_{y}^{\ast}Y)^{\Z}\setminus\{ 0\}$.

 Let us fix a \emph{generic} 
smooth metric on $B$. By the pull-back it gives a metric on 
$\widetilde{Y}$. Notice that there is a canonical closed $1$-form $\beta$
on $\widetilde{Y}$ such that $\beta_{|(y,\mu)}=\mu$, where
$y\in Y,\,\, \mu\in (T_{y}^{\ast}Y)^{\Z}$.
Using the metric we obtain dual to $\beta$ gradient vector field $v$ on $\widetilde{Y}$.

For any $s\in B^{sing}$ and a choice of 1-form $\alpha(0)=\pm dy$ in local coordinates,
 we take the unique integral line of $v$ starting at $(s,\alpha(0))$.
Set ${\cal L}_{in}$ will be the set of all  lines obtained in this way.
Each line $l\in {\cal L}_{in}$ carries a covariantly constant closed $1$-form $\alpha_l$.
Using Axiom 2 as a definition, we obtain a
canonical parametrization of each line by the time parameter $t$.
 Since the metric is generic, a line cannot return to a point in $B^{sing}$.

Then we proceed inductively. If two already constructed lines
 $l_1,l_2\in {\cal L}$ meet at $x\in Y$ we produce a new
integral line $l$ of $v$ with the direction covector satisfying 
the condition {\bf 3b)}  for any
pair of coprime positive integers $n_1, n_2$.
In this way we construct a set of lines
${\cal L}$ satisfying all the axioms. The only non-trivial thing
to check is that for each line values of the parameter $t$ 
are in one-to-one correspondence with the interval $(0,+\infty)$.
In order to see this we observe that the length of each line is infinite. 
Indeed, an integral curve of $v$ cannot have a limiting point in $Y$ (since
the flow generated by $v$ is smooth, and the lengths of tangent
vectors are bounded from below because of the integrality
of $1$-forms).

We conclude that there exists a set ${\cal L}$ of lines satisfying
Axioms 1-7.

\section{Groups and symplectomorphisms}

In this section we are going to discuss the sheaf of
groups of symplectomorphisms $Symp:=Symp({\cal O}_Y^{can})$
of the sheaf ${\cal O}_Y^{can}$. Let $U\subset Y$ be an open convex
subset. By definition,
a {\it symplectomorphism} of ${\cal O}_Y^{can}(U)$
is an automorphism of  $K$-algebra ${\cal O}_Y^{can}(U)$ preserving projection
 to $Y$ and 
  the canonical symplectic form 
$\Omega={d\xi\wedge d\eta\over{\xi\eta}}$ (the latter is understood
as an element of the algebra of K\"ahler differential forms).
To each line $l$ we will assign a symplectomorphism of the restriction
of ${\cal O}_Y^{can}$ to $l$,
so that the assignment will be compatible with the collision
of lines. Then we are going to modify the sheaf ${\cal O}_Y^{can}$
using symplectomorphisms, associated with lines and obtain
the sheaf ${\cal O}_Y^{modif}$. This sheaf will be glued 
with the sheaf ${\cal O}^{model}_{\R ^2}$ near each point of $B^{sing}$.

\subsection{Pro-nilpotent Lie algebra}

Here it will be  convenient to work in local coordinates
$(x,y)=(\log|\xi|, \log|\eta|)$ on $Y$.

Let $(x_0,y_0)\in \R^2$ be a point, 
$\alpha_1,\alpha_2\in (\Z^2)^\ast$ be $1$-covectors such that
$\alpha_1\wedge\alpha_2>0$. 
Denote by $V=V_{(x_0,y_0),\alpha_1,\alpha_2}$ the closed angle
$$\{(x,y)\in \R^2|\langle \alpha_i, (x,y)-(x_0,y_0)\rangle\ge 0, i=1,2\,\}\,\,.$$

Let ${\cal O}(V)$ be a $K$-algebra consisting
of series $f=\sum_{n,m \in \Z}c_{n,m}\xi^{n}\eta^{m}$,
such that $c_{n,m}\in K$ and
for all $(x,y)\in V$ we have:

\begin{enumerate}
\item if $c_{n,m}\ne 0$ then $\langle (n,m),(x,y)-(x_0,y_0)\rangle\le 0$, 
where we identified $(n,m)\in \Z^2$
with a covector in $(T_{p}^{\ast}Y)^{\Z}$;
\item $\log|c_{n,m}|+nx+my\to -\infty$ as long as $|n|+|m|\to +\infty$.
\end{enumerate}

Algebra ${\cal O}(V)$ is  a Poisson algebra with respect to
the bracket $\{\xi,\eta\}=\xi\eta$.

For an integer covector $\mu=adx+bdy\in (\Z^2)^*$ we denote
by $R_{\mu}$ the monomial $\xi^a\eta^b$.

Let us consider a pro-nilpotent Lie algebra 
$\g:=\g_{\alpha_1,\alpha_2,V}\subset {\cal O}(V)$ consisting of series
$$f=
\sum_{n_1,n_2\ge 0,n_1+n_2>0}c_{n_1,n_2}R_{\alpha_1}^{-n_1}R_{\alpha_2}^{-n_2} $$
satisfying the condition
$$\log|c_{n,m}|-n_1\langle \alpha_1, (x,y)\rangle-n_2\langle \alpha_2, 
(x,y)\rangle\le 0\,\,\,\,\forall \,(x,y)\in V\,\,.$$
The latter condition is equivalent to 
$\log|c_{n,m}|-\langle n_1 \alpha_1+n_2\alpha_2, (x_0,y_0)\rangle\le 0$.

Lie algebra $\g$ admits a filtration by Lie subalgebras 
$\g^{\ge k}, k\in \Z, k\ge 1$, $\g=\g^{\ge 1}$, such that
$\g^{\ge k}$ consists of the above series which satisfy the condition
$n_1+n_2\ge k$. Clearly 
$[\g^{\ge k_1},\g^{\ge k_2}]\subset \g^{\ge k_1+k_2}$, and
$\g=\varprojlim_{k\to +\infty}\g/\g^{\ge k}$.

Thus, $\g$ is a topological complete
pro-nilpotent Lie algebra over $K$. We denote by $G$ 
the corresponding pro-nilpotent Lie group $\exp(\g)$. It inherits the filtration
by normal subgroups $G^{\ge k}$ obtained from the corresponding
Lie algebras.

\subsection{Lie groups $G_{\lambda}$}

For each $\lambda \in [0,+\infty]_{\Q}:=\Q_{\ge 0}\cup\infty$ we define a Lie subalgebra 
$$\g_{\lambda}=
\left\{\sum_{n_1,n_2}c_{m,n}R_{\alpha_1}^{-n_1}R_{\alpha_2}^{-n_2}\in \g
\,|\,\,c_{n_1,n_2}\in K, \,\,{n_2\over{n_1}}=\lambda\,\right\}\,\,.$$
Each $\g_{\lambda}$ is an abelian Lie algebra. It carries
the induced filtration by Lie algebras 
$\g_{\lambda}^{\ge k}=\g_{\lambda}\cap \g^{\ge k}$. Denote by $G_{\lambda}=\exp(\g_{\lambda})$
 the corresponding pro-nilpotent group.

\begin{lmm}
For any given  $k\ge 1$ there exist finitely
many $\lambda_1<\lambda_2<\dots<\lambda_{N_k}$ such that
$\g_{\lambda}/\g^{\ge k}_\lambda=0$ 
for $\lambda\ne \lambda_i, 1\le i\le N_k$.

\end{lmm}

{\it Proof.} Indeed, for the
monomial $R_{\alpha_1}^{-n_1}R_{\alpha_2}^{-n_2}\in \g_\lambda$ which maps non-trivially 
 to the quotient $\g_{\lambda}/\g^{\ge k}_\lambda$
we have: $n_1+n_2\le k, n_1/n_2=\lambda$,
where $n_1,n_2$ are non-negative integers. There are finitely many
such non-negative integers $n_1$ and $n_2$. $\blacksquare$

It follows from the Lemma that we have a natural isomorphism of vector spaces
$\prod_{\lambda \in [0,+\infty]_{\Q}}\g_{\lambda}/\g_{\lambda}^{\ge k}\to
\g/\g^{\ge k}$, hence the map 
$$(f_{\lambda})_{\lambda\in [0,+\infty]_{\Q}}\mapsto \sum_{\lambda}f_{\lambda}=\sum_{i=1}^{N_i} f_{\lambda_i},\,\,\,\,
\mbox{ where }f_\lambda\in \g_{\lambda}/\g_{\lambda}^{\ge k}\,\,\,\forall \lambda \in [0,+\infty]_{\Q} $$
is well-defined and 
gives rise (after taking the projective limit as $k\to +\infty$)
to the isomorphism 
$\g\simeq\prod_{\lambda \in [0,+\infty]_{\Q}}\g_{\lambda}$.

In a similar way we define the map
$\prod_\to:\prod_{\lambda \in [0,+\infty]_{\Q}}G_{\lambda}\to G$,
the product is taken with respect to the natural
order on $\Q$. Namely, for any $k\ge 1$ we define
$${\textstyle \prod_\to^{(k)}}:\prod_{i=1}^{N_k} G_{\lambda_i}/G_{\lambda_i}^{\ge k} 
\to G/G^{\ge k} \,,\, (g_1,\dots,g_{N_k})\mapsto g_1\dots g_{N_k},\,
\mbox{ for } g_i\in G_{\lambda_i}/G_{\lambda_i}^{\ge k}$$ 
and then set $\prod_\to:=\varprojlim_k \prod_\to^{(k)}\,$.
\begin{thm} Map $\prod_\to$ is a bijection of sets.

\end{thm}

{\it Proof.} Let $k\ge 1$ be an integer.  We claim that $\prod_\to^{(k)}$ is a bijection
of sets (this implies the proposition by taking the projective limit
as $k\to +\infty$). We will prove the bijection by induction in $k$ . 
Case $k=1$ is obvious because all the groups under considerations
 are trivial.

 We would like
to prove that  $\prod_\to^{(k+1)}$ is a bijection assuming that $\prod_\to^{(k)}$ is a bijection.
 Let $h$ be an element of $ G/G^{\ge k+1}$ and 
$\overline{h}$ its image  in $G/G^{\ge k}$. 
By the induction assumption there exist  unique 
$\overline{h}_i\in G_{\lambda_i}/G_{\lambda_i}^{\ge k+1}, 1\le i\le N_{k+1}$
such that $\overline{h}_1\dots\overline{h}_{N_{k+1}}=\overline{h}$.
Let $h_i,1\le i\le N_{k+1}$ be any liftings of $\overline{h}_i$
to $G_i/G_i^{\ge k}$. Then $h_1\dots h_{N_{k+1}}=h\pmod{G^{\ge k}}$, hence
$c:=h_1\dots h_{N_{k+1}}h^{-1}$ belongs to $G^{\ge k}/G^{\ge k+1}\subset Center(G/G^{\ge k+1})$.
The last inclusion holds because 
$[\g,\g^{\ge k}]=[\g^{\ge 1},\g^{\ge k}]\subset \g^{\ge k+1}$.

Next we observe that the isomorphism of abelian Lie algebras

$$\bigoplus_{1\le i\le {N_{k+1}}}\g_{\lambda_i}^{\ge k}/\g_{\lambda_i}^{\ge k+1}\simeq \g^{\ge k}/\g^{\ge k+1}$$
implies an isomorphism of the corresponding abelian groups

$$\prod_{1\le i\le {N_{k+1}}}G_i^{\ge k}/G_i^{\ge k+1}\simeq G^{\ge k}/G^{\ge k+1}\,\,.$$

Hence we can write uniquely $c=c_1\dots c_{N_{k+1}}$, where $c_i\in G^{\ge k}/G^{\ge k+1}\subset Center(G/G^{\ge k+1})$.
It follows that $\prod_\to^{(k+1)}\left((h_ic_i^{-1})\right)=h$. Also it is now 
clear that this decomposition of $h$ 
is unique. This concludes the proof. $\blacksquare$

\subsection{Function $ord_l$}

For $l\in {\cal L}$ we will  define an {\it order} function
$$ord_l \in \Gamma((0,+\infty), f_l^{\ast}(Aff_{\Z,Y}))$$
 (its meaning will become clear later) 
 by the following inductive procedure:
\begin{enumerate}
\item Let $l\in {\cal L}_{in}$ and $t>0$ be sufficiently small.
Then in the standard affine coordinates near 
$s=f_l(0)$ one has $\alpha_l=\pm f_l^{\ast}(dy)$.
We define $ord_l=\pm f_l^{\ast}(y)$. Then $d(ord_l)=\alpha_l$, and we can extend
uniquely $ord_l$ for all $t\in (0,+\infty)$.
\item Let $l\in {\cal L}_{com}$ and $l_1, l_2$ be parents of $l$. In the notation of Axiom 3
 we have  
 $f_{l_1}(t_1)=f_{l_2}(t_2)=f_l(0)$ and $\alpha_l(0)=
n_1\alpha_{l_1}(t_1)+n_2\alpha_{l_2}(t_2)$. Then we define
$ord_l(0):=n_1ord_{l_1}(t_1)+n_2ord_{l_2}(t_2)$. Again, using
the condition $d(ord_l)=\alpha_l$ and the knowledge of $ord_l(0)$
we can extend $ord_l$ for $t>0$.
\end{enumerate}
Notice that $ord_l(t)$ can be thought of as affine function
on the tangent space $T_{f_l(t)}Y$ (in the induced integral
affine structure). In particular, we have a half-plane
$P_{l,t}\subset T_{f_l(t)}Y$ defined by the inequality
$ord_l(t)>0$. 
The family of half-planes $P_{l,t}$ is covariantly
constant with respect to $\nabla^{aff}$.

Each half-plane $P_{l,t}$ contains $0\in T_{f_l(t)}Y$ strictly in its interior.
Recall that at the end of Section 9.1 we defined another half-plane
$P_{l,t}^{(0)}\subset T_{f_l(t)}Y$. It is easy to see that $P_{l,t}^{(0)}$ is the half-plane 
parallel to $P_{l,t}$ such that $0\in T_{f_l(t)}Y$ is on the boundary of $P_{l,t}^{(0)}$.

\subsection{Symplectomorphisms assigned to lines}

In this section we are going to assign to each line 
$l\in {\cal L}$ a symplectomorphism $$\varphi_l\in \Gamma\left((0,+\infty),f_l^* \left(Symp\right)\right)$$ 
 giving for each $t>0$ a transformation 
$\varphi_l(t): {\cal O}^{can}_{Y,f_l(t)}\to {\cal O}^{can}_{Y,f_l(t)}$. 
This symplectomorphism in local coordinates will belong to the subgroup $G_\lambda$ where $\lambda$ is the slope
 of $\alpha_l(t)$. More precisely, we demand that $\varphi_l(t)$ is of the form
$$\varphi_l(t)=\exp\{F_{l,t}(\xi^{-a}\eta^{-b}),\cdot\}\,\,,$$
where $\alpha_l(t)=a dx+ bdy$, operation $\{\cdot,\cdot\}$ is the  Poisson bracket on ${\cal O}^{can}_{Y,f_l(t)}$
 and $F_{l,t}(z)\in zK[[z]]$ is an 
analytic function of \emph{one} variable  satisfying the following condition.
  Let us consider  the pullback (by the exponential map) of the function $F_{l,t}(\xi^{-a}\eta^{-b})$
 to a section of the sheaf ${\cal O}^{can}$ on vector space
 $T_{f_l(t)}Y\simeq \R^2$ considered as a manifold with $\Z$-affine structure. 
Then this pullback should  admit an analytic continuation from $0\in T_{f_l(t)}$ to the half-plane
 $P_{l,t}$, and obey there the  bound
 $$|F_{l,t}(\xi^{-a}\eta^{-b})|\le \exp(-ord_l(t))\,\,.$$

Let us explain the construction of $\varphi_l(t)$, leaving the justification
for the next sections.

Symplectomorphisms $\varphi_l$ are constructed by an inductive procedure.
Let $l=l_+\in {\cal L}_{in}$ be (in standard affine coordinates)
a line in the half-plane $y>0$ emerging from $(0,0)$
(there is another such line $l_{-}$
in the half-plane $y<0$). Assume that $t$
is sufficiently small. Then we define
$\varphi_l(t)\in Symp_{f_l(t)}$ 
on topological generators $\xi,\eta$ by the formula (as in Section 8)

$$\varphi_l(t)(\xi,\eta)=(\xi(1+1/\eta),\eta)\,\,.$$

Notice that $\varphi_l(t)=\exp\{F(\eta^{-1}),\cdot\}$, where
$F(z)=\sum_{n>0}(-1)^n z^n/n^2$ is convergent for $|z|<1$.

In order to extend $\varphi_l(t)$ to the interval $(0,t_0)$, where
$t_0$ is not small, we cover the corresponding segment of $l$ by
open charts. Notice that change of affine coordinates transforms
$\eta$ into a monomial multiplied by a constant from $K^{\times}$.
Therefore $\eta$ extends analytically in a unique way to a global section over $(0,+\infty)$
 of the sheaf $f_l^{\ast}(({\cal O}^{can})^{\times})$.
Moreover, the norm $|\eta|$ strictly decreases as $t$ increases,
and remains strictly smaller than $1$. Hence $F(\eta)$
can be canonically extended  for all $t>0$.

Each symplectomorphism $\varphi_l(t)$ is defined by a series
which converges in the half-plane
$P_{l,t}$. Using the exponential map associated with the affine structure
as well as estimates of $ord_l(t)$, we can
extend analytically $\varphi_l(t)$ into a neighborhood
of $f_l(t)$. 

Let us now assume that $l_1$ and $l_2$ collide at 
$p=f_{l_1}(t_1)=f_{l_2}(t_2)$,
generating the line $l\in {\cal L}_{com}$. 
Then $\varphi_l(0)$ is defined with the help of factorization theorem
in the group $G$. More precisely, we set $\alpha_i:=\alpha_{l_i}(t_i),\,\,i=1,2$ and the angle
 $V$ to be the intersection of half-planes $P_{l_1,t_1}\cap P_{l_2,t_2}$.
By construction  elements
 $g_0:=\varphi_{l_1}(t_1)$  and $g_{+\infty}:=\varphi_{l_2}(t_2)$ belong respectively 
to $G_0$ and $G_{+\infty}$. Then we can use the
factorization Theorem 6 and write down the formula

$$g_{+\infty} g_0= {\textstyle \prod_\to}\left( 
(g_\lambda)_{\lambda\in [0,+\infty]_\Q}\right)=
g_0\dots g_{1/2}\dots g_1 \dots
g_{+\infty}\,\,,$$
where $g_{\lambda}\in G_\lambda$ and the product on the right
 is in the \emph{increasing} order.
 There is no clash of notations because it is easy to see that the boundary factors in the decomposition
 from above are indeed equal to $g_0$ and $g_{+\infty}$.
 Each term $g_\lambda$ with $0<\lambda= n_1/n_2 <+\infty$ 
corresponds to the newborn line $l$ with the 
direction covector
$n_1\alpha_{l_1}(t_1)+n_2\alpha_{l_2}(t_2)$.
Then we set  $\varphi_l(0):=g_{\lambda}$. This 
transformation is defined by a series
which is convergent in a neighborhood of $p$, and using the analytic continuation 
as above, we obtain $\varphi_l(t)$ for $t>0$.
The decomposition identity can be rewritten as
$$g_0\dots g_{1/2}\dots g_1 \dots
g_{+\infty}g_0^{-1}g_{+\infty}^{-1}=id$$
where each factor corresponds to half-lines at the collision point
 (see Figure 5), and the meaning of the identity is that
 the  infinite composition of symplectomorphisms in the natural
 cyclic order on half-lines, is trivial.

\section{Modification of the sheaf ${\cal O}^{can}$}

\subsection{Pieces of lines and convergence regions}

\begin{defn}  A  neighborhood $U$ of a point $x\in Y$ is \emph{convex} 
 if there exists an open convex $U_1\in T_xY,
0\in U_1$ which is isomorphic to $U$ by means of the exponential map
$\exp_x:T_xY\to Y$  associated with the affine structure on $Y$.

\end{defn}

For $x\in Y$ let $U\subset U'$ be convex neighborhoods of $x$ such that
$U$ is relatively compact in $U'$. Let $l\in {\cal L}$. Then there is a natural embedding 
$f_l^{-1}(U)\to f_l^{-1}(U')$.

\begin{defn} A \emph{piece} of $l$ defined by the pair $(U,U')$ is an element of the
image of the set of connected components
$\pi_0(f_l^{-1}(U))$ into $\pi_0(f_l^{-1}(U'))$ under
the above embedding.

\end{defn}

In plain words a piece $L$ of $l$ is an equivalence class of
a connected interval of  $l\cap U$. Two connected intervals
are equivalent if they are contained in a larger connected
interval of $l\cap U'$. The sole purpose of the introduction of the notion
 of a piece is to avoid some pathology. Namely, for any pair $(U,U')$ as above,
  any $l\in{\cal L}$ and any $T\in\R_{>0}$, 
 there is only a finite  number of pieces of $l$ in $(U,U')$ which have points with time parameter
 $t\in (0,T)$.

Let $L$ be a piece of $l$ defined by a pair $(U,U')$. Then one
can define an affine function $ord_L\in Aff_{\Z,Y}(U')$ in the following way.
Let $t>0$ be such that $f_l(t)$ belongs to $L$. Since $U'$ is
convex, there is a unique continuation of $ord_l(t)$ to $U'$.
This is an affine function which does not depend on the choice of $t$.
We will denote it by $ord_L$.

For any germ of a symplectomorphism $\varphi\in Symp_p$ at a point $p\in Y$ we  define
its convergence region as the maximal convex subset
$\Omega(\varphi)\subset T_p Y$ such that
 the pullback 
$\exp_p^{\ast}(\varphi)$ extends to $\Omega(\varphi)$.
Since the definition of $\varphi_l$ (and hence its convergence
region) is covariant with respect to the affine connection we
have the following result:

\begin{prp} Let $p=f_l(t)$ belongs to a line $l$.
Then the convergence
region of $\varphi_l(t)$ at $p$ contains an open half-plane $P_{l,t}$.

\end{prp}
It is clear that one can define convergence regions for symplectomorphisms associated with pieces of lines,
 and a similar property holds for them.

\subsection{Main assumptions, and an apology}

Let us suppose that our collection of lines  satisfies the following
assumptions:

\begin{description}
\item[\bf Assumption A1]
 There is a smooth metric $g=g_B$ and a collection of balls $D(s,r_s)$
with centers at $s\in B^{sing}$ such that each ball $D(s,r_s)$ contains exactly
two lines $l_{\pm}\in {\cal L}_{in}$ outcoming of $s$.

\item[\bf Assumption A2] There exists $\varepsilon>0$ such that
for any $p=f_l(t)\in Y^{\prime}:=B\setminus \cup_{s\in B^{sing}}D(s,r_s)$
the distance
in $T_pY$ between $0\in T_pY$ and the boundary  of $P_{l,t}$
is greater  or equal to $\varepsilon$.
\end{description}
We are going to show that such a collection does exist in Section 11.

Assumptions {\bf A1} and {\bf A2} are very artificial, they do not hold in physical picture
 which is the main motivation for the construction. It is quite possible that they can be weakened
 or even omitted. The main purpose of introducing them here is the possibility to
 define the sheaf
 of analytic functions by simple gluing. In complex geometry it is similar to the  gluing
of closed Riemann surfaces with boundaries by the mean of real-analytic
identifications of the boundaries. It is well-known that
one can replace real-analytic maps by smooth ones (or even by
quasi-symmetric continuous maps). Maybe the rest of this section is unnecessary, and 
 unpleasant technical arguments in Section 11.5 can be avoided.

\subsection{Infinite product and its convergence}

Denote by $W_{\cal L}:=\cup_{l\in {\cal L}}f_l([0,+\infty))$ the set of all points of all lines.
 It has measure zero.
Let $p$ be a point of $Y$.
We consider two convex neighborhoods $U\subset U'$  of $p$ such  that $U$ 
is relatively compact in $U'$.
 
For any two points $x,y$ belonging to  $U\setminus W_{\cal L}$,
and a path
$\gamma$ joining $x$ and $y$ in $U$, we would like to define an 
infinite ordered product $i_{x,y}^{\gamma}$
 of transformations $\varphi_L^{\pm 1}$, where factors
 correspond to the
 intersection points of $\gamma$ with all possible pieces $L$
 relative to $(U,U')$. Factors in the infinite product are ordered
 according to the time parameter of $\gamma$, the sign corresponds to 
the mutual position of orientations of $\gamma$ and a piece $L$ at
 the intersection point.

In order to give a precise meaning to the infinite product the 
neighborhood $U$ of $p$ should be  sufficiently small. 
Then we will have an analytic continuation
 of 
symplectomorphisms $\varphi_L$  to $U$, and  the convergence of
  the infinite product.
We are also going  to prove that the product  is independent of the
 choice of path $\gamma$.
In order to achieve these goals it suffices to assume:
\begin{description}
\item[\bf C1] for any $l,t$ such $f_l(t)\in U$ the set $(\exp_{f_l(t)})^{-1}(U)$ is contained in
   $P_{l,t}\,$;
\item[\bf C2] for any $C\in \R$ there is only a finite number of pieces
  $L$ of lines in $U$ such that $\inf_{x\in U} ord_L(x)< C$.
\end{description}

\begin{thm} Assume two above conditions.
Then the product defining $i_{x,y}^{\gamma}$ converges at every point of $U$
and in fact gives an element of $Symp(U)$.
Moreover, the product does not depend on the choice of path $\gamma$,
and for any $x,y,z\in U\setminus W_{\cal L}$ 
satisfies the relation $i_{x,y}i_{y,z}=i_{x,z}$.

\end{thm}

{\it Proof.} 
 Condition {\bf C1} implies that all transformations $\varphi_L$ admit an analytic continuation to $U$.
 Let us introduce a decreasing filtration by positive real numbers
  $Symp^{\ge r}(U),\,\,r\in \R,r\ge 0$ on group $Symp(U)$ by the formula
 $$\left\{g\in Symp(U)\,|\,\,\log|\xi'/\xi-1|,\log|\eta'/\eta-1|< -r\,\mbox 
{ where }(\xi',\eta')=g((\xi,\eta))
  \,\right\}$$
This is a complete filtration, and condition {\bf C2} implies that in any quotient $Symp(U)/Symp^{\ge r}(U)$
 only a finite number of elements $\varphi_L$ are non-trivial.
 Therefore we can define the product in the quotient group. 

In order to prove independence of $\gamma$, we consider the quotient group 
$Symp(U)/Symp^{\ge r}(U)$, and the finite $1$-dimensional CW-complex
(graph) consisting of finitely many pieces $L$,
such that $\varphi_L\ne 1$ in 
the quotient. For each vertex $v$ of the graph
there is a natural cyclic order on the edges incident to $v$.
The product $\varphi_v=\prod_L \varphi_L^{\pm 1}$ taken in the cyclic order
over the set of edges incident to $v$ is equal to $id$ (this follows from
the construction of $\varphi_l$ via factorizations). Since $U$ is simply-connected, 
we conclude that
the image of $i_{x,y}^{\gamma}$ in $Symp(U)/Symp^{\ge r}(U)$
does not depend on $\gamma$. Using completeness of the filtration we
see that $i_{x,y}:=i_{x,y}^{\gamma}$ does not depend on $\gamma$. Proof
of the identity $i_{x,y}i_{y,z}=i_{x,z}$ is similar.
$\blacksquare$

\begin{thm} Assumptions {\bf A1} and {\bf A2} imply that for any $p\in Y$ there exist
 neighborhood $U$ (and also $U'$) satisfying conditions {\bf C1} and {\bf C2}.

\end{thm}

{\it Proof.} Assumption {\bf A1} implies that the result near any singular point $s\in B^{sing}\,$,
 as there are only two lines near $s$. If we are far from $B^{sing}$ then obviously {\bf A2} implies {\bf C1}.

In order to check {\bf C2} we prove the following lemma
\begin{lmm} Under Assumptions {\bf A1} and {\bf A2}, for any $C>0$ the set 
$$\{(l,t)| \,\,ord_l(t)(f_l(t))<C\,\}\subset {\cal L}\times (0,+\infty)$$
 consists of a finite number of intervals.
\end{lmm}
{\it Proof:}
We proceed by induction in ``complexity
of the line". Let $\delta \in \R_{>0}$ be the infimum of 
$ord_l(t)(f_l(t))$ where $l \in  {\cal L}_{in}$ has a collision at time  $t$.
 This number is strictly positive because the number of initial lines is finite, and by {\bf A1}
 there is no collisions at small times. 
 Observe that the value of $ord_l(0)$ at
the beginning of any composite line $l$ is greater or equal
to the sum $ord_{l_1}(t_1)+ord_{l_2}(t_2)$.
 Therefore the inequality in the lemma 
implies that the number of collisions is bounded from above by $C/\delta$.
 Also we have an upper bound on integer coefficients $(n_1,n_2)$ in each collision 
 (see Axiom $\bf 3b)$ in Section 9.2).
Let us observe that the length of each edge of the ansector tree of $l$
is also bounded from above by $A\,ord_l$,
for some absolute constant $A>0$. Hence we have only finitely many possibilities for intersections.
$\blacksquare$

For point $p\in Y$ which is far from $B^{sing}$ we chose as $U$ a neighborhood of radius $\epsilon'\ll\epsilon$ where
 $\epsilon>0$  is constant from Assumption {\bf A2}. Then for any point of a line $f_l(t)\in U$ we will have the inclusion
 $$U\subset \exp_{f_l(t)} \left(\frac{1}{2} P_{l,t}\right)\,\,.$$
This implies that $ord_L$ in $U$ for the corresponding piece $L$ is bounded below by 
$$\frac{1}{2} ord_l(t)(f_l(t))\,\,.$$
Since (by the last lemma) there exists only a finite number of pieces $L$ intersecting such $U$, we obtain convergence condition
{\bf C2}. $\blacksquare$

\subsection{Construction of the modified sheaf ${\cal O}^{modif}_B$}

For any point $p\in Y$ and a neighborhood $U$ satisfying conditions {\bf C1} and {\bf C2}
 we define the sheaf ${\cal O}^{modif}_U$ as the result of the identification of copies of the sheaf
 $\left({\cal O}^{can}_Y\right)_{|U}$ labeled by points $x\in U\setminus W_{\cal L}$,
 by isomorphisms $i_{x,y}$.
 It  follows from formulas in Section 8 that near singular points one can identify canonically this sheaf with the restriction 
of the sheaf ${\cal O}^{model}_{\R^2}$ to a punctured neighborhood of $(0,0)\in \R^2$.

\begin{prp} For the modified sheaf ${\cal O}^{modif}$ one has
a canonical nowhere vanishing section $\Omega$ of the associated sheaf 
of $K$-analytic $2$-forms. 

The $K$-affine structure $Aff^\Omega_{K,Y}$ on $Y$ associated with $\Omega$ coincides
with the initial one $Aff_{K,Y}$.

\end{prp}

{\it Proof.} Existence of $\Omega$ follows from the fact that all
modifications associated with lines are symplectomorphisms.
In order to finish the proof it suffices to check that the modification
associated with a line does not change the $K$-affine structure on $Y$.
In local coordinates we may assume that 
$\Omega={d\xi\over \xi}\wedge{d\eta\over \eta}$ and the modification
is of the form $\varphi(\xi,\eta)=(\xi f(\eta^{-1}),\eta)$, where
$f(z)=1+\sum_{n\ge 1}c_n z^n \in K[[z]]$ is convergent
in an appropriate domain.
We need to check that the automorphism $\varphi$ acts
trivially on the quotient sheaf 
$\pi_{\ast}({\cal O}_X^{\times})/\ker\,p_{\Omega}$ (see Section 7.2
for the notation). This check reduces to the calculation
of
$$p_{\Omega}\left({\xi f(\eta)\over \xi}\right)=
\exp\left({Res(\Omega \log(\xi f(\eta)/\xi)))\over Res(\Omega)}\right)\,\,.$$
The latter is equal to 
$\exp\left(Res(\Omega \log(f(\eta)))\right)=1$ because 
$\log(f(\eta^{-1}))$ belongs to $\eta^{-1} K[[\eta^{-1}]]$ and therefore has no constant term.
$\blacksquare$

Thus, we have a solution of the Lifting Problem under Assumptions {\bf A1} and {\bf A2}.

\subsection{Construction of the collection of lines}

We would like to show that there exists a smooth metric $g$ 
and a collection of lines satisfying the Assumptions {\bf A1} and {\bf A2}.

Let $g_0$ be an arbitrary smooth metric, flat near
singular points. We define germs of lines $l\in {\cal L}_{in}$
in such a way that for each $s\in B^{sing}$ in local coordinates 
these lines are given by $\{(0,y)|y>0\}$ and $\{(0,y)|y<0\}$.
The metric $g$ will coincide with $g_0$ in a sufficiently small
neighborhood $U=\cup_{s\in B^{sing}}D(s,r_s)$ of the singular set.
Hence Assumption {\bf A1} will be satisfied.
 
In order to construct the whole family of lines we introduce
a $3$-dimen\-sional manifold ${\cal M}$ consisting of
pairs $(x,P)$ where $x\in B\setminus \overline{U}_1$ and $P$
is a half-plane in $T_xB$ whose boundary contains zero. Here $U_1:=\cup_{s\in B^{sing}}D(s,2r_s)$
 is  a larger neighborhood  of $B^{sing}$.

We would like to construct a smooth section $v: (x,P)\mapsto v_{(x,P)}\in T_xB$
of the pull-back to ${\cal M}$ of the tangent bundle $TB$ satisfying
the following conditions:
\begin{enumerate}
\item for any $(x,P)\in {\cal M}$ one has $v_{(x,P)}\in int(P)\,\,$;
\item for any $x\in B\setminus \overline{U}_1$ the map
$(x,P)\mapsto {\R}_{>0}^{\times}\cdot v_{(x,P)}$ is an
orientation-preserving diffeomorphism
$$ S^1\simeq (T_xB^{\ast}\setminus\{0\})/{\R}_{>0}^{\times}\to S^1\simeq 
(T_xB\setminus \{0\})/{\R}_{>0}^{\times}\,\,\,;$$
\item for every $l\in {\cal L}_{in}$ there exists a smooth extension
of the piece of $l$ in $U_1$ to a larger piece intersecting $\partial U_1$
such that
$$\dot{f}_l(t)\in  {\R}_{>0}^{\times}\cdot v_{(f_l(t),P_{l,t})}\,\,,$$
for such $t>0$ that $f_l(t)\in B\setminus \overline{U}_1\,\,$;
\end{enumerate}

Let us associate with the section $v$ a nowhere vanishing 
vector field $\hat{v}$ on
$T^{\ast}(B\setminus \overline{U}_1)\setminus (Zero\,\,\,Section)$
in the following way:
\begin{itemize}
\item For each $(x,\alpha)\in T_x^{\ast}B$ the 
vector $\hat{v}(x,\alpha)$ is tangent to the horizontal distribution
associated with the flat connection $\nabla$ (the one which
defines the affine structure on $B\setminus B^{sing}$).
 \item Projection of $\hat{v}(x,\alpha)$ to $B$ coincides with
$v_{(x,P_\alpha)}$, where $P_\alpha=\{\gamma|(\alpha,\gamma)>0\}$.
\end{itemize}
Clearly these conditions determines $\hat{v}$ uniquely.
Now we formulate last condition:
\begin{enumerate}
\item[4.] there exits $r_s'>2 r_s$ such that for almost all (in the sense of Baire category) initial
values $(x_0,P_0)\in {\cal M}$  the integral curve
of $\hat{v}$ starting at $(x_0,P_0)$ reaches the pullback of
$B\setminus \cup_{s\in B^{sing}} D(s,r_s')$ in finite time.
\end{enumerate}

Using the vector field $\hat{v}$ we will construct
(under certain genericity assumptions) a set ${\cal L}$
of lines satisfying Assumption {\bf A1}. Namely, the data
consisting of a line $l$ and an integer-valued $1$-form $\alpha_l$
(see Section 9) will be
 an integral line of $\hat{v}$.

We are going to construct lines by induction by the number of collisions.
 Lines $l\in {\cal L}_{in}$ will be constructed using condition 3.
The genericity assumption mentioned after the condition 4 is the
assumption that no more than two lines collide  and that initial values
 for newborn lines will be sufficiently generic. Conditions 1 and 4 plus genericity 
imply that one can parametrize any line $l\in {\cal L}$ by
the new ``time" $t>0$ such that the Axiom 2 is satisfied.
Axiom 6 follows from the condition 2. Other axioms and
the Assumption {\bf A1} will be satisfied automatically.

Now we would like to discuss Assumption {\bf A2}.

\begin{prp} Suppose that the metric $g$ and field $v$ described above
are such that for any $(x,P)\in {\cal M}$ there exists $C>0$
such that 
$$(\nabla_{v_{(x,P)}}\,g)(n_P,n_P)\le C\,g(n_P,v_{(x,P)})\,\,,$$
where $n_P$ is the normal unit vector to $P$ directed inside
and $\nabla_{v_{(x,P)}}\,g$ is the covariant derivative
of the metric $g$ considered as a symmetric tensor on the cotangent
bundle.

Then the Assumption {\bf A2} is satisfied.

\end{prp}

{\it Proof.} In order to satisfy Assumption {\bf A2}
it suffices to find such $\varepsilon>0$
that for any $x\in B\setminus \overline{U}_1$ and any half-plane 
$P_x\subset T_xB,\,\, 0\in int(P_x)$ with
the distance $dist_{g_x}(0,\partial P_x)=\varepsilon$, and another half-plane
$P_x^{\prime}\subset T_x B$  parallel to $P_x$ such that $0\in \partial P_x^{\prime}$,
 one has
the following property: 
if $P_{x+\delta tv_{(x,P_x')}}$ is the half-plane
obtained from $P_x$ by a small covariant (with respect to the \emph{affine} 
connection $\nabla^{aff}$) shift $\delta t$
in the direction of 
$v_{(x,P_x^{\prime})}$, then 
$$dist_{g_{x+\delta tv_{(x,P_x^{\prime})}}}(0,P_{x+\delta tv_{(x,P_x^{\prime})}})\ge 
dist_{g_x}(0,\partial P_x)\,\,.$$ 
Here $g_x$ etc. denotes the  induced flat
metric on the tangent space $T_xB$. This property  guarantees that
the condition $dist_{g_x}(0,\partial P_{l,t})\ge \varepsilon$ will
propagate along the line. For a new line obtained as a result
of collision of $l_1$ and $l_2$ at the times $t_1$ and $t_2$ respectively
one has 
$$dist_{g_x}(0,\partial P_{l,0})\ge \min\{dist_{g_x}(0,\partial P_{l_1,t_1}),
dist_{g_x}(0,\partial P_{l_2,t_2})\}$$ since $\partial P_{l,0}$ contains
the intersection point $\partial P_{l_1,t_1}\cap \partial P_{l_2,t_2}$, see Figure 7.

\begin{figure}\label{figure7}
\centerline{\epsfbox{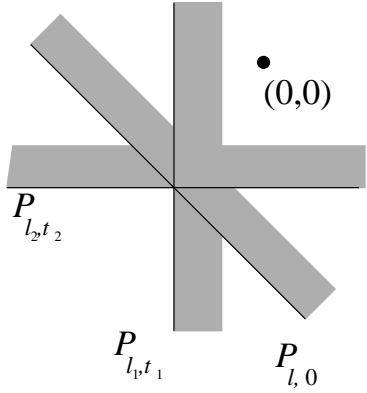}}
\caption{Three half-planes containing zero.}
\end{figure}

One can easily see that the infinitesimal inequality from above is equivalent
to
$$\delta t g_x(v_{(x,P_x^{\prime})},n_{P_x^{\prime}})+\varepsilon/2(g_{x+\delta tv_{(x,P_x^{\prime})}}-g_x)(n_{P_x^{\prime}},n_{P_x^{\prime}})\ge 0$$
(the change of the distance consists of two summands: one corresponds
to the shift along $\delta t v_{(x,P_x)}$ with the fixed metric, and
the other one corresponds to the change of the metric).
Taking the limit $\delta t\to 0$ we arrive to the inequality
for the covariant derivative of the metric with $C=2/\varepsilon$.
$\blacksquare$

Now our goal is to construct the field of directions $v$ and the metric $g$
satisfying the conditions 1--4 and the inequality from the last
Propostion. This will conclude the construction of the
set ${\cal L}$ of lines satisfying the Assumptions {\bf A1} and {\bf A2}.

Since $\partial U_1$ is a boundary of the convex set, we can locally model
it by the graph of function $y=f(x)$ such that $f^{\prime\prime}(x)>0$,
$f^{\prime}(x_0)=0$. We may assume that $P=P_0$ is the upper half-plane. Then we take 
$$v_{\left((x,y),P\right)}=\partial/\partial y+{(f(x)-f(x_0))/f^{\prime}(x)\over
{f(x)-f(x_0)+f(x)-y}}\,\,\partial/\partial x\,\,.$$

We extend this local model of $v$ near $\partial U_1$ to 
$B\setminus \overline{U}_1$ in such a way that conditions 1 and 2
are satisfied. It is clear that we can satisfy  conditions 3,4 as well
by taking a small perturbation of $v$. On Figure 8 there is a picture of the field
$(x,y)\mapsto v_{\left((x,y),P_0\right)}$.

\begin{figure}\label{figure8}
\centerline{\epsfbox{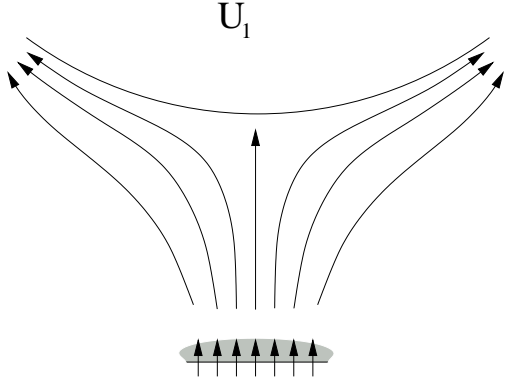}}
\caption{Vector field near $\partial U_1$ for $P=\mbox{ the upper half-plane}$.}
\end{figure}

For an arbitrary choice of the metric $g$ we have $g(n_P,v_{z,P})>0$
for all $(z,P)\in {\cal M}$. The problem with inequality

$$(\nabla_{v_{(z,P)}}\,g)(n_P,n_P)\le Cg(n_P,v_{(z,P)})$$
arises only as the point $z$ approaches $\partial U_1$. Indeed, 
in this case the vector $v_{(z,P)}$ can be very close to the tangent vector to 
$\partial P_z\subset T_zB$.

\begin{lmm} With the above choice of $v$ assume that the metric
satisfies for any $z\in \partial U_1$ the condition

$$(\nabla_{e_z}\,g)(n_z,n_z)=0\,\,,$$
where $e_z\in T_zB$ is the unit tangent vector to $\partial U_1$
and $n_z$ is the normal vector to $\partial U_1$ (all scalar
products and lengths are taken with respect to the metric $g$).

Then there exists $C>0$ such that

$$(\nabla_{v_{(z,P)}}\,g)(n_P,n_P)\le Cg(n_P,v_{(z,P)})$$
for all $(z,P)\in {\cal M}$.

\end{lmm}

{\it Proof.} We need to check that the ratio

$${(\nabla_{v_{(z,P)}}\,g)(n_P,n_P)}\over {g(n_P,v_{(z,P)})}$$
is bounded for $(z,P)\in {\cal M}$.

It suffices to prove the Lemma assuming that $U_1$ is the parabolic domain
$\{(x,y)\in {\R}^2| y>x^2\}$ and $P$ is the upper half-plane.
The vector field $v_{(z,P)}$ is given for $z=(x,y)$ by
the formulas

$$v_{(z,P)}=\partial/\partial y+{x\over{4x^2-2y}}\,\partial/\partial x\,\,.$$

The denominator is equal to
$g(n_P,v_{(z,P)})=\langle dy, v_{(z,P)}\rangle \cdot
\sqrt{g(\partial/\partial y, \partial/\partial y)}=
\sqrt{g(\partial/\partial y, \partial/\partial y)}=\exp(O(1))$
near $(0,0)$.

The numerator is equal to 

$${x\over{4x^2-2y}}f_1(x,y)+f_2(x,y)\,\,,$$
where $f_1(x,y)=(\nabla_{\partial/\partial x}\,g)(n_P,n_P)$
and $f_2(x,y)=(\nabla_{\partial/\partial y}\,g)(n_P,n_P)$
are two $C^{\infty}$-functions.

By assumption of the Lemma we have $f_1(0,0)=0$. Therefore
$|f_1(x,y)|\le const\, \max\{|x|, |w|\}$ where
$w=x^2-y$ is a convenient local coordinate near
the point $(0,0)$. Notice also that $f_2(x,y)=O(1)$.

Now we can estimate first summand of the numerator
assuming that $|x|$ and $|w|$
are sufficiently small. As we have seen, it is bounded by

$$I:={x\over{x^2+w}}O(\max\{|x|,|w|\}\,\,.$$

There are three cases which we need to consider.

a) If $0<w<x^2$ then $I={x\over{x^2}}O(|x|)=O(1)$.

b) if $x^2\le w<x$ then $I={x\over w}O(|x|)=O(1)$.

c) If $x\le w\le 1$ the $I={x\over w}O(|w|)=O(1)$.

We see that the numerator is bounded. This concludes
the proof of Lemma. $\blacksquare$

Finally, we have the following result.

\begin{lmm} There exists  metric $g$ satisfying the conditions
of  Lemma 6.

\end{lmm}
{\it Proof:} First of all, the condition on $g$ from Lemma 6 is the condition on a loop
 $g_{|T_z B}$ of scalar products on 2-dimensional spaces, here $z\in \partial U_1\simeq S^1$.
We can write $g=\exp(\psi)g_0$ where $\det(g_0)=1$
and $\psi$ is a smooth function. 
Then we have
$$\nabla_{e_z}(\exp{\psi}g_0)=
\exp({\psi})\nabla_{e_z}g_0+
\exp(\psi)\partial_{e_z} (\psi)\, g_0\,\,.$$
The equation
of Lemma 6 gives $\partial_{e_z} \psi=-
(\nabla_{e_z}g_0)(n_z,n_z)/g_0(n_z,n_z)$.
The RHS of this expression is known as long as we know $g_0$.
Hence we can say that
$d\psi=\beta_{g_0}$, where $\beta_{g_0}$ is a 1-form depending on 
the restriction $(g_0)_{|\partial U_1}$.
We see that it suffices to find such $g_0$ that
$\int_{\partial U_1}\beta_{g_0}=0$ (then $\psi$ and hence $g$ does exist).

Let us consider the functional $I(g_0)=\int_{S^1}\beta_{g_0}$. 
We can interpret a metric $g_0$ as a point in the Lobachevsky plane
${\cal H}=SL(2,\R)/SO(2)$. More precisely, let us consider
the  space $S$ of pairs $(g_0,P)$ where $g_0$ is a positive quadratic form   on $\R^2$ such
that $\det(g_0)=1$ and $P$ is a half-plane in $\R^2$ (the meaning of $P$ is the inward oriented 
tangent half-plane to $\partial U_1$ at point $z\in \partial U_1$). This 
space is naturally diffeomorphic to $S^{\ast}(\R^2)\times {\cal H}$. The latter
manifold can be identified in $SL(2,\R)$-equivariant way with the manifold consisting of pairs $(x,y)$,
where $x\in {\cal H}$ and $y$ belongs to the absolute. Hence $(g_0)_{|\partial U_1}$
is (locally) a non-parametrized path in $S$
(it would be a global path, if the bundle over $S^1$ given by the
all metrics on $S^1$ with the determinant $1$ was trivial).

Next we observe that the variation $\delta I(g_0)=\int_{N}\omega$,
where $N$ is a $2$-dimensional surface bounded by the paths
defined by $g_0$ and $g_0+\delta g_0$, and $\omega$ is a
canonical $SL(2,\R)$-invariant $2$-form on $S$. One can show that even by a small variation
of the path defined by $g_0$ we can make $I(g_0)$ an arbitrary
real number.
In particular, we can find $g_0$ such that $I(g_0)=0$.
This concludes the proof of Lemma 6. $\blacksquare$

Summarizing, we have constructed a set of lines satisfying
the Assumptions {\bf A1} and {\bf A2}. This concludes the proof
of Theorem 5. Thus we have obtained a solution
of the Lifting Problem, which is a $K$-analytic K3 surface.

\subsection{Independence and uniqueness}

It is natural to ask how the above construction of the $K$-analytic K3 surface
$(X^{an}, \Omega)$ depends on the choice of the set ${\cal L}$
of lines. We know that the ``periods" of $\Omega$ (they are encoded
in the initial $K$-affine structure) do not depend on ${\cal L}$
(see Sections 7.3, 10.4). In the light of Torelli theorem
(see Appendix B) it is natural to formulate the following
conjecture.

\begin{conj} The isomorphism class
of the pair $(X^{an}, \Omega)$ does not depend on the
choice of the set ${\cal L}$ of lines.

\end{conj}

More precisely, the change of ${\cal L}$
corresponds to the change of the projection
$\pi:=\pi_{\cal L}: X^{an}\to B$ (see Section 7.3).

\begin{rmk} For $B=S^2$ and $B^{sing}=\{x_1,\dots,x_{24}\}$
with the standard singular $\Z$-affine structure we have
constructed a $K$-analytic K3 surface depending on $20$ parameters
in $K^{\times}$. More precisely, we have a $20$-dimensional
$K$-analytic space of
conjugacy classes of representations

$$\pi_1(S^2\setminus B^{sing})\to SL(2, {\Z})\ltimes (K^{\times})^2$$
such that the monodromy around each singular point
is conjugate to the pair $(A,(1,1))$ where
$A\in SL(2,{\Z})$ is equal to
$$\left(\begin{array}{cc} 1 & 1\\0& 1
\end{array}\right)\,\,.$$
(compare with Section 3.3).

\end{rmk}

\subsection{Remark on the case of positive and mixed characteristic}

Our construction of $(X^{an}, \Omega)$
works even without the assumption $char\,k=0$ where $k$ is the
residue field of $K$.
This can be explained from the point of view of factorization
theorem (see Section 10.4). It turns out that symplectomorphisms
which appear in the infinite product in the RHS of the factorization
theorem are infinite series whose coefficients  are
integer polynomials in the coefficients of the ``parent"
symplectomorphisms.

For example, let $f_0(z)=1+\sum_{n\ge 1}c_n z^n$ and
$f_{\infty}(z)=1+\sum_{n\ge 1}d_n z^n$ be two  power series convergent when $|z|<1$. Let us consider
two symplectomorphisms:
$F_0(\xi,\eta)=(\xi, \eta f_0(\xi^{-1}))$ and 
$F_{\infty}(\xi,\eta)=(\xi f_{\infty}(\eta^{-1}), \eta)$
and decompose $F_{\infty}\circ F_0$ into the infinite
ordered product
$\prod_\to(F_{\lambda})$.
Here $$F_{p/q}(\xi, \eta)=(\xi f_{p/q}(\xi^{-p}\eta^{-q})^q, 
\eta f_{p/q}(\xi^{-p}\eta^{-q})^{-p})$$ where
$f_{p/q}(z)=1+\sum_{n\ge 1}c_n^{p/q}z^n$.
Then one can check that for any coprime $p,q\in {\bf Z}_{>0}$ and any $ n\ge 1$
one has $$c_n^{p/q}\in {\Z}[c_1,c_2,\dots ,d_1,d_2,\dots ]\,\,.$$

This implies that our construction works when one replaces
 $K$ by arbitrary commutative ring $R$ endowed with a complete
 non-trivial valuation $val:R\to (-\infty,+\infty]$.

\subsection{Further generalizations}

First of all, one can introduce a small parameter $\hbar\in K, |\hbar|<1$ of noncommutativity
 in the picture, coordinates $\xi, \eta$ will not commute but instead satisfy
the relation
$$\eta\xi=\xi\eta \exp(\hbar)\,\,.$$
For such a noncommutative analytic torus one can still define sheaf ${\cal O}^{can}_\hbar$ on $\R^2$ by
 the ``same'' formula as in the commutative case: 
$${\cal O}^{can}_\hbar(U)=\left\{ 
\sum_{n,m\in \Z}c_{n,m}\xi^n \eta ^m\,|\,\forall (x,y)\in U\,\,\,\sup_{n,m}\left(
\log |c_{n,m}| +nx+my\right)<\infty\right\}$$
where $U\subset \R^2$ is connected. Also one can construct a non-commutative deformation
 of the model sheaf near the singular point.
All arguments with the groups work as well. In this way we will obtain a kind of quantized 
K3 surface over a non-archimeden field.

Secondly, we believe that one can generalize our construction to higher dimensions.
 Instead of lines there will be codimension one walls which should be flat hypersurfaces
 with respect to $\Z$-affine structure and carry  foliations by parallel lines.
 Generically on the intersection of two such foliated hypersurfaces one can ``separate'' variables 
into the product of a purely 2-dimensional situation studied in the present paper, and $n-2$ dummy variables.
 Presumably everywhere except a countable union of codimension 2 subsets
 one can use 2-dimensional factorization and define gluing volume preserving maps.
 One can hope that by a kind of Hartogs principle the sheaf will have a canonical extension to the 
 whole space $B$.

\appendix

\section{Analytic geometry}

In this section we collect several facts and definitions 
about rigid analytic spaces and Clemens polytopes. 
Some of them are well-known,
the rest is borrowed from [KoT].

We always work over a complete non-archimedean local field $K$.
The field $K$ carries a valuation map 
$val_K:=val: K\to {\R}\cup \{+\infty\}$ such that $val(0)=+\infty,
val(1)=0, val(xy)=val(x)+val(y), val(x+y)\ge \min(val(x),val(y))$.

We will assume that the valuation is non-trivial.
The ring 
$$\O_K=val_K^{-1}({\R}_{\ge 0}\cup \{+\infty\})$$
is called the ring of integers of $K$. The residue field is defined as
$k=\O_K/m_K$, where $m_K=val_K^{-1}({\R}_{>0}\cup \{+\infty\})$
is the maximal ideal in $\O_K$.

Our main example is the field $K={\C}((t))$ of Laurent series
in one variable. In this case $val_K(\sum_{n\ge n_0}c_nt^n)=n_0$,
as long as $c_{n_0}\ne 0$.

\subsection{Berkovich spectrum}

We refer the reader to [Be1] for the general definition
of an analytic space and more details.
In this Appendix  we restrict ourselves to
analytic spaces associated with algebraic varieties (although
we use the general definition in the paper as well).

Let $R=R/K$ be a commutative unital finitely generated $K$-algebra. 
The underlying set of the 
Berkovich spectrum $Spec^{an}(R):=Spec^{an}(R/K)$ 
can be defined in two  ways. First one uses valuations
(or, equivalently, multiplicative seminorms).

\begin{defn} (Valuations)
A point $x$ of $Spec^{an}(R/K)$ is an additive valuation 
$$
val_x: R\to \R\cup\{+\infty\}
$$
extending $val:=val_K$, i.e. it is a map satisfying the conditions
\begin{itemize}
\item $val_x(r+r')\le \max(val_x(r),val_x(r'))$;
\item $val_x(rr')=val_x(r)+val_x(r')$; 
\item $val_x(\lambda)=val_K(\lambda)$
\end{itemize}
for all $r,r'\in R$ and all $\lambda\in K$.
\end{defn}

Having a valuation and a real number $q_0\in (0,1)$ 
one can
define the {\it multiplicative seminorm} $|a|=q_0^{val_K(a)}, a\in R$. 
In particular,
in the previous definition 
one can take seminorms $|\cdot|_x$ instead of valuations
$val_x(\cdot)$. The reader has noticed that in the main
body of the paper, for $R=K$
we often took $|a|=e^{-val(a)}$. It is easy to translate
the definition of Berkovich spectrum to the language of multiplicative
seminorms. We use it freely in the paper. 

The second way to define $X^{an}$ uses evaluations
(characters).

\begin{defn} 
(Evaluation maps)
A point  $x$   of $Spec^{an}(R/K)$ is an equivalence class of
homomorphisms of $K$-algebras 
$$
eval_x\,:\, R\to K_x\,\,,
$$
where $K_x\supset K$ is a
complete field equipped with a non-archimedean valuation, which
extends the valuation $val_K$, and such that  
$K_x$ is generated by the closure of the image of $eval_x$.
\end{defn}
The field $K_x$ is determined by $x\in X^{an}$ in a canonical way. We define for $r\in R$ and $x\in X^{an}$
 the ``value'' $r(x)\in  K_x$ as the image $eval_x(r)$.
 
In order to pass from the first description of 
$Spec^{an}(R/K)$
to the second, starting with a valuation
$val_x$ one defines the field $K_x$
as the completion of the field of fractions of $R/I_x$, 
where $I_x=(val_x)^{-1}(\{+\infty\})$.

\begin{defn} 
The topology on $Spec^{an}(R/K)$
is the weakest topology such that for all $r\in R$ the map
$$
\begin{array}{ccc}
Spec^{an}(R/K) & \to   & \R\cup\{+\infty\},\\
        x             & \mapsto & val_x(r)
\end{array}
$$
is continuous. 
\end{defn}

An element $f\in R$ defines a function $f: Spec^{an}(R)\to K_x$,
where $K_x$ is the non-archimedean valuation field, which is
the completion of the
field of fractions of the domain $R/\ker(val_x)$. 
Since each $K_x$ carries a seminorm, we obtain a function
$|f|: Spec^{an}\to {\R}_{\ge 0}, x\mapsto |f(x)|$.

A fundamental system of neighborhoods $U=U_x\subset Spec^{an}(R)$ 
of a point $x$ is parametrized by the following data: 
a finite collections of functions
$$
(f_i)_{i\in I},\,\, (g_j)_{j\in J}\,\,\,\in R
$$
and numbers 
$$
\beta^+_i,\beta^-_i,\gamma_j\in \R_{>0} 
$$
such that  $\beta_i^- <|{f_i(x)}|< \beta_i^+,\,\, |{g_j(x)}|=0$,
 The corresponding neighborhood consists of points $x'$ such that
 $\beta_i^- <|{f_i(x')}|< \beta_i^+,\,\, 
|{g_j(x')}|<\gamma_j$
for all $i\in I, j\in J$ and $x'\in U$.

Let us assume that elements $(f_i)_{i\in I},\, (g_j)_{j\in J}$ generate $R$, i.e.
  $$R=K[(f_i)_{i\in I},\, (g_j)_{j\in J}]/I$$ where $I$ is an ideal.
Let us consider the algebra of series
$$
s=\sum_{n_I\in \Z^I,\,m_J\in \N^J} c_{I,J}f_I^{n_I}g_{J}^{m_J}
$$
with constants $c_{I,J}\in K$,
absolutely convergent when variables $(f_i)_{i\in I},\, (g_j)_{j\in J}$ satisfy the above inequalities.
The quotient of this algebra by the topological closure of ideal $I$ is the algebra ${\cal O}_{Spec^{an}(R/K)}(U)$.

As in the case of schemes we can glue $Spec^{an}(R/K)$
into ringed spaces called {\it analytic spaces} (or rigid
analytic spaces).
Moreover we get a functor
$$
\begin{array}{ccc}
(Schemes/K) &  \to    &  
(K-{\rm analytic}\,\, {\rm spaces}) \\
     X   & \mapsto &  (X^{an},{\cal O}_{X^{an}}).
\end{array}
$$

\begin{prp}
The space $X^{an}$

a) is a locally compact Hausdorff space as long as $X$ is separated;

b) has the homotopy type of a finite $CW$-complex;

c) is contractible if $X$ has good reduction with
irreducible special fiber.

\end{prp}

\begin{exa}
Let $X={\bf A}^1=Spec(K[x])$ be the affine line.
The analytic space $X^{an}$ contains, among others, 
points of the following types:
\begin{itemize}
\item $X(K)\hookrightarrow   X(\overline{K})/Gal(\overline{K}/K)\hookrightarrow X^{an}$;
\item for  $r\in {\R}_{\ge 0}$ define
$$
 |\sum_{j=0}^d c_jz^j|_r := \max_j(|{c_j}|r^j)\,\,.
$$
This gives an embedding $ \R_{\ge 0}\hookrightarrow X^{an}$.
\end{itemize}
We see that $X^{an}$ contains, in a sense, both $p$-adic and real
points.

\end{exa}

Define the \emph{cone} over $X^{an}$ as 
$$
C_{X^{an}}(\R):=X^{an}\times \R_{>0}\,\,.
$$
We interpret a point 
${\bf x}=(x,\lambda)$ of $C_{X^{an}}(\R)$ as a 
$K_x$-point of
$X$,  where $K_x\supset K$ is a  
complete field with the $\R$-valued valuation 
$$
val_{\bf x}:=\lambda\,val_x\,\,,
$$ 
whose restriction to $K$ is proportional to $val_K$.
The set of points ${\bf x}\in C_{X^{an}}(\R)$ such that
the valuation $val_{\bf x}$ is $\Z$-valued is denoted by 
$C_{X^{an}}(\Z)$.

\subsection{Algebraic torus and the logarithmic map}

Here we will describe explicitly the main example for our paper.
 Let  $X={\bf G}_m^n=Spec(K[z_i^{\pm 1}]), 1\le i\le n$ be an algebraic torus.
and $X^{an}=({\bf G}_m^{an})^n$ the corresponding analytic space.

Firstly, we define an embedding $i_{can}:\R^n\mono X^{an}$. 
 For real vector $(x_i)_{1\le i\le n}\in \R^b$
 the corresponding point
 $p:=i_{can}(x_1,\dots,x_n)\in X^{an}$
 will be described in terms of valuations.

 For every Laurent polynomial $f=\sum_{I\in {\Z}^n}c_I z^I,\,\, c_I\in K$ we set 
$$val_p(f):=\min_{I\in \Z^n}\left(val(c_I)-\sum_{i=1}^n x_i I_i\right)\,\,.$$

Secondly, we define a projection $\pi_{can}: X^{an}\to \R^n$ by formula
$$\pi_{can}(y)=\left(-val_y(z_1),\dots,-val_y(z_n)\right)=\left(\log|z_1|_y,\dots,\log|z_n|_y\right)\,\,.$$
 The fiber over
a point $(x_1,\dots ,x_n)\in {\R}^n$ can be identified with
the set of
such seminorms $|\cdot|_y$ that $|z_i|_y=\exp(x_i), 1\le i\le n$.
We see $\pi_{can}$ is a kind of torus fibration\footnote{This is the origin of the term
 ``analytic torus fibration'' introduced in Section 4.1.}. Moreover, $\pi_{can}\circ i_{can}=id_{\,\R^n}$.

For any open connected $U\in ({\R})^n$
the $K$-algebra of analytic functions on $\pi_{can}^{-1}(U)$
consists of series
$f=\sum_{I\in {\Z}^n}c_I z^I$ with coefficients $c_I\in K$
such that for any $p=(x_1,\dots,x_n)\in U$ 
we have
$\log |c_I| +\sum_{i=1}^n x_i I_i\to +\infty$ when
$|l|\to +\infty$. It is easy to see that $\pi_{can}^{-1}(U)=\pi_{can}^{-1}(Conv(U))$
where $Conv(U)$ is the convex hull of $U$.

The sheaf 
$(\pi_{can})_{\ast}({\cal O}_{X^{an}}):={\cal O}^{can}_{{\R}^n}$
(canonical sheaf)
plays an important role in the paper 
(see Sections 4.1,  7.3, 8).

\subsection{Clemens polytopes}

Let $X$ be a smooth proper scheme over the non-archimedean field $K$.
We assume that $K$ carries a discrete valuation $val$ such that
$val(K^{\times})={\Z}$.

\begin{defn}
A {model} of $X$ is a  scheme of finite type  
${\cal X}/{\O}_K$ flat and proper over
${\O}_K$, together with an isomorphism
${\cal X}\times_{Spec(K)} Spec({\O}_K)\simeq X$.
Denote the special fiber of ${\cal X}$ by 
$$
{\cal X}^0:={\cal X}\times _{Spec(K)} Spec(k)\,\,.
$$ 
\end{defn}

A model has no nontrivial automorphisms. 
Thus, the stack of equivalence 
classes of models is in fact a set,  
which we denote by $Mod_X$. 
It carries a natural partial order. Namely, we say that
${\cal X}_1 \ge {\cal X}_2$
if there exists a map ${\cal X}_1\to {\cal X}_2$ over 
$Spec({\O}_K)$. 
Such a map is automatically unique. 

\begin{defn} A model 
${\cal X}$ has \emph{ normal crossings}
if the scheme ${\cal X}$ is regular and    
the reduced  subscheme
${\cal X}^0_{red}$ is a divisor with normal
crossings.
\end{defn}
By the resolution of singularities, in the case $char\,k=0$
we know that every model is 
dominated by a model with normal crossings.
\begin{defn} A model ${\cal X}$    
has  {simple  normal crossings}
({snc model} for short)     
if
\begin{itemize}
\item 
it has normal crossings;
\item  
all irreducible components of   
${\cal X}^0_{red}$ are smooth and
\item 
all intersections of irreducible components of  ${\cal X}^0_{red}$
 are either empty or irreducible.
\end{itemize}
\end{defn}
The set of equivalence classes of 
snc models will be denoted by  $Mod_X^{snc}$.
It is a filtered partially ordered set.
The order is given by dominating maps of models
which give the identity automorphism on the generic fiber.

It is easy to show that starting with 
any model with normal crossings
and applying blow-ups centered at 
certain self-intersection loci
of the special fiber we can get a snc model. 
In what follows we use snc models
only. This choice is dictated by 
convenience and not by necessity.
Working with 
snc models has the advantage that 
all definitions and calculations
can be made very transparent. The reader can consult [Be2]
 for the approach in the general case, without the use of the resolution of singularities.

Let ${\cal X}$ be  an  snc model and  
$I=I_{{\cal X}}$ the  set of irreducible 
components of ${\cal X}^0_{red}$.
Denote by $D_i\subset {\cal X}$ the divisor 
corresponding to $i\in I$.       
For any finite non-empty subset $J\subset I$ put
$$
D_J:=\bigcap_{j\in J} D_j\,\,.
$$ 
By the snc property the set $D_J$ is either empty 
or is a smooth connected proper variety over $k$ of dimension 
$\dim(D_J)=(n-|J| +1)$.
For a divisor $D_i\subset {\cal X}^0$ 
we denote by $d_i\in \Z_{>0}$ the 
order of vanishing of $u$  at $D_i$, 
where $u\in K$ is an uniformizing element, $val_K(u)=1$.
Equivalently,   
$d_i$ is  the multiplicity of $D_i$ in ${\cal X}^0$.

\begin{defn}    
The {Clemens polytope} $S_{\cal X}$ is 
the finite simplicial subcomplex
of the simplex $\Delta^I$ such that $\Delta^J$ is a  face
of $S_{\cal X}$ iff $D_J\ne\emptyset$.
\end{defn}

Clearly, $S_{\cal X}$ is a nonempty connected CW-complex.
We will also consider the cone over $S_{\cal X}$:
$$
C_{\cal X}(\R):=\left\{\sum_{i\in I} a_i \langle D_i\rangle |\,
 a_i\in {\R}_{\ge 0},\,\,\, \bigcap_{i:\,a_i>0}D_i\ne\emptyset\right\}
\setminus\{0\}\subset \R^I\,\,.
$$ 
Analogously, we can define    
$C_{\cal X}(\Z),\,$.

We identify $S_{\cal X}$ with the following subset of 
$C_{\cal X}(\R)$:
$$
\left\{ \sum_{i\in I} a_i \langle D_i\rangle\in C_{\cal X}(\R)| \,\sum_i a_i 
d_i=1\right\}\,\,.
$$
Obviously, we can also describe 
$S_{\cal X}$ as a quotient of $C_{\cal X}(\R)$:
$$
S_{\cal X}=C_{\cal X}(\R)/\R^{\times}_+\,\,.
$$

\subsection{Simple blow-ups}

Let ${\cal X}$ be an snc model, 
$J\subset  I_{\cal X}$ a non-empty subset and   
$Y\subset D_J$ a smooth  irreducible variety 
of dimension less or equal than $n$. Let us assume that
$Y$ intersects transversally (in $D_J$) 
all subvarieties $D_{J'}$ of $D_J$
(for $J'\supset J$),
and that all intersections $Y\cap D_J$ are
either empty or irreducible.   
It is obvious that the blow-up ${\cal X}':=Bl_Y({\cal X})$ of    
${\cal X}$ with the center at $Y$ is again a snc model.

\begin{defn}  
For a pair of snc models ${\cal X}'\ge   {\cal X}$   as above
we say that ${\cal X}'$ is obtained from 
${\cal X}$ by a simple blow-up.
If $Y=D_J$ we say that we have  a  
simple blow-up of the first type.
Otherwise (when $\dim(Y)<\dim(D_J)$),
we have a simple blow-up of the second type.
\end{defn}

Let us describe the behavior  of $S_{\cal X}$ under simple blow-ups.
To the set of vertices we add a 
new vertex corresponding to the divisor 
$\widetilde{Y}$ obtained from $Y$:
$$
I_{{\cal X}'}=I_{\cal X}\sqcup \{new\},\,\,\, D_{new}:=\widetilde{Y}\,\,.
$$
The degree of the new divisor is      
(for both the first and the second type)
$$
d_{new}:=\sum_{i\in J} d_j\,\,.
$$

For blow-ups of the first type 
we have automatically $\# J>1$.
Here is the list of faces of $S_{\cal {X}'}$:

1) $I'$ for $I'\in Faces(S_{\cal X}),\,I'\not\subset J$; 

2) $I'\sqcup \{new\}$ for $I'\in 
Faces(S_{\cal X}),\,I'\ne J,\,I'\cup J\in Faces(S_{\cal X})$;

3) the vertex $\{new\}$. 
\newline For blow-ups of the second type the list of faces 
of $S_{{\cal X}'}$ is

1) $I'$ for $I'\in Faces(S_{\cal X})$;

2) $I'\sqcup \{new\}$ for $I'\in
Faces(S_{\cal X}),\,I'\supset J,\,Y\cap D_{I'}\ne \emptyset$;

3) the vertex $\{new\}$. 

On can deduce from results [AKMW] the following 
\begin{thm}(Weak factorization) Assume that $char\, k=0$. 
Then for any two snc models ${\cal X},\,{\cal X}'$ 
there   exists a finite 
alternating sequence
of simple blow-ups
$$
{\cal X}<{\cal X}_1>{\cal X}_2<\dots <{\cal X}_{2m+1}>{\cal X}'\,\,.
$$
\end{thm}

\begin{cor}               
Simple homotopy type of 
$S_{\cal X}$ does not depend on the choice of a snc
model ${\cal X}$.
\end{cor}

\subsection{Clemens cones and valuations}

Let ${\cal X}$ be a snc model of $X$.
We define a map
$$
i_{\cal X}:C_{\cal X}(\R)\to C_{X^{an}}(\R)
$$
such as follows.
For $J=\{j_1,\dots,j_k\}\subset I_{\cal X}$ such that $D_J\ne \emptyset$ let us
consider a point $x\in C_{\cal X}(\R)$
$$
x=\sum_{i=1}^k a_i \langle D_{j_i}\rangle,\,\,\,a_i\in {\R}_{>0}\,\,
\forall i\in\{1,\dots,k\}
$$
and an affine Zariski open subset
$U\subset {\cal X}$ containing the generic point of $D_J$.
 One can embed ${\cal O}(U)$ into the algebra
 of formal series $K_J[[z_1,\dots,z_k]]$ where $K_J$ is the field of rational functions on $D_J$
 and $z_i=0$ are equations of divisors $D_{j_i},\,\,i=1,\dots,k\,$.
We define a valuation $v_x$ of ${\cal O}(U)$ by the formula
$$
v_x\left(\sum_{n_1,\dots,n_k\ge 0} c_{n_1,\dots,n_k}\prod_{i=1}^k z_i^{n_i}\right)=
\inf\left\{\sum a_i n_i\,|\,c_{n_1,\dots,n_k}\ne 0\right\}\,\,.
$$

We define $i_{\cal X}(x)$ 
to be the image of the point $v_x\in 
Spec^{an}({\cal O}(U)/K)$ in $X^{an}$. It is easy to check 
that the element 
$i_{\cal X}(x)$ does not depend on the choice of the open subset $U$.

The following proposition is obvious:
\begin{prp}
The map $i_{\cal X}^\R$ is an embedding.
\end{prp}

We will denote also by $i_{\cal X}$ the induced embedding $S_{\cal X}\mono X^{an}$.

\subsection{Clemens cones and paths}

For a model ${\cal X}$ we can interpret elements of 
$C_{X^{an}}(\Z)$
as \emph{paths} in ${\cal X}$, i.e. equivalence
classes of maps 
$$
\phi: Spec({\O}_L)\to {\cal X},
$$ 
where ${\O}_L$ is the 
ring of integers in a field $L$ 
with discrete valuation in $\Z$, such that the image of 
$\phi$ does not lie in ${\cal X}$.
We define the map
$$
p_{\cal X}^\Z\,:\, C_{X^{an}}(\Z)\to C_{\cal X}(\Z)
$$
as 
$$
p_{\cal X}^\Z([\phi]):=\sum_i a_i\langle D_i\rangle,
$$
where $a_i\in\Z_{\ge 0}$ is the multiplicity
of the intersection of the path $\phi$ with the divisor
$D_i,\,\,i\in I_{\cal X}$.

The following proposition can be derived from [Be1].

\begin{prp} 
The map 
$p_{\cal X}^\Z$ extends uniquely to a 
continuous $\R_+^{\times}$-equiva\-riant
map $p_{\cal X}^\R\,:\, C_{X^{an}}(\R)\to C_{\cal X}(\R)$. 
The map $p_{\cal X}^\R$ is a surjection. 
\end{prp}

We denote by $p_{\cal X}:X^{an}\to  S_{\cal X}$ 
the map induced by $p_{\cal X}^\R$. 

Let $f: {\cal X}^{\prime}\to {\cal X}$ be a dominating
map of models. Let us denote by $m_{i,i^{\prime}}\in {\Z}_{\ge 0}$
the multiplicity of a divisor 
$D_{i^{\prime}}, i^{\prime}\in I_{{\cal X}^{\prime}}$ in the
proper pull-back of $D_i, i\in I_{\cal X}$. 
We define 
$p_{{\cal X}^{\prime}, {\cal X}}^{\Z}: C_{X^{an}}(\Z)\to C_{\cal X}(\Z)$
by the formulas $\sum_{i^{\prime}}a_{i^{\prime}}\langle D_{i^{\prime}}\rangle \mapsto
\sum_i m_{i,i^{\prime}}a_{i^{\prime}} \langle D_{i} \rangle$. Let  
$p_{{\cal X}^{\prime}, {\cal X}}: S_{X^{an}}(\R)\to S_{\cal X}(\R)$
 be the corresponding by map of Clemens polytopes.

Then we have the following result, which is easy to prove.

\begin{lmm}
For any dominating map of models ${\cal X}'\to {\cal X}$ 
we have
$$
p_{{\cal X}}^\Z=p_{{\cal X}',{\cal X}}^\Z\circ p_{{\cal X}'}^{\Z}\,\,.
$$
\end{lmm}

\begin{cor}
For dominating maps ${\cal X}''\ge {\cal X}'\ge {\cal X}$
we have
$$
p_{{\cal X}'',{\cal X}}=
p_{{\cal X}',{\cal X}}\circ p_{{\cal X}'',{\cal X}'}\,\,.
$$
\end{cor}

\begin{thm}
For any algebraic $X$ the analytic space $X^{an}$ 
is a projective
limit over the partially ordered  
set of snc models ${\cal X}$ of
Clemens polytopes $S_{\cal X}$. The connecting maps are 
$p_{{\cal X}',{\cal X}} $.
\end{thm}

With any meromorphic at $t=0$ family of smooth complex projective varieties $X_t,\,\,\,0<|t|<\epsilon$ 
one can associate a
 variety $X$ over the field $\C((t))$. It is easy to see that for any snc model $\cal X$ one can 
canonically complete the family $X_t$ by adding $S_{\cal X}$ as the fiber over $t=0$.
 The total space is not a complex manifold by just a Hausdorff locally compact space
 which maps properly to the dick $\{t\in \C\,|\,|t|<\epsilon\}$. Passing to the projective limit
we see that one can compactify the family $X_t$ at $t=0$ by $X^{an}$.

\section{Torelli theorem for K3 surfaces}

Here we recall the  classification theory of complex K3 sufaces 
(see [PSS] and its extension to non-algebraic case in [LP]). 
Let $X$ be a complex K3 surface, i.e. smooth connected
complex manifold with $\dim_{\C}X=2$ which admits a nowhere
vanishing holomorphic $2$-form $\Omega$, and such that
$H^1(X,{\Z})=0$.

It is known that the group $H^2(X,{\Z})$ endowed with the Poincare
pairing $(\cdot,\cdot)$ is isomorphic to the lattice
$$\Lambda_{K3}=\left(\begin{array}{cc} 0 & 1 \\ 1 & 0\end{array}\right)
\oplus
\left(\begin{array}{cc} 0 & 1 \\ 1 & 0\end{array}\right)
\oplus
\left(\begin{array}{cc} 0 & 1 \\ 1 & 0\end{array}\right)
\oplus \left(-E_8\right)\oplus \left(-E_8\right)$$
of signature $(3,19)$.

Complex $1$-dimensional vector space 
$H^{2,0}(X)={\C}\cdot [\Omega]\subset H^2(X,{\Z})\otimes {\C}$
satisfies the condition $(v,v)=0, (v,\overline{v})>0$ for
any non-zero vector $v$.
Finally, it is known that $X$ admits a K\"ahler metric, and
K\"ahler cone ${\cal K}_X\subset H^2(X,{\R})$ of all 
K\"ahler metrics on $X$ is an open subset of
$C_X:=\{[\omega]\in H^2(X,{\Z})\otimes {\R}|([\omega],[\Omega])=0,
([\omega],[\omega])>0 \}$.
In fact ${\cal K}_X$ is a connected component of the set
$C_X\setminus \cup_{v\in H^2(X,{\Z}), (v,v)=-2, (v,[\Omega])=0}H_v$,
where $H_v$ is the hyperplane orthogonal to $v$.

Axiomatizing these data we arrive to the following definition.

\begin{defn} K3 period data is a quadruple
$(\Lambda, (\cdot,\cdot), H^{2,0}, {\cal K})$ consisting
of a free abelian group $\Lambda$, a symmetric pairing
$(\cdot,\cdot): \Lambda\times \Lambda\to {\Z}$, a
$1$-dimensional complex vector subspace  $H^{2,0}\subset\Lambda\otimes {\C}$
and a set ${\cal K}\subset \Lambda\otimes {\R}$ satisfying the
following conditions:
\begin{enumerate}
\item $rk\,\Lambda =22$;

\item $(\Lambda, (\cdot,\cdot))$ is isomorphic to $\Lambda_{K3}$;

\item for any $v\in H^{2,0}\setminus \{0\}$ one has $(v,v)=0$ and 
$(v,\overline{v})>0$;

\item the set ${\cal K}$ is a connected component of
$C\setminus \cup_{v\in \Lambda, (v,v)=-2, (v,H^{2,0})=0}H_v\,$,
where $C=\{w\in \Lambda\otimes {\R}|(w,H^{2,0})=0, (w,w)>0\}$ and
$H_v$ is the hyperplane orthogonal to $v$.
\end{enumerate}

\end{defn}

The K3 period data form a groupoid. On the other hand, K3 surfaces 
 also form a groupoid (morphisms are isomorphisms of K3 surfaces).
Then classical global Torelli theorem can be formulated in
the following way.

\begin{thm} Groupoid of K3 surfaces is equivalent to the
groupoid of K3 period data.

\end{thm}

In particular the automorphism group of a K3 surface is isomorphic
to the automorphism group of its period data.

More generally one can speak about holomorphic
families of K3 surfaces over complex analytic spaces.
For a K3 surface over an analytic space $M$ the period
data consist of a local system of integral
lattices $(\Lambda, (\cdot,\cdot))$ pointwise isomorphic
to $\Lambda_{K3}$, a holomorphic line 
subbundle $H^{2,0}$ of $\Lambda\otimes_{\Z}{\cal O}_M$ 
which is isotropic
with respect to the symmetric pairing $(\cdot,\cdot)$,
and satisfies pointwise the condition $(v,\overline{v})>0,
v\in H^{2,0}_x\setminus \{0\}, x\in M^{red}$, and an open subset
of the total space of the bundle over $M^{red}$ with
the fibers $\Lambda_x\otimes {\R}\cap (H^{2,0})^{\perp}$
($(H^{2,0})^{\perp}$ is the orthogonal complement) satisfying 
pointwise the condition 4) from the definition of K3 period data.
Then Torelli theorem holds for families as well.

\vspace{5mm}

{\bf References}

\vspace{3mm}

[AKMW], D.~Abramovich, K.~Karu, K.~Matsuki, J.~Wlodarczyk,
 {\it Torification and factorization of birational maps},  
J. Amer. Math. Soc.  {\bf 15}  (2002),  no. 3, 531--572, and preprint 
math.AG/9904135.

\vspace{2mm}

[Ar] V.~Arnold, {\it Mathematical methods of classical mechanics},
Springer, 1997.

\vspace{2mm}

[Au] M.~Audin, {\it Spinning tops. A course on integrable systems}, 
Cambridge Studies Adv. Math., vol. 51, Cambridge University Press, 1996.

\vspace{2mm}

[Be1] V.~Berkovich, {\it 
Spectral theory and analytic geometry over
non-archi\-medean fields}, AMS Mathematical Surveys and Monographs,
n. 33, 1990.

\vspace{2mm}

[Be2] V.~Berkovich, {\it Smooth p-adic analytic spaces are locally contractible},
 Inv.Math., {\bf 137}, 1--84 (1999).

\vspace{2mm}

[Be3] V.~Berkovich, {\it Smooth p-adic analytic spaces are locally contractible, II},
 to appear.



\vspace{2mm}

[GS] M.~Gross, B.~Siebert, {\it Mirror Symmetry via Logarithmic degeneration data, I},
 preprint math.AG/0309070.

\vspace{2mm}

[GW] M.~Gross, P.~M.~H.~Wilson, {\it Large complex structure limits
of K3 surfaces},  J. Differential Geom.  {\bf 55}  (2000),  no. 3, 475--546, and
 preprint math.DG/0008018.

\vspace{2mm}

[HZh] C.~Haase, I.~Zharkov, 
    {\it Integral affine structures on spheres and torus 
fibrations of Calabi-Yau toric hypersurfaces I}, preprint 
math.AG/0205321.

\vspace{2mm}
[KN] S.~Kobayashi, K.~Nomizu, {\it Foundations of differential
geometry}, vol. 1, John Wiley and Sons, 1963.

\vspace{2mm}

[Ko] M.~Kontsevich, {\it Homological algebra of mirror symmetry},
Proc. ICM Z\"urich, vol.1, 1994, and preprint math.AG/9411018.

\vspace{2mm}

[KoSo] M.~Kontsevich, Y.~Soibelman, {\it Homological mirror symmetry
and torus fibrations},  in  {\it Symplectic geometry and mirror symmetry}
 (Seoul, 2000),   World Sci. Publishing, River Edge, NJ, 2001, 203--263,
 and  math.SG/0011041.

\vspace{2mm}

[KoT] M.~Kontsevich, Yu.~Tschinkel, {\it Non-archimedean K\"ahler
geometry}, in preparation.

\vspace{2mm}

[LeS] N.~C.~Leung, M.~Symington, 
{\it Almost toric symplectic four-manifolds}, preprint
math.SG/0312165.

\vspace{2mm}

[LYZ] J.~Loftin, S.-T.~Yau, R.~Zaslow, {\it Affine manifolds, 
SYZ geometry and the ``Y'' vertex}, preprint
 math.DG/0405061.

\vspace{2mm}
[LP] E.~Looijenga, C.~Peters, {\it Torelli theorems for K\"ahler K3 surfaces}, Comp. Math. {\bf 42} (1981),  
 145--186.
\vspace{2mm}

[Mi] G.~Mikhalkin, {\it Amoebas of algebraic varieties and tropical
geometry}, preprint math.AG/0403015.

\vspace{2mm}

[Mor] D.~R.~Morrison, {\it Mathematical aspects of mirror symmetry}, in 
 {\it Complex algebraic geometry} (Park City, UT, 1993),  
IAS/Park City Math. Ser., 3, Amer. Math. Soc., Providence, RI, 1997, 265--327,
 and preprint alg-geom/9609021.
\vspace{2mm}

[PSS] I.~I.~Pjatecki\u{i}-\u{S}apiro, I.~R.~\u{S}afarevi\u{c}, {\it A Torelli theorem for algebraic surfaces
 of type K3}, Math. USSR Izvestija {\bf 5} (1971), No. 3, 547--588.

\vspace{2mm}

[SYZ] A.~Strominger, S.-T.~Yau, E.~Zaslow, {\it Mirror symmetry
is T-duality}, Nucl. Phys. B479 (1996), 243-259.

\vspace{2mm}

[Tyu] A.~Tyurin, {\it On Bohr-Sommerfeld bases},   Izv. Math.  {\bf 64}  (2000), 
 no. 5, 1033--1064,
 and  preprint math.AG/9909084. 

[Zu] N.~Zung, {\it Symplectic topology of integrable Hamiltonian
systems, II: Topological classification} , Compositio Math.,
138:{\bf 2} (2003), 125-156.

\vspace{5mm}

Addresses:

M.K.: IHES, 35 route de Chartres, F-91440, France

{maxim@ihes.fr}

\vspace{5mm}

Y.S.: Department of Mathematics, KSU, Manhattan, KS 66506, USA

{soibel@math.ksu.edu}

\end{document}